\documentclass[12pt,letterpaper]{amsart}
\usepackage{fullpage}
\usepackage{amsthm}
\usepackage{amsmath}
\usepackage{setspace}
\usepackage{amsfonts}
\usepackage{amssymb}
\newtheorem{theorem}{Theorem}[section]
\newtheorem*{proposition*}{Proposition}
\newtheorem*{theorem*}{Theorem}

\newtheorem{lemma}[theorem]{Lemma}
\newtheorem{proposition}[theorem]{Proposition}
\newtheorem{corollary}[theorem]{Corollary}
\newtheorem{conjecture}[theorem]{Conjecture}
\theoremstyle{remark}
\newtheorem*{remark}{Remark}
\usepackage{tikz}
\usetikzlibrary{matrix,arrows,decorations.pathmorphing}
\usepackage{tikz-cd}
\usepackage{comment}
\usepackage{subfiles}
\usepackage{enumitem}
\usepackage{graphicx}
\usepackage{import}
\usepackage{tocbasic}
\usepackage[colorlinks=true,linkcolor=black,anchorcolor=black,citecolor=black,filecolor=black,menucolor=black,runcolor=black,urlcolor=black]{hyperref}
\hypersetup{linktocpage}

\DeclareMathOperator{\Exp}{Exp}
\usepackage{bbding}
\newcommand{\flower}{\text{\FourClowerOpen}}
\title{An Optimal Yoccoz Inequality for near-parabolic quadratic polynomials}
\author{Alex Kapiamba}

\begin{document}
	\maketitle
	
	Abstract:
	The Pommerenke-Levin-Yoccoz inequality gives a $O(1/q)$ bound on the size of the $p/q$-limbs of the Mandelbrot set. For over thirty years it has been conjectured that this bound can be improved to $O(1/q^2)$. By studying the Hausdorff limits of external rays when a periodic cylce becomes parabolic, we provide the first examples of limbs satisfying the $O(1/q^2)$ bound. We also show that $O(1/q^2)$  is optimal for these examples.

	\section{Introduction}

	We consider the quadratic polynomials
	$$Q_\alpha(z):= e^{2\pi i\alpha} z+z^2$$
	with $\alpha\in \mathbb{C}$, and let $\mathcal{M}$ be the set of all $\alpha$ so that the orbit  under $Q_\alpha$  of the critical point $z = -e^{2\pi i\alpha}/2$stays bounded. Understanding the geometry of the Mandelbrot set, or equivalently $\mathcal{M}$, has been a central pillar of holomorphic dynamics since the seminal work of Douady and Hubbard in \cite{Orsay1} and \cite{Orsay2}.
	While much of the structure of $\mathcal{M}$ is now understood there is one fundamental problem  which remains unanswered: is $\mathcal{M}$ locally connected (MLC)? 
	
	One of the first major breakthroughs towards MLC was made by Yoccoz, who reduced the problem to understanding the self-similarity of the Mandelbrot set. A point $\alpha\in \mathcal{M}$ is called \textit{parabolic} if $Q_\alpha$ has a parabolic periodic cycle.
	For any parabolic parameter $\alpha$ there is a unique bounded component $\mathcal{L}_\alpha$ of $\mathcal{M}\setminus\{\alpha\}$, called the $\alpha$-\textit{limb} of $\mathcal{M}$.
	Inside the limb there is a \textit{small Mandelbrot copy} $\mathcal{M}_\alpha\subset \mathcal{L}_\alpha$ which admits a dynamically defined homeomorphism $\tau_\alpha: \overline{\mathcal{M}_0}\to \overline{\mathcal{M}_\alpha}$. 	
	For any set $X\subset \mathbb{C}$, let $\text{Diam} \,X$ denote the Euclidean diameter of $X$.
	Yoccoz showed that to prove MLC it suffices to control the sizes of nested small Mandelbrot copies.
	\begin{theorem}[Yocooz \cite{hubbard}]\label{yoccoz infinitely renormalizable}
		If $\alpha\in \mathcal{M}$ is not contained in an infinite nested sequence 
		$$\mathcal{M}_{\alpha_1}\supsetneq \mathcal{M}_{\alpha_2}\supsetneq \cdots$$
		of small Mandelbrot copies, then $\mathcal{M}$ is locally connected at $\alpha$. If $\alpha$ is contained in an infinite nested sequence of small Mandelbrot copies as above and if $\emph{Diam}\, \mathcal{M}_{\alpha_n}\to 0$ when $n\to 0$, then $\mathcal{M}$ is locally connected at $\alpha$. 
	\end{theorem}

	\begin{figure}
		\begin{center}
			\def\svgwidth{6.4in}
\begingroup%
  \makeatletter%
  \providecommand\color[2][]{%
    \errmessage{(Inkscape) Color is used for the text in Inkscape, but the package 'color.sty' is not loaded}%
    \renewcommand\color[2][]{}%
  }%
  \providecommand\transparent[1]{%
    \errmessage{(Inkscape) Transparency is used (non-zero) for the text in Inkscape, but the package 'transparent.sty' is not loaded}%
    \renewcommand\transparent[1]{}%
  }%
  \providecommand\rotatebox[2]{#2}%
  \newcommand*\fsize{\dimexpr\f@size pt\relax}%
  \newcommand*\lineheight[1]{\fontsize{\fsize}{#1\fsize}\selectfont}%
  \ifx\svgwidth\undefined%
    \setlength{\unitlength}{1243.94923908bp}%
    \ifx\svgscale\undefined%
      \relax%
    \else%
      \setlength{\unitlength}{\unitlength * \real{\svgscale}}%
    \fi%
  \else%
    \setlength{\unitlength}{\svgwidth}%
  \fi%
  \global\let\svgwidth\undefined%
  \global\let\svgscale\undefined%
  \makeatother%
  \begin{picture}(1,0.33003458)%
    \lineheight{1}%
    \setlength\tabcolsep{0pt}%
    \put(0,0){\includegraphics[width=\unitlength]{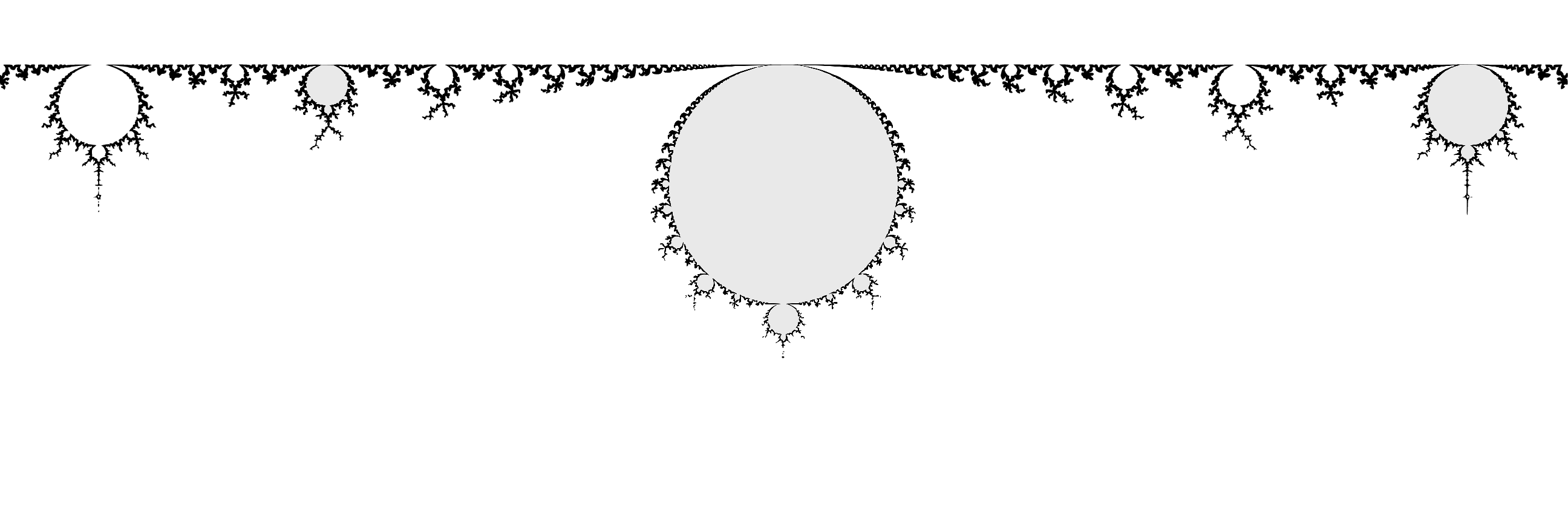}}%
    \put(0.52947723,0.08702299){\makebox(0,0)[lt]{\lineheight{1.25}\smash{\begin{tabular}[t]{l}$\mathcal{L}_{0/1}$\end{tabular}}}}%
    \put(0.87911191,0.17459556){\makebox(0,0)[lt]{\lineheight{1.25}\smash{\begin{tabular}[t]{l}$\mathcal{L}_{1/2}$\end{tabular}}}}%
    \put(0.21394097,0.2041176){\makebox(0,0)[lt]{\lineheight{1.25}\smash{\begin{tabular}[t]{l}$\mathcal{L}_{-1/3}$\end{tabular}}}}%
  \end{picture}%
\endgroup%

			\caption{The boundary of $\mathcal{M}$ with some limbs shaded.}
			\label{M}
		\end{center}
	\end{figure}

	One of the key tools in proving Theorem \ref{yoccoz infinitely renormalizable} is the Pommerenke-Levin-Yoccoz (PLY) inequality, which bounds the multipliers of periodic points of polynomials. For the family of maps $Q_\alpha$, the PLY inequality bounds the size of rational limbs of $\mathcal{M}$:
	\begin{theorem}[Pommerenke \cite{pommerenke},  Levin\cite{levin}, Yoccoz \cite{hubbard}]\label{PLY}
		For any reduced rational number $p/q\in \mathbb{Q}$,
		$$\emph{Diam}\,\mathcal{L}_{p/q}\leq \frac{\log 2}{\pi q}.$$
	\end{theorem}
	
	In \cite{problems}, Milnor conjectured that the $O(1/q)$ bound in Theorem \ref{PLY} can be improved to $O(1/q^2)$. While there has been some improvement on controlling the size of limbs of $\mathcal{M}$, see for example \cite{Levin_multipliers1}, \cite{Levin_multipliers2}, and \cite{pacman},  
	there have been no success in proving an $O(1/q^2)$ PLY inequality for \textit{any} infinite subset of $\mathbb{Q}$. In this paper, we provide the first such examples. Before stating the main result, we need a few definitions.

	For any integer $N\geq 1$ we denote by 
	$\mathbb{Q}_N$ the set of all reduced rational numbers $p/q$ that can be written as a continued fraction
	\begin{equation}\label{continued fraction}
		p/q= \cfrac{\varepsilon_1}{a_1+\cfrac{\varepsilon_2}{\ddots + \cfrac{\varepsilon_N}{a_N}}},
	\end{equation}
	where $\varepsilon_n= \pm 1$ and $a_n\geq 2$ is an integer for all $1\leq n \leq N$. 
	Our main result on the PLY inequality is the following
	\begin{theorem}\label{quadratic PLY}
		For any $N\geq 1$, there exists a constant $C_N>1$ such that 
		$$\frac{C_N^{-1}}{q^2}< \emph{Diam}\, \mathcal{L}_{p/q}<\frac{C_N}{q^2}$$
		for all $p/q\in \mathbb{Q}_N$.
	\end{theorem}
	Note that Theorem \ref{quadratic PLY} also shows that $O(1/q^2)$ is the optimal bound in these cases. 
	Our proof of Theorem \ref{quadratic PLY} will allow us to control the size of even more limbs of $\mathcal{M}$. 	
	We will say that a parabolic parameter $\alpha$ has \textit{combinatorics} $\langle \alpha_j\rangle_{j=1}^N$ if
	$$\mathcal{M}_\alpha =  \tau_{\alpha_1}\circ \cdots \circ \tau_{\alpha_N}(\mathcal{M}_0).$$ 
	We will prove the following generalization of Theorem \ref{quadratic PLY}:
	\begin{theorem}\label{main}
		For any positive integers $N_1$ and $N_2$ there exists $ \epsilon>0$ and $C>1$ such that if $\alpha$ has combinatorics $\langle p_n/q_n\rangle_{n=1}^{N_1}$ with $p_n/q_n\in \bigcup_{m=1}^{N_2}\mathbb{Q}_{m}$ for all $n\geq 1$ and $|p_n/q_n|< \epsilon$ for all $n\geq 2$, 
		then 
		$$ \frac{C^{-1}}{\prod_{n=1}^{N_1}q_n^2}< \emph{Diam}\, \mathcal{L}_{\alpha} < \frac{C}{\prod_{n=1}^{N_1}q_n^2}.$$
	\end{theorem}
	\noindent Note that when $N_1 = 1$, Theorem \ref{main} reduces to Theorem \ref{quadratic PLY}. When $N_1>1$, Theorem \ref{main} bounds the size of nested small satellite Mandelbrot copies in $\mathcal{M}$. The precise control given by Theorem \ref{main} provides a powerful tool in studying the geometry of $\mathcal{M}$; these bounds have recently been used by Lomonaco and Petersen \cite{compatability} to complete a program on the quasiconformal compatibility of of small Mandelbrot copies.

	Our approach to proving Theorem \ref{main} relies on careful analysis of bifurcation of parabolic cycles. 
	Let $f_0$ be a holomorphic function with a parabolic periodic cycle and let $h_0$ be a perturbation of $f_0$. When $h_0$ is close to $f_0$, high iterates of $h_0$ may approximate certain transcendental functions which are called   \textit{Lavaurs maps} for $f$. This phenomenon, called \textit{parabolic implosion}, was first studied by  Douady, Hubbard, and Lavaurs in \cite{Orsay2} and \cite{Lavaurs} and has had numerous applications in holomorphic dynamics, see for example \cite{douady}, \cite{shishikura_boundary},\cite{positive_area}, \cite{MMYconjecture}, and \cite{CS_satellite}.
	In \cite{shishikura_boundary} and \cite{shishikura_1}, Shishikura introduced an alternative description of parabolic implosion in terms of two renormalization operators. 
	Instead of a Lavaurs map $L_1$ for $f_0$, there is a \textit{parabolic renormalization} of $f_0$: a function $f_1$ which is semi-conjugate to $L_1$. Instead of high iterates of $h_0$ approximating $L_1$ there is a \textit{near-parabolic renormalization} of $h_0$: a function $h_1$ which is semi-conjugate to iterates of $h_0$ and which approximates $f_1$. 
	Working with these renormalization operators has a notable upshot: the parabolic renormalization $f_1$ is a  simpler function then the Lavaurs map $L_1$. 
	In particular, there are classes of maps which are invariant under these renormalization operators, see for example \cite{shishikura}, \cite{Yang}, \cite{Cheritat22}, and \cite{nondegenerate}. 
	
	If $f_1$ is a near-parabolic renormalization of $f_0$ and $f_1$ also has a parabolic periodic cycle, then there similarly exist parabolic renormalizations of $f_1$ and near-parabolic renormalizations of $h_1$.  A \textit{parabolic tower} is a (possibly infinite) sequence 
	$\langle f_n \rangle_{n=0}^N$ such that $f_n$ is a parabolic renormalization of $f_{n-1}$ for all $n>0$. 
	We can similarly construct \textit{near-parabolic towers}, sequences $\langle h_n \rangle_{n=0}^N$ such that $h_n$ is close to $f_n$ for all $n\geq 0$ and $h_n$ is a near-parabolic renormalization of $h_n$ for all $n\geq 1$.
	Parabolic towers were first introduced, with a slightly different definition, by Epstein in \cite{epstein_thesis}. 
	
	\begin{figure}
		\begin{center}
			\def\svgwidth{5in}
			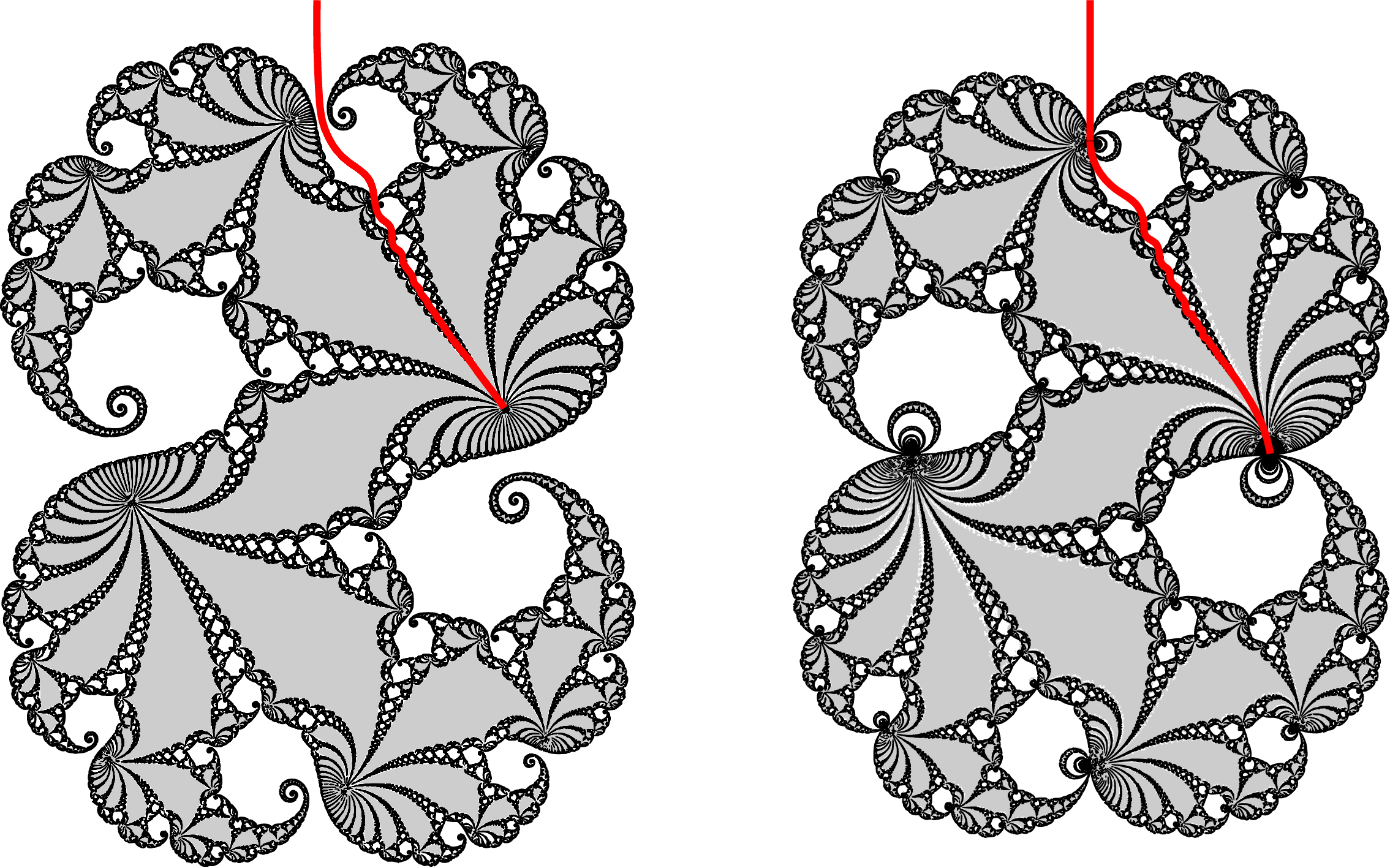
			\caption{\textit{Left:} The filled Julia set for $Q_\alpha$ for $\alpha$ close to one with an external ray drawn in red. \textit{Right:} The analogue of the filled Julia set for a Lavaurs map of $Q_1$ with an analogue of an external ray drawn in red.}
			\label{filled Julia set}
		\end{center}
	\end{figure} 

\begin{figure}
		\begin{center}
			\def\svgwidth{6.4in}
			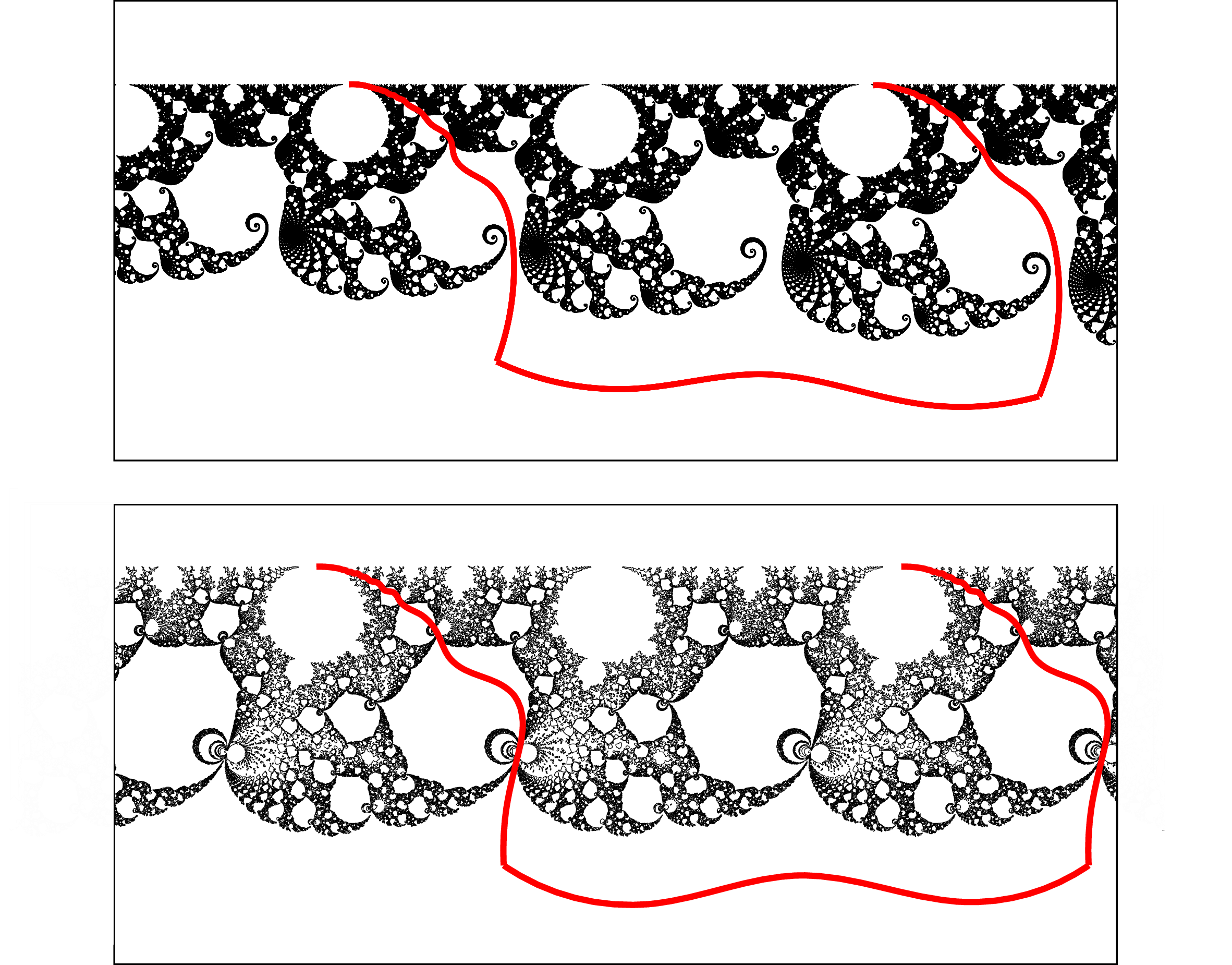
			\caption{A zoomed-in subset of the boundary of $\mathcal{M}$ near $1/n$ for some large $n$  (top) and the boundary of the analogue of $\mathcal{M}$ in the space of Lavaurs maps for $Q_1$ (bottom). In red, certain subsets of external rays of $\mathcal{M}$ which approximate analogues of external rays in the space of Lavaurs maps are drawn in. The bound on the size of these subsets of external rays of $\mathcal{M}$ induces the bound on the size of limbs in Theorem \ref{main}.}
			\label{bounds fig}
		\end{center}
	\end{figure}

	Lavaurs and Douady showed in \cite{Lavaurs} and \cite{douady} respectively that for polynomials, if $h$ is a perturbation of a parabolic map $f$ and if high iterates of $h$ approximate a Lavaurs map $L$ of $f$, then the filled Julia set of $h$ approximates the analog of the filled Julia set for $L$; see for example Figure \ref{filled Julia set}. 
	Our main analysis here is to push this control further: when $f= Q_0$ and $h = Q_\alpha$, the external rays of $h$ approximate analogues of external rays for  $L$. 
	To prove Theorem \ref{main}, we first extend this control on the geometry of external rays to towers.  More precisely, let $\langle f_n\rangle_{n=0}^N$ be a parabolic tower with $f= Q_{\alpha_0}$  and let $\langle h_{n}\rangle_{n=0}^N$ be a near-parabolic tower with $h_0 = Q_\alpha$. While $h_N$ is not a polynomial and has no canonical definition of external rays, the semi-conjugacy from iterates of $h_0$ to $h_N$ induces analogues of external rays for $h_N$. 
	We show that if the towers have {finite height}, that is $N$ is finite, then the external rays of $h_N$ approximate the external rays of Lavaurs maps for $f_N.$
	This bound on the asymptotic geometry of the external rays of $h_N$ then induces a bound on the geometry of the external rays of $\mathcal{M}$ near $\alpha_0$. Using these external rays we produce curves which surround limbs of $\mathcal{M}$ and which have controlled size, leading to Theorem \ref{main} (see Figure \ref{bounds fig}).
	
	It is important in our argument that we only consider finite height parabolic towers;   the integers $N_1$ and $N_2$ in Theorem \ref{main} bound the height of parabolic towers under consideration. However, we expect that our analysis  can be generalized to some infinite height  parabolic towers. If possible, then we could remove the dependence on $N_1$ and $N_2$ for some  combinatorics in Theorem \ref{main}. For any $M\geq 2$, the \textit{high-type numbers} $HT_M$ are the rational numbers $p/q$ with continued fraction expansion as in (\ref{continued fraction}) with $a_n\geq M$ for all $n$. Successfully applying our methods to infinite height towers would prove the following conjecture:
	\begin{conjecture}\label{main conjecture}
		There exists some $M\gg 2$  and $C>1$ such that if $\alpha$ has combinatorics $\langle p_n/q_n\rangle_{n=1}^N$ with $p_n/q_n\in HT_M$ for all $n$, then
		$$ \frac{C^{-1}\cdot 2^{-N}}{\prod_{n=1}^{N}q_n^2}< \emph{Diam}\, \mathcal{L}_{\alpha} < \frac{C\cdot 2^N}{\prod_{n=1}^{N}q_n^2}.$$
	\end{conjecture}
	Conjecture \ref{main conjecture} represents the unification of a quadratic PLY inequality and MLC for satellite combinatorics in the high-type setting. 
	While Conjecture \ref{main conjecture} does not hold without the high-type assumption, indeed the size of the limbs with combinatorics $\langle 1/2, \dots, 1/2\rangle$ can be uniformly bounded away from zero, we may wonder if Conjecture \ref{main conjecture} holds in general for small  Mandelbrot copies instead of limbs. 
	 We can compare Conjecture \ref{main conjecture} with Theorem D in \cite{CS_satellite}, where it is shown that there exists $M\geq 2$, $C>0$, and $\lambda\in (0, 1)$  such that if $\alpha$ has combinatorics  $\langle 1/q_n\rangle_{n=1}^N$ for some $N\geq 1$ with $q_1\geq M$ and $q_{n+1}\geq q_{n}^2$ for all $n\geq 1$, then 
	$$\text{Diam}\, \mathcal{M}_{\alpha}< C\cdot \lambda^N.$$
	While the analysis in \cite{CS_satellite} is related to ours,  both rely on near-parabolic renormalization, Cheraghi and Shishikura use the PLY inequality to prove Theorem D in \cite{CS_satellite}; it is the suboptimality of the PLY inequality which forces the quadratic growth of denominators there. Our approach does not make use of the PLY inequality,  it instead produces a version of the PLY inequality with optimal asymptotics, so no such growth condition appears in Conjecture \ref{main conjecture}.

 	This paper is organized as follows. In $\S$\ref{implosion section} we review the theory of parabolic implosion. In $\S$\ref{lavaurs section} we study the geometry and dynamics of Lavaurs maps and their parameter spaces. In $\S$\ref{rays section} we show how the asymptotic geometry of external rays for maps near a parabolic map are controlled by Lavaurs maps. In $\S$\ref{bounds section} we prove Theorem \ref{main} and discuss possible generalizations such as Conjecture \ref{main conjecture}.\\
 	
 	\noindent \textbf{Acknowledgments:}
	  I am grateful to Xavier Buff, Dzmitry Dudko, Adam Epstein, Sarah Koch, Luna Lomonaco, Mikhail Lyubich, John Hubbard,  and Carsten Petersen for many insightful discussions of elephants and parabolic implosion. 	  
	  Dudko, Lyubich, and Petersen have recently conducted a study of families of horn maps in \cite{taming} which is closely related to $\S$\ref{lavaurs section} of this paper. While this research has been conducted independently, the definition of escaping sets in $\S$\ref{lavaurs section} is very much inspired by the definition of ``Lagoas" in \cite{taming}. Part of this project was carried out during Spring 2022 special semester in Complex Dynamics at the Mathematical Sciences Research Institute, and parts of this project were supported by the Ford Foundation and an NSF postdoctoral fellowship (DMS-2303168).  
	  The results in this paper also appear in the PhD thesis 
	  \cite{kapiamba_thesis}.\\

	\section{Parabolic Implosion}\label{implosion section}
	
	A holomorphic function $f$ has a \textit{$(p/q-)$parabolic $(k-)$periodic} point at zero if $k\geq 1$ is the minimal integer such that $f^k(0) =0$ and if $(f^{k})'(0) = e^{2\pi ip/q}$ for some  reduced rational number $p/q$.	In this section we will review some of the theory on the local dynamics near parabolic periodic points, for more details see \cite{shishikura_1}, \cite{Oudkerk}, or \cite{nondegenerate}.

	We denote by $\mathcal{F}$ the class of holomorphic functions $f:\hat{\mathbb{C}}\dashrightarrow \hat{\mathbb{C}}$ such that:
	\begin{enumerate}
		\item $f(0) = 0$ and $f(\infty)= \infty$.
		\item The restriction $f: f^{-1}({\mathbb{C}}^*)\to {\mathbb{C}}^*$ is a branched covering map with a unique critical value $cv^f$ and all critical points are of local degree 2.
	\end{enumerate}
	For example, $\mathcal{F}$ contains the quadratic polynomial $Q_\alpha(z)= e^{2\pi i\alpha}z+z^2$ for all $\alpha\in \mathbb{C}$.
	We denote by $\mathcal{F}^\flower$ the class of functions $f:\hat{\mathbb{C}}\dashrightarrow \hat{\mathbb{C}}$ such that:
	\begin{enumerate}
		\item $f$ is M\"obius conjugate to an element of  $\mathcal{F}$.
		\item $0$ is a $p/q$-parabolic $k$-periodic point of $f$ for some $p/q$ and $k$.
		\item $f^{nk}(cv^f)\to 0$ when $n\to \infty$, where $cv^f$ is the unique critical value of $f$.
	\end{enumerate}
	For example, if $Q_\alpha$ has a parabolic periodic cycle then $Q_\alpha\in \mathcal{F}^\flower$.
	We will usually focus on the case where $f\in \mathcal{F}^\flower$ has a parabolic \textit{fixed} point, that is where $k=1$ above. The general case of periodic points can be studied similarly by replacing $f$ with $f^k$.

	\subsection{Modified continued fractions}\label{modified continued fractions section}

	For any reduced rational number $p/q\in (-1/2, 1/2]$ there is a unique finite sequence 
	$\langle (a_n, \varepsilon_n)\rangle_{n=1}^N$, called the \textit{modified continued fraction expansion} of $p/q$, such that $$p/q=\cfrac{\varepsilon_1}{a_1+\cfrac{\varepsilon_2}{\ddots+\cfrac{\varepsilon_N}{a_N}}}$$
	and $\varepsilon_n=\pm 1$, $a_n\in \mathbb{Z}_{\geq 2}$ for all $n$. Moreover if $a_n = 2,$ then $\varepsilon_n=+1$. The modified continued fraction of $p/q$ plays an important role for perturbations of $p/q$-parabolic fixed points; let us record  some facts about these sequences. 
	
	We denote by $\mathbb{Q}_N$ the set of all rational numbers whose modified continued fraction expansion has length $N$; we will take the empty sequence as the modified continued fraction expansion of $0/1$ and set $\mathbb{Q}_0= \{0/1\}.$
	For any $p/q\in \mathbb{Q}_N$ with $N>0$, we  will say that the \textit{parent} of $p/q$ is the rational number with modified continued fraction expansion $\langle (a_n, \varepsilon_n)\rangle_{n=1}^{N-1}$. We define the parent of $0/1$ to be $1/0$; here we will abuse notation and not distinguish between the fraction $1/0$ and the pair $(1, 0).$
	For $p/q\neq 0/1$, we can recursively compute $p/q$ from the modified continued fraction expansion using the formula 
	\begin{equation*}
	p= a_Np'+\varepsilon_Np'' \text{ and }q= a_Nq'+\varepsilon_Nq'',
	\end{equation*}
	where $p''/q''$ is the parent of $p'/q'$.
	\begin{proposition}\label{inductive denominators}
		$q\geq 2q'$
	\end{proposition}
	\begin{proof}
		If $p/q=0/1$ then $q' = 0$ so the proposition holds. So we assume that $p/q\neq 0/1$ and that the proposition holds for $p'/q'$.
		If $a_N = 2$, so $\varepsilon_N= +1$, then 
		$$q \geq 2q' +q''\geq 2q'.$$
		If $a_N>2$, then $q\geq 3q' - q'' \geq \frac{5}{2}q' > 2q'.$
	\end{proof}

	For any sequences $\langle s_n\rangle_{n=1}^{N}$ and $\langle t_m\rangle_{m=1}^{M}$ with $N_1\neq \infty$, we denote the concatenation of the sequences by 
	$$\langle s_n\rangle_{n=1}^{N}\oplus\langle t_m\rangle_{m=1}^{M}:= \langle s_1, \dots, s_{N_1}, t_1, \dots t_{M}\rangle.$$
	For any $p_1/q_1$ and $p_2/q_2$ with modified continued fractions expansions $\omega_1$ and $\omega_2$ respectively, we define $p_1/q_1\oplus p_2/q_2$ to be the rational number with modified continued fraction expansion $\omega_1\oplus \omega_2$. 
	
	For any $p/q$ with modified continued fraction $\langle (a_n, \varepsilon_n)\rangle_{n=1}^N$.
	We define the \textit{signature} of $p/q$ to be the value 
	$$\mathfrak{S}(p/q):= (-1)^N\prod_{n=1}^N\varepsilon_n,$$
	and we define the M\"obius transformation
	\begin{equation}\label{definition of Mobius}
		\mu_{p/q}(z):= \frac{p+ \mathfrak{S}(p/q)\cdot p'\cdot z}{q+ \mathfrak{S}(p/q)\cdot q'\cdot z}= \cfrac{\varepsilon_1}{a_1+\cfrac{\varepsilon_2}{\ddots+\cfrac{\varepsilon_N}{a_N+\mathfrak{S}(p/q)\cdot z}}} .
	\end{equation}
	From this definition, we make the following two observations.
	\begin{proposition}\label{adding fractions}
		For any rational $p_1/q_1, p_2/q_2$ in $(-1/2, 1/2]$, if $p/q= p_1/q_1\oplus p_2/q_2$ then 
		$$\mu_{p/q}(z) = \mu_{p_1/q_1}\left(\mathfrak{S}(p_1/q_1)\cdot\mu_{p_2/q_2} (\mathfrak{S}(p_1/q_1)\cdot z)\right).$$
	\end{proposition}

	\begin{proposition}\label{Mobius geometry}
		$$\left|\mu_{p/q}(z)- \mu_{p/q}(w)\right| = \left|\frac{z-w}{(q+ \mathfrak{S}(p/q)q'z)(q+ \mathfrak{S}(p/q)q'w)}\right|.$$
	\end{proposition}

	\subsection{Petals for parabolic maps}

	For  any $z_1, z_2\in \mathbb{C}$ satisfying $\text{Re}(z_2-z_1)> \text{Im}(z_2-z_1)$ we set
	$$\Omega(z_1, z_2):= \{w\in \mathbb{C}: \text{Re}\,z_1-|\text{Im}(w-z_1)|< \text{Re}\,w< \text{Re}\,z_2+|\text{Im}(w-z_2)|\}.$$
	We also allow $z_1$ or $z_2$ to be $\infty$, in those cases we ignore the inequality containing $z_1$ or $z_2$ respectively. 
	We will say that a \textit{strip} in $\mathbb{C}$ is a connected set bounded by two parallel lines, and we will say that a strip contained in $\Omega(z_1, z_2)$ is \textit{maximal} if both components of its boundary non-trivially intersect $\partial \Omega(z_1, z_2)$.
	For any $w\in \mathbb{C}$ and $\delta\in \mathbb{C}$ we set  $T_\delta(w)= w+\delta$.
	
	Fixing now some map $f\in \mathcal{F}^\flower$ with a $p/q$-parabolic fixed point at zero, the local dynamics of $f$ near zero are well-understood:
	\begin{theorem}\label{flower parabolic}
		For any sufficiently large $\mathfrak{M}>0$, there exist Jordan domains $P_{att}^f$ and $P_{rep}^f$ inside $Dom(f^q)$ and analytic isomorphisms
		$$\varphi_{att}^f:\Omega(-1, \infty)\to P_{att}^f \text{ and }\varphi_{rep}^f:\Omega(-\infty, -\mathfrak{M})\to P_{rep}^f $$
		which satisfy:
		\begin{enumerate}
			\item The following geometric conditions:
			\begin{enumerate}
				\item $\varphi_{att}^f(0) = cv^f\notin P_{rep}^f$.
				\item $\varphi_{att}^f(w)\to 0$ and $\varphi_{rep}^f(w)\to 0\to 0$ when $w\to \infty$. 
			\end{enumerate}
		\item  The following dynamical conditions:
			\begin{enumerate}
				\item\label{Fatou uniqueness} Up to post-composition by a translation, $(\varphi_{att}^{f})^{-1}$ and $(\varphi_{rep}^{f})^{-1}$ are the unique univalent maps on $P_{att}^f$ and $P_{rep}^f$ respectively which conjugate $f^q$ to $T_1$.
				\item\label{petal adjacency} For any $w\in \mathbb{C}$ with $|\emph{Re}\,w|< 1$, if $|\emph{Im}\,w|$ is sufficiently large then $f^m\circ \varphi_{rep}^f(w) \in P_{att}^f$, where 
				$0\leq m < q$ is an integer satisfying
				$$mp\equiv \begin{cases}
				-1\mod q &\text{ if }\emph{Im}\,w>0,\\
				0\mod q &\text{ if }\emph{Im}\,w<0.
				\end{cases}$$
				Moreover, 
				$$(\varphi_{rep}^f)^{-1}\circ \varphi_{rep}^f(w)-w\to \frac{1-\mathfrak{S}(p/q)}{2}= \begin{cases}
				0 &\text{ if }\mathfrak{S}(p/q)= +1, \\
				1 &\text{ if }\mathfrak{S}(p/q)= -1.
				\end{cases}$$
				when $\emph{Im}\,w\to +\infty.$
			\end{enumerate}
		\end{enumerate}
	\end{theorem}

	We will fix some choice of $\mathfrak{M} =\mathfrak{M}^f$  for $f$ and denote $\Omega_{att}^f= Dom(\varphi_{att}^f), \Omega_{rep}^f= Dom(\varphi_{rep}^f)$.
	It follows from Theorem \ref{flower parabolic} that there is a unique critical point of $f^q$ on $\partial P_{att}^f$, we will call this critical point $cp^f.$ We define the \textit{lower component} of $P_{att}^f\cap P_{rep}^f$ to be the component containing $\varphi_{rep}^f(w)$ for $|\text{Re}\,w|< 1$ and $-\text{Im}\,w>0$ sufficiently large. 
	
	For any $z$ such that $f^{mq+j}(z) \in P_{att}^f$ for some integers $m\geq 0$ and $0\leq j<q$,  we can define
	$$\rho^{f}(z):=(\varphi_{att}^f)^{-1}\circ f^{mq+j}(z)- m.$$
	The \textit{parabolic basin} $U^f$ of $f$ is defined to be the domain of $\rho^f$, and the \textit{immediate basin} $U_0^f$ of $f$ is the component of $U^f$ containing $cv^f$. The dynamics of $f$ on the immediate basin are well understood.
	\begin{proposition}\label{basin rigidity}
		The restriction of $f^q$ to $U_0^f$ is analytically conjugate to the restriction of $Q_0$ to $U_0^{Q_0}.$ 
	\end{proposition}

	For any $w\in \Omega_{att}^{ f}$ such that $z=\varphi_{rep}^{f}(w)\in Dom(h^{mq})$ for some $m\geq 0$, we can define
	$$\chi^{f}(w+m):= f^{mq}(z).$$
	The \textit{horn map} for $f$ is defined to be the map 
	$$H^{f}(w):= \rho^{f}\circ \chi^{f}.$$
	\begin{proposition}\label{horn map parabolic}
		The horn map $H^{ f}$ is analytic, its domain contains both an upper and a lower half-plane, and $$H^{ f}\circ T_1= T_1\circ H^{ f}.$$
		There exist constants $\mathfrak{c}_\pm^{ f}$ such that $H^{f}(w) - w\to \mathfrak{c}_\pm^{f}$ when $\emph{Im}\, w\to \pm \infty$ respectively, and
		$$\mathfrak{c}_+^{f}=\frac{1-\mathfrak{S}(p/q)}{2}.$$
	\end{proposition}

	Setting $\Exp^\pm(w) = e^{\pm 2\pi i w}$, a \textit{parabolic renormalization} of $f$ is defined to be  
	a function of the form
	\begin{align*}
	\mathcal{R}_\delta^\pm f &:= \Exp^\pm \circ T_{\pm\delta- \mathfrak{c}_{\pm}^{f}}\circ H^{f}\circ (\Exp^\pm)^{-1}= \Exp(\delta)\cdot \mathcal{R}^\pm_0 f
	\end{align*}
	with $\delta\in \mathbb{C}$.
	It follows from Proposition \ref{horn map parabolic} that the domain of $\mathcal{R}_\delta^\pm f$ contains punctured neighborhoods of $0$ and $\infty$, and we can analytically extend
	$\mathcal{R}^\pm_\delta f$ by setting $\mathcal{R}^\pm_\delta f(0) = 0$ and
	$\mathcal{R}^\pm_{\delta} f(\infty) = \infty$. Moreover, this extension satisfies $(\mathcal{R}^\pm_\delta  f)'(0) = \Exp(\delta)$.
	\begin{proposition}
		$\mathcal{R}^\pm_\delta f\in \mathcal{F}$ and $cv^{\mathcal{R}^\pm_\delta f}=\Exp^\pm(\pm\delta- \mathfrak{c}_\pm^{ f})$ for  all $\delta\in \mathbb{C}$. The connected components of $Dom(\mathcal{R}^\pm_\delta f)$ containing $0$ and $\infty$ are simply connected subsets of $\mathbb{C}$ and $\hat{\mathbb{C}}\setminus \{0\}$ respectively; the restriction of $\mathcal{R}^\pm_\delta f $ to either component also belongs to $\mathcal{F}$.
	\end{proposition}
	
	The parabolic renormalizations $\mathcal{R}_\delta^+f$ and $\mathcal{R}_\delta^-f$ are called \textit{top} and \textit{bottom} parabolic renormalizations of $f$ respectively. 
	We will focus our attention on the top parabolic renormalizations and denote $\Exp= \Exp^+$, $\mathcal{R}_\delta f= \mathcal{R}_\delta^+f$.

	\subsection{Petals for near-parabolic maps}
	
	Let us keep $f$ as in the previous subsection.
	Throughout this article we will consider holomorphic maps in the compact-open topology with domains, so a neighborhood of $f$ is of the form
	$$\left\lbrace h:Dom(h)\to \mathbb{C} \;\middle|\;
	\begin{tabular}{@{}l@{}}
	$h \text{ is analytic on }Dom(h)\supset K, \text{ and }$\\
	$|f(z) - h(z)|< \epsilon \text{ for all }z\in K$
	\end{tabular}
	\right\rbrace$$
	for some $\epsilon>0$ and compact set $K\subset Dom(f).$

	Fixing some $\mathfrak{C}\gg1$ we set  $$A=A(\mathfrak{C}):= \left\{x+iy:  0< x< \frac{1}{2\mathfrak{C}}, |y|< \mathfrak{C}x^{2}\right\}.$$
	Let $h$ be a holomorphic function satisfying $h(0) = 0$ and $h'(0) = e^{2\pi i\mu_{p/q}(\alpha)}$ for some $\alpha\in A$.
	If $h$ is sufficiently close to $f$ then there is a unique critical point of $h^q$ close to $cp^f$; we denote this critical point by $cp^{h, f}$ and define $cv^{h, f}= h^q(cp^{h, f}).$
	When $h$ is close to $f$, the local dynamics of $h$ near zero mimics that of $f$.

	\begin{theorem}\label{flower perturbed}
		If $h$ is sufficiently close to $f$, then there exists a  holomorphic function $$\varphi_{att}^{h, f}: \Omega\left(-1, \frac{1}{\alpha}- \mathfrak{M}^f\right)\to Dom(h^q)$$ which satisfies:
		\begin{enumerate}
				\item The following geometric conditions:
				\begin{enumerate}
					\item  $\varphi_{att}^{h, f}(0) = cv^{h, f}$.
					\item  $\varphi_{att}^{h, f}(w)$ tends to $0$ or to a non-zero fixed point $x^{h, f}$ of $h^q$ when $\emph{Im}\, w\to +\infty$ or  $\emph{Im}\, w\to -\infty$ respectively.
					\item $\varphi_{att}^{h,f}$ univalently maps any maximal strip in its domain onto a Jordan domain.
				\end{enumerate}
			\item The following dynamical condition:
			\begin{enumerate}
				\item $\varphi_{att}^{h, f} \circ T_1 = h^q\circ \varphi_{att}^{h, f}.$
			\end{enumerate}
			\item The following continuity conditions:
			\begin{enumerate}
				\item $\varphi_{att}^{h, f}$ depends continuously and holomorphically on $h$.
				\item $\varphi_{att}^{h, f}\to \varphi_{att}^f$ and $\varphi_{rep}^{h, f}:= \varphi_{att}^{h, f}\circ T_{1/\alpha}\to \varphi_{att}^f$ when $h\to f$. 
				\item if $w^h\in Dom(\varphi^{h, f})$ satisfies $w^h\to \infty$ and $w^h-\frac{1}{\alpha}\to \infty$ when $h\to f$, then $\varphi^{h, f}(w^h)\to 0$.
			\end{enumerate}
		\end{enumerate}
	\end{theorem}

	We denote $\Omega_{att}^{h, f}= Dom(\varphi_{att}^{h, f})$ and $\Omega_{rep}^{h, f}= Dom(\varphi_{rep}^{h, f})$. Just as for $f$, we can define partial extensions of $(\varphi_{att}^{h, f})^{-1}$ and $\varphi_{rep}^{h, f}$. For any open set $X$ compactly contained in $Dom(\rho^f)$ there exist integers $m\geq 0$ and $0\leq j < q$ such that 
	$f^{mq+j}(X)\subset P_{att}^f$; if $h$ is sufficiently close to $f$ then
	$$\rho^{h, f}(z):=(\varphi_{att}^{h, f})^{-1}\circ h^{mq+j}(z)-m$$
	is defined for all $z\in X$. Similarly, for any open set $Y$ compactly contained in $Dom(\chi^f)$ there exists an integer $m\geq 0$ such that $T_{-m}(X)\subset \Omega_{rep}^f$; if $h$ is sufficiently close to $f$ then 
	$$\chi^{h, f}(w):= h^{mq}\circ \varphi_{rep}^{h, f}(w)$$
	is defined for all $w\in Y$. 
	\begin{proposition}\label{continuity of horn maps}
		The maps $\rho^{h, f}$ and $\chi^{h, f}$ are analytic and converge to $\rho^f$ and $\chi^f$ respectively when $h\to f$. 
	\end{proposition}
	As a consequence of Proposition \ref{continuity of horn maps}, we can control the dynamics of high iterates of $h$. 
	\begin{corollary}\label{convergence to Lavaurs maps}
		For any $\delta\in \mathbb{C}$, $h^{nq}\to \chi^{ f}\circ T_\delta\circ \rho^{f}$ when $h\to f$ and $n- \frac{1}{\alpha}\to \delta.$
	\end{corollary}
	Maps of the form $\chi^f\circ T_\delta\circ \rho^f$ for some $\delta\in \mathbb{C}$ are called \textit{Lavaurs maps}; they are the primary subject of Section \ref{lavaurs section} below.

	The \textit{horn map for $h$ relative to $f$} is defined by $H^{h, f}:= \rho^{h, f}\circ \chi^{h, f}$ and extended using the function equation
	$$H^{h, f}\circ T_1 = T_1\circ H^{h, f}.$$

\begin{proposition}\label{horn map perturbed}
	The horn map $H^{h, f}$ is well-defined and analytic, depends continuously and holomorphically on $h$, converges $ H^f$ when $h\to f$, and its domain contains both an upper and lower half-plane.  There exist constants $\mathfrak{c}_\pm^{h, f}$ such that $H^{h, f}(w) - w\to \mathfrak{c}_\pm^{h, f}$ when $\emph{Im}\, w\to \pm \infty$ respectively.
	Moreover $\mathfrak{c}_+^{h, f} = \mathfrak{c}_+^f$
	and 
	$$\mathfrak{c}_-^{h, f}-\frac{1}{\alpha}= \frac{2\pi i}{\log (h^q)'(x^{h, f})},$$
	taking the branch of $\log$ with imaginary part in $(-\pi/2, \pi/2].$
\end{proposition}

For $\alpha\in A_+$, a \textit{near-parabolic renormalization} of $h$ relative to $f$ is defined to be a function of the form
$$\mathcal{R}^\pm_f h := \Exp^\pm\circ T_{-1/\alpha}\circ H^{h, f}\circ (\Exp^\pm)^{-1}.$$
Just as the parabolic renormalization, we can analytically extend $\mathcal{R}^\pm_fh$ to $0$ and $\infty$ by setting $\mathcal{R}^\pm_fh(0) =0$ and $\mathcal{R}^\pm_fh(\infty) =\infty$; we can compute the derivative at zero as $$(\mathcal{R}^\pm_fh)'(0) = \Exp^\pm(\mathfrak{c}_\pm^{h,f}-1/\alpha).$$
\begin{proposition}\label{near parabolic renormalization}
	The near-parabolic renormalization $\mathcal{R}^\pm_{f} h$ depends continuously and holomorphically on $h$. Moreover, $\mathcal{R}^\pm_fh\to \mathcal{R}^\pm_{\delta}f$ when $h\to f$ and $\mathfrak{c}_{\pm}^{h, f}-1/\alpha\to \pm\delta \mod 1$.
\end{proposition}

One of the key features of the near-parabolic renormalization of $h$ is its semi-conjugacy to high iterates of $h$.
\begin{proposition}\label{lifting dynamics}
	Fix some $\zeta, \zeta'\in \mathbb{C}$, let $w, w'$ be elements of $\Exp^{-1}(\zeta)$ and $\Exp^{-1}(\zeta')$ respectively inside $\Omega_{rep}^{h, f}$, and set  $z= \varphi_{rep}^{h, f}(w) , z'= \varphi_{rep}^{h, f}(w')$, $\zeta = \Exp(w)$.  If $\mathcal{R}_f^\pm h(\zeta) = \zeta'$, then there exist non-negative integers $m, m'$ such that 
	$h^{m}(w) = h^{m'} (w').$ If additionally $\zeta$ is sufficiently close to either zero or infinity, then we can pick $w$ and $w'$ so that $m'= 0.$
\end{proposition}

The maps  $\mathcal{R}_f^+h$ and $\mathcal{R}_f^-h$ are called the \textit{top} and \textit{bottom} near-parabolic renormalizations of $h$ relative to $f$ respectively. 
Just as for the parabolic renormalizations, we will focus on the top near-parabolic renormalization and denote $\mathcal{R}_fh = \mathcal{R}_f^+h.$

While we have only considered so far perturbations with $\alpha\in A$, similar implosive phenomena occur when $\alpha\in -A$. 
Let us denote by $f^*$ and $h^*$ the conjugates of $f$ and $h$ respectively by $z\mapsto \overline{z}$, so 
$$(f^*)'(0) = \Exp(-p/q) \text{ and }(h^*)'(0) = \Exp\circ \mu_{-p/q}(-\overline{\alpha}).$$
When $\alpha \in -A$ and $h$ is close enough to $f$, we can therefore define the near-parabolic renormalization
$\mathcal{R}_f^\pm h := \mathcal{R}_{f^*}^\pm h^*.$ 
Thus all of our analysis in the $\alpha\in A$ case applies to the $\alpha\in -A$ case; however we note that $h^*$ and $\mathcal{R}_{f^*}^\pm h^*$ depend anti-holomorphically on $\alpha$ when $h$ depends holomorphically on $\alpha$.
	
	\section{Lavaurs maps}\label{lavaurs section}
	
	Let us fix a map $f$ in $\mathcal{F}^\flower$. As in the previous section we will assume that $0$ is a fixed point of $f$, so there is some  rational  $p/q\in [-1/2, 1/2]$ such that $f'(0) = e^{2\pi ip/q}.$ 
	We will say that $f$ has \textit{Jordan  basin} if $U_0^f$ is a Jordan domain. 
	\begin{proposition}\label{Jordan basins}
		If $f$ has Jordan  basin, then the components of  $Dom(\mathcal{R}_0f)$ containing $0$ and $\infty$ are both Jordan domains.
	\end{proposition}
	\begin{proof}
		This proposition is proved in \cite[Theorem 2.31]{Yampolsky} for $p/q= 0/1$, the same argument can be applied in the general case.
	\end{proof}
	
	Let us now assume that $f$ has Jordan  basin. We define $W_+^f$ to be the component of $Dom(H^f)$ which contains an upper half-plane, so $\Exp(W_+^f)\cup\{0\}$ is the component of $Dom(\mathcal{R}_0 f)$ containing zero. It follows from Proposition \ref{Jordan basins} that $\partial W_+^f$ is $T_1$-invariant and homeomorphic to $\mathbb{R}$. 	Let $0 \leq k_+^f< q$ be the integer such that $k_+^fp\equiv -1 \mod q$, it follows from part (\ref{petal adjacency}) of Theorem \ref{flower parabolic}  that
	$f^{k_+^f}\circ \chi^f(W_+^f) = U_0^f.$
	Denoting $\chi^f_+ := f^{k_+^f}\circ \chi^f$, we define the \textit{(upper) $\delta$-Lavaurs map} for $f$ to be 
	$$L_\delta^f:=  \chi_+^f\circ T_\delta \circ \rho^f.$$
	As $0 \leq k_+^f < q$, it follows from the definition of $\rho^f$ that we have the following commutative diagrams:
	\begin{center}
		\begin{tikzcd}
		U_0^f \arrow[rr, "L_\delta^f", dashed] \arrow[d, "T_\delta\circ \rho^f"] &  & U_0^f \arrow[d, "T_\delta\circ \rho^f"] \\
		\mathbb{C} \arrow[d, "\Exp"] \arrow[rr, "T_\delta\circ H^f", dashed]   &  & \mathbb{C} \arrow[d, "\Exp"]          \\
		\mathbb{C} \arrow[rr, "\mathcal{R}_\delta f", dashed]                  &  & \mathbb{C}                           
		\end{tikzcd}
		\, and \,
		\begin{tikzcd}
		W_+^f \arrow[rr, "T_{\delta}\circ H^f", dashed] \arrow[d, "\chi_+^f"] & & W_+^f\arrow[d, "\chi_+^f"] \\
		U_0^f \arrow[rr, "{L_\delta^f}", dashed]                                   &  & U_0^f                         
		\end{tikzcd}
	\end{center}
	Let us note that we could instead define lower $\delta$-Lavaurs map for $f$ to be $\chi^f\circ T_{-\delta- \mathfrak{c}_-^f}\circ \rho^f$, so the Lavaurs map is semi-conjugate to both $\mathcal{R}_\delta^-f$ and the restriction of $T_{-\delta- \mathfrak{c}_-^f}\circ H^f$ to the component of $Dom(H^f)$ which contains a lower half-plane. All of our analysis in this section can be applied similarly to these lower Lavaurs maps.

	For the rest of this section we will usually suppress  dependence on $f$ and $\delta$ in the notation when the choices are clear, so for example $\chi_+ = \chi_+^f$ and $L= L_{\delta}^f$.	It follows from the definition that $$L_\delta\circ f^q = L_{\delta+1} = f^q\circ L_\delta.$$
	So for any (possibly negative) integer $m$, we can define
	$$L_\delta\circ f^{qm}:= L_{\delta+m}.$$
	We equip $\mathbb{Z}^2$ with the lexicographic ordering;  so $(m, n)> (m', n')$ if and only if either $m> m'$ or $m =m'$ and $n> n'$. The function $L_\delta^{d}\circ f^{qm}$ is therefore defined if and only $(d, m)\geq (0, 0)$.
	While $L_\delta$ is not necessarily defined on all of $U_0^f$,  as $Dom(\chi_+)$ contains a left half-plane any $z\in U_0^f$ lies in the domain of $L_\delta\circ f^{-qm}$ when $m\geq 0$ is sufficiently large.

	\subsection{Escaping sets}
	
	We define  $\mathcal{E}_1^f:= \mathbb{C}\setminus \overline{W_+^f}.$ As $\partial W_+$ is locally connected and  $\chi_+(W_+) = U_0$,  we can continuously extend $\chi_+$ to $\partial W_+ = \partial \mathcal{E}_1$ such that 
	$$\chi_+(\partial \mathcal{E}_1) = \partial U_0.$$
	We will say that a point in $\hat{\mathbb{C}}$ is \textit{$0$-nonescaping}, \textit{$0$-Julia}, or \textit{$0$-escaping} for $L$ if it belongs to $U_0$, $\partial U_0$, or $\mathbb{C}\setminus U_0$ respectively.	For any $d\geq 1$, we inductively define a point $z$ to be \text{$d$-nonescaping}, \text{$d$-Julia}, or \text{$d$-escaping} for $L$
	if it is $(d-1)$-nonescaping and 
	$$T_\delta\circ \rho\circ L^{d-1}(z)$$
	belongs to $W_+, \partial W_+,$ or $\mathcal{E}_1$ respectively.
	For all $d\geq 0$, we denote by $K_d^{L}$, $J_d^{L}$, and $E_d^{L}$ the set of all $d$-nonescaping,  $d$-Julia, and $d$-escaping points respectively. 
	We will say that $L$ is $d$-nonescaping,  $d$-Julia, or $d$-escaping  when $cv^f$ is the same respectively for $L$. 
	
	\begin{figure}
		\begin{center}
			\def\svgwidth{2.5in}
			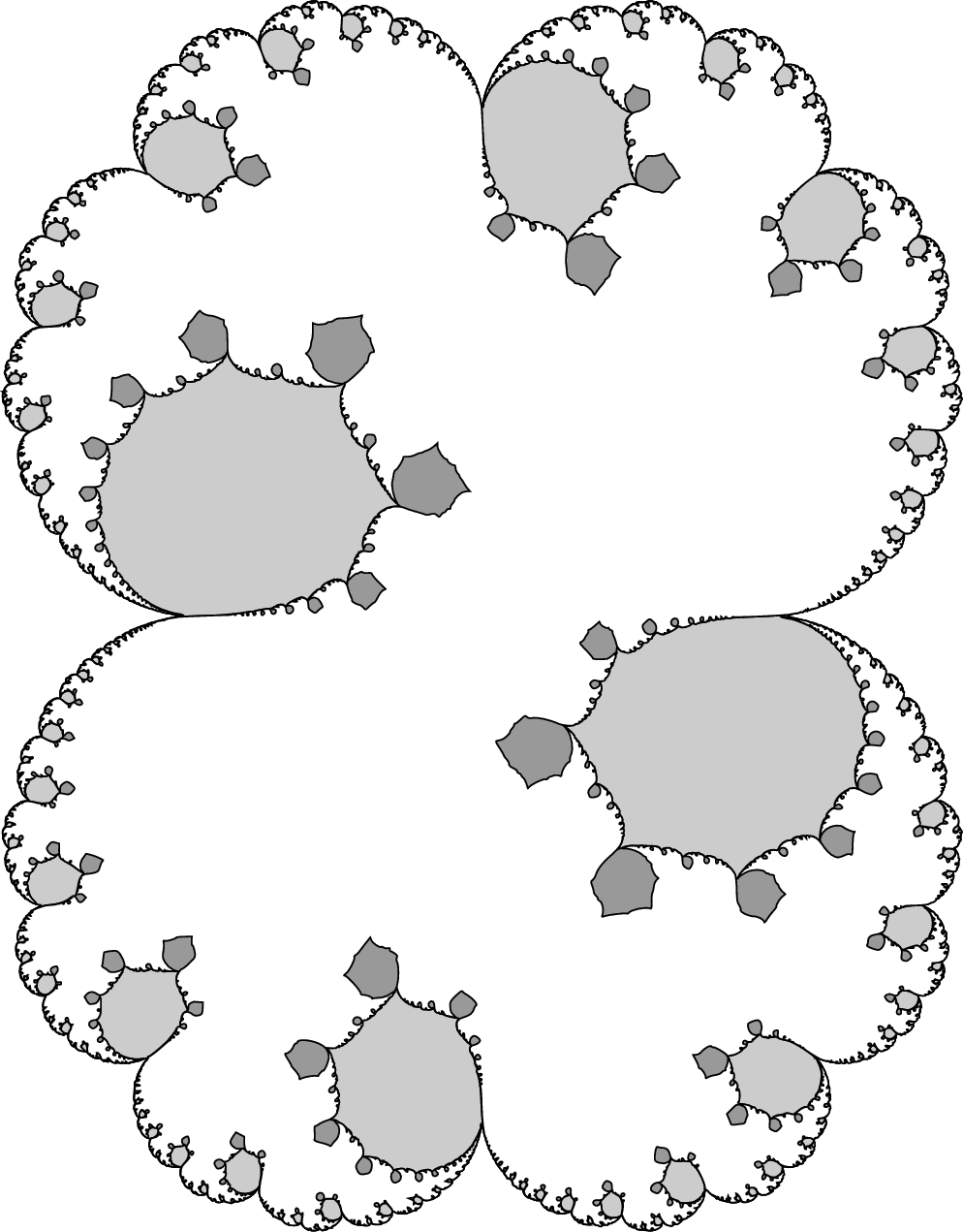
			\caption{The sets $E_1$ (in light gray) and $E_2$ (in dark gray) for
				a $2$-nonescaping Lavaurs map of $Q_1$.}
			\label{escape_regions}
		\end{center}
	\end{figure}
	
	To simplify our notation, we will usually suppress the dependence on $L$ when the choice is clear; for example we will write $K_d= K_d^L$ and $E_d= E_d^L.$ For nonescaping points we can completely characterize the critical points of $L$.
	
	\begin{proposition}\label{critical values}
		For all $d\geq 0$, $z\in K_d$ is a critical point of $L^d\circ f^{mq}$ for some $(d, m)>(0, 0)$ if and only if there exists some $(0, 0)< (d', m')< (d, m)$ such that $L^{d'}\circ f^{m'q}(z) = cv^f$. 
	\end{proposition}
	
	\begin{proof}
		If $L^{d'}\circ f^{m'q}(z) = cv^f$ for some $z\in K_d$, then $L^{d'}\circ f^{(m'-1)q}(z) = cp^f$; thus $z$ is a critical point of $L^d\circ f^{mq}.$ 
		If $z\in K_d$ is a critical point of $L^d\circ f^{mq},$ then either there exists some maximal
		$(0, 0)\leq (d', m')< (d, m)$ such that $z$ is not a critical point of 
		$L^{d'}\circ f^{m'q}$, so in particular $L^{d'}\circ f^{m'q}(z)\in U_0$ is a critical point of $f^q$. Hence $L^{d'}\circ f^{(m'+1)q}(z)=cv^f$.
	\end{proof}

	\begin{remark}
		Our definition of the escaping sets above ignores some of the dynamics of $L$ as some points in $E_1^L$ are mapped into $U_0^f$ by $L$. The upshot of the definition we use here is that it  significantly simplifies the topology of the escaping sets and still allows us to completely describe the dynamics of $L$ that are conjugate to $\mathcal{R}_\delta^f$ near zero. If we wanted to  describe the dynamics of $L$ which are conjugate to $\mathcal{R}_\delta f$ near infinity, or similarly the dynamics of $L$ which are conjugate to a lower parabolic renormalization of $f$ near zero, then we would instead define the escaping sets using the the component of $Dom(H^f)$ which contains a lower half-plane; with this adjustment the rest of our analysis in this section proceeds with minimal change. If we wanted to study the dynamics of $\mathcal{R}_\delta f$ away from zero or infinity, then we could  define the escaping sets of $L$ to be the iterated pre-images of $\hat{\mathbb{C}}\setminus \overline{U^f}$. With such an approach, our arguments in the next section would control the Hausdorff limits of even more general collections of external rays. However we do not need that much control to produce an optimal PLY inequality; we keep with the current definition of escaping sets for simplicity.
	\end{remark}

	\begin{proposition}
		There exists a unique homeomorphism $$\Upsilon^f_0:  \overline{\mathbb{H}}\to \overline{\mathcal{E}_1^f}$$
		which satisfies:
		\begin{enumerate}
			\item $\Upsilon_0$ is univalent on $\mathcal{E}_1$. 
			\item $\Upsilon_0 \circ T_1=  T_1\circ \Upsilon_0$.
			\item $t=0$ is the minimal real number satisfying $ \chi_+\circ \Upsilon_0(t) = 0.$
		\end{enumerate}
	\end{proposition}
	\begin{proof}
		As $Dom_0(\mathcal{R}_\delta f)$ is a Jordan domain in $\mathbb{C}$ containing $0$, there exists an analytic isomorphism from $\mathbb{D}$ to $Dom_0(\mathcal{R}_\delta f)$ which is unique up to pre-composition by a rotation. This isomorphism lifts by $\Exp$ to an analytic isomorphism $\Upsilon_0:  \mathbb{H}\to \mathcal{E}_1$ which is unique up to pre-composition by a translation and  which commutes with $T_1$. As  $\partial\mathcal{E}_1$ is  homeomorphic to $\mathbb{R}$, $\Upsilon_0$ extends continuously to a homeomorphism from $\overline{\mathbb{H}}$ to $\overline{\mathcal{E}_1}$. As $\chi_+$ maps a left half-plane to $P_{rep}^f$  and $\chi_+(\partial \mathcal{E}_1)= \partial U_0$, there exists some minimal $t\in \mathbb{R}$ such that $\chi_+\circ \Upsilon_0(t) = 0$. We can uniquely post-compose $\Upsilon_0$ by a real translation so that $t = 0.$
	\end{proof}

	For any $d\geq 0$, we will say that a point $z\in J_d$ is \textit{$d$-asymptotic} for $L$ if 
	there is some $(d, m)\geq (0, 0)$ such that $L^d\circ f^{mq}(z)= 0$, using the continuous extension of $f^q$ to $\overline{U_0}.$
	We have the following description of  escaping sets (see Figure \ref{escape_regions}).
	
	\begin{proposition}\label{escaping structure}
		If $L$ is $d$-nonescaping for some $d\geq 1$, then:
		\begin{enumerate}
			\item $E_d$ is a countable union of pairwise disjoint Jordan domains, and 
			$T_\delta\circ \rho\circ L^{d-1}$ maps each component univalently onto $\mathcal{E}_1$.
			\item Every component of $E_d$ has a $(d-1)$-asymptotic point on its boundary, moreover that point is the unique  point on the boundary which is not $(d-1)$-nonescaping.
			\item Every $(d-1)$-asymptotic point lies on the boundary of a component of $E_d$.
		\end{enumerate}
	\end{proposition}

	\begin{proof}
		This proposition essentially follows from  the following lemma:
		\begin{lemma}\label{covering of rho}
			For any Jordan domain $X$ in $\hat{\mathbb{C}}$ which has $\infty$ on its boundary and which avoids the negative integers, there exists a unique branch of $\rho^{-1}$ defined on $X$ such that $\rho^{-1}(w)\to 0$ when $w\to \infty$, moreover the image of $X$ is a Jordan domain in $U_0$ whose boundary intersects $\partial U_0$ only at $0$. If $0\in X$, then this branch sends $0$ to $cv^f.$
		\end{lemma}
		\begin{proof}
			For all integers $n$ we define $$\tilde{Y}_n:=\{x+iy: x> {-n}\}\setminus \mathbb{R}_{\leq 1}.$$
			Up to a homotopy rel $\mathbb{Z}$ we may assume that $X$ is contained in $\tilde{Y}:= \bigcup_{n\in \mathbb{Z}}\tilde{Y}_n.$
			For all $n\leq 1$ we can define $Y_n:= \varphi_{att}(\tilde{Y}_n)$.  For any $n>1$ we define $Y_n$ to be the unique component of $\rho^{-1}(\tilde{Y}_n)$ which contains $Y_{n-1}$. Thus $Y_n$ is the unique component of $f^{-q}(Y_{n-1})$ containing $Y_{n-1}$. 
			It follows that  $Y:= \bigcup_{n\in \mathbb{Z}}Y_n$ is the unique $f^q$-invariant component of $\rho^{-1}(\tilde{Y})$, $0\in \partial Y$, and $cv^f\in Y.$
			As every other branch of $\rho^{-1}$ defined on $Y$ corresponds to a component of $f^{-kq}\circ f^{kq}(Y)$, to finish proving the lemma it suffices to show that zero is the unique point in $\partial Y\cap \partial U_0$.
			
			As $Y$ and $\partial U_0$ are $f^q$-invariant, $\partial Y\cap \partial U_0$ must also be $f^q$-invariant. Additionally, $\partial Y\cap \partial U_0$ is connected; indeed if $\langle z_n\rangle_{n=1}^\infty$ and $\langle z_n'\rangle_{n=1}^\infty$ are sequences in $Y$ which tend towards $\partial U_0$ then for all $n$ we can find a curve $\gamma_n$ in $Y$ connecting $z_n$ to $z_n'$ such that any sequence of points in $\gamma_n$ tends to $\partial U_0$ when $n\to \infty$. It is shown in \cite{Orsay1} that the dynamics of the polynomial $Q_0$ on $\partial U_0^{Q_0}$ are topologically conjugate to the doubling map on the circle, so Proposition \ref{basin rigidity} implies that the only $f^q$-invariant connected subsets of $\partial Y\cap \partial U_0$ are $\{0\}$ and $\partial U_0$. 
			Let $z_{-1}$ be the non-zero pre-image under $f^q$ of $0$ in $\partial U_0$.
			As $\gamma:={\varphi_{att}((1, +\infty))}$ is a curve in $Y$ whose closure connects $f^q(cv^f)$ and zero, there are two components $\gamma_+$ and $\gamma_-$ of $f^{-2q}(\gamma)$ which each connect one of the two pre-images of $z_{-1}$ in $\partial U_0$ to $cp^f$. The union $\gamma_+\cup \gamma_-$ is a curve in $U_0$ which avoids $Y$ and  whose closure connects two distinct points in $\partial U_0$, hence $\partial Y\cap \partial U_0\neq \partial U_0$.
		\end{proof}
		If $L$ is $1$-nonescaping, then by definition $\delta\in W_+$. As $\mathcal{E}_1$ is a Jordan domain in $\hat{\mathbb{C}}$, it follows from Lemma \ref{covering of rho} that there is a unique component $B$ of $E_1$ which has zero on its boundary, and $0$ is the point in $J_0$ on the boundary of $B$. Pulling back by $f^q$ completes the proof in the $d=1$ case. For  $d>1$, we observe that pulling back $E_{d-1}$ by $\chi_+$ produces a countable union of pairwise disjoint Jordan domains in $\mathbb{C}$. As $L$ is $d$-nonescaping these domains all avoid $\delta$. As $ \rho$ is a  covering map branched over the negative integers we can pull back by $T_\delta\circ \rho$ and apply Lemma \ref{covering of rho} to complete the proof by induction.
	\end{proof}

	When $L$ is $d$-nonescaping for some $d\geq 1$, we will call the unique $(d-1)$-asymptotic point on the boundary of a component of $E_d$ the \textit{parent} of that component. 
	We will say that that a component $B$ of $E_d$ is a \textit{descendant} of a component $B'$ of $E_{d'}$ with $1\leq d' \leq d$ if either $B = B'$ or the parent of $B$ belongs to the closure of a descendant of $B'$.

	Let $\mathcal{A}= \mathcal{A}^f$ denote the set of all $x\in \mathbb{R}$ such that $\chi_+\circ \Upsilon_0(x)$ is $0$-asymptotic.
	For  $d\geq 0$  we define an \textit{enriched angle of depth $d$} for $f$ to be a sequence $\Theta = \langle \theta_n\rangle_{n=0}^d$ with $\theta_n\in \mathcal{A}$ for all $n$. By abuse of notation, when $d = 0$ we will not differentiate between  $\langle \theta_0\rangle$ and $\theta_0$.
	We will say that another enriched angle $\Theta' = \langle\theta_n'\rangle_{n=0}^{d'}$ is  \textit{equivalent} to $\Theta$, and write $\Theta\sim \Theta'$, if and only if
	\begin{enumerate}
		\item $d = d'$, 
		\item $\theta_n = \theta_n'$ for all $0< n \leq d$, and
		\item $\chi_+\circ \Upsilon_1^{-1}(\theta_0) = \chi_+\circ \Upsilon_1^{-1}(\theta_0')$.
	\end{enumerate}
	For any integer $m$, we also define the following operations on enriched angles:
	\begin{enumerate}
		\item $\Theta\oplus \Theta' := \langle \theta_0, \dots, \theta_d, \theta_0', \dots, \theta_{d'}'\rangle$.
		\item $\Theta+m := \langle \theta_n+m\rangle_{n=0}^d$.
		\item $\Theta*m := \langle \theta_n + m\cdot n\rangle_{n=0}^d$.
		\item $\lfloor \Theta\rfloor:= \langle \theta_n\rangle_{n=0}^{d-1}$ if $d>0$.
		\item $\lceil \Theta\rceil:= \langle \theta_n\rangle_{n=1}^d$ if $d>0$.
	\end{enumerate}
	 We will say that $\Theta$ is \textit{basic} if ${\theta_1}\sim{ 0 }$. 
	 We define the \textit{principle integer} of $\Theta$ to be the minimal integer $j$ so that $(d, j)\geq (0, 0)$ and $\theta_d+j \sim 0$.
	
	For all $\theta\in \mathcal{A}$, we set $z_\theta^f = z_\theta^L:= \chi_+\circ \Upsilon_0(\theta)$. If $L$ is $d$-nonescaping for some $d\geq 1$, then for any enriched angle $\Theta= \langle \theta_n\rangle_{n=0}^d$ we define $B_{\lfloor\Theta\rfloor}^L$ to be the unique component of $\mathcal{E}_d$ which has $z_{\lfloor\Theta\rfloor}^L$ as its parent and and we define $\upsilon_{\lfloor\Theta\rfloor}^L: \mathbb{H}\to B_{\lfloor\Theta\rfloor}^L$ to be the corresponding branch of $(T_\delta \circ \rho\circ L^{d-1})^{-1}\circ \Upsilon_0$. 
	It follows from 
	Proposition \ref{escaping structure} that $\upsilon_{\lfloor\Theta\rfloor}^L$ extends continuously to a homeomorphism from $\overline{\mathbb{H}}\cup\{\infty\}$ to $\overline{B_{\lfloor\Theta\rfloor}^L}$, we define
	$$z_\Theta^L = \upsilon_{\lfloor\Theta\rfloor}^L(\theta_d).$$
	For an enriched angle $\Theta$ of depth $0$, while $\lfloor\Theta\rfloor$ is not defined we will nevertheless denote
	$B_{\lfloor\Theta\rfloor}^L:= \hat{\mathbb{C}}\setminus U_0^f$, $\upsilon_{\lfloor\Theta\rfloor}:= \chi_+\circ \Upsilon_0$, and $z_{\lfloor\Theta\rfloor}^L:= 0$ to simplify some statements.
	
	The two following propositions show how our labeling respects the dynamics of $f$ and $L$.
	\begin{proposition}\label{angle dynamics}
		Assume that $L$ is $d$-nonescaping for some $d\geq 0$. For any enriched angle $\Theta$ of depth $d$,  $f^q(z_\Theta) = z_{\Theta+1}$.  If $d>0$, then  $L(z_\Theta) = z_{\lceil \Theta\rceil}$.
	\end{proposition} 	
	\begin{proof}
		Set $\Theta = \langle \theta_n\rangle_{n=0}^d$.
		If $d =0$, then 
		$$f^q(z_{\Theta} )= f^q\circ \chi_+\circ \Upsilon_0(\theta_0)= \chi_+\circ \Upsilon_0\circ T_1(\theta_0) = z_{\Theta+1}.$$
		If $d>0$ and the proposition holds for smaller values of $d$, then it follows from the definition that
		$f^q\circ \upsilon_{\lfloor\Theta\rfloor} = \upsilon_{\lfloor\Theta\rfloor+1}\circ T_1$. Hence
		$$f^q(z_\Theta) = f^q\circ\upsilon_{\lfloor\Theta\rfloor}(\theta_d) = \upsilon_{\lfloor \Theta\rfloor+1}(\theta_d+1) = z_{\Theta+1}.$$
		If $d = 1$, then 
		$$L(z_\Theta) = \chi_+\circ   { \Upsilon}_{0}(\theta_1)=  z_{\lceil\Theta\rceil}.$$
		If $d>1$ and the proposition holds for smaller values of $d$, then it follows from the definition that $L\circ \upsilon_{\lfloor\Theta\rfloor} = \upsilon_{\lceil\lfloor\Theta\rfloor\rceil} = \upsilon_{\lfloor\lceil\Theta\rceil\rfloor}$. Hence
		$$L(z_\Theta) = L\circ \upsilon_{\lfloor\Theta\rfloor}(\theta_d) = \upsilon_{\lfloor\lceil\Theta\rceil\rfloor}(\theta_d)=z_{\lceil\Theta\rceil}.$$
	\end{proof}

	\begin{proposition}\label{angle dynamics invertible}
		Assume that $L$ is $1$-nonescaping and fix some $\theta\in \mathcal{A}$. For any sufficiently large integers $m, m'\geq 0$,
		$$f^{k_+^f}\circ\varphi_{rep}^f\circ T_{\delta-m'}\circ (\varphi_{att}^f)^{-1}(z_{\langle 0, \theta+m\rangle}) = z_{\theta+m-m'}$$
	\end{proposition}
	\begin{proof}
		It follows from the definitions that $\chi= \varphi_{rep}$ on $\Omega_{rep}$ and $\varphi_{att}$ is a branch of $\rho^{-1}$ defined on $\Omega_{att}$. As $\varphi_{att}(w)\to 0$ when $w\to \infty$, it follows that $\varphi_{att}$ maps $T_{-\delta}(\mathcal{E}_1)\cap \Omega_{att}$ into $B_0$.
		
		Let $m\geq 0$ be a large enough integer so that $z_{\langle 0, \theta+m\rangle} = f^{mq}(z_{\langle 0, \theta\rangle})\in P_{att}.$ It follows from Proposition \ref{angle dynamics} and the definitions of $\rho$ and $\chi_+$ that for any sufficiently integer $m'\geq 0$,
		$$z_{\theta+m-m'}= \chi_+\circ  T_{\delta-m }\circ \rho(z_{\theta+m-m'}) =f^{k_+^f}\circ\varphi_{rep}^f\circ T_{\delta-m'}\circ (\varphi_{att}^f)^{-1}(z_{\langle 0, \theta+m\rangle}).$$
	\end{proof}
	
	When we modify $\delta$ modulo $\mathbb{Z}$, the following propositions show how the labeling changes.
	
	\begin{proposition}\label{translation label change}
		If $L_\delta$ is $d$-nonescaping for some $d\geq 0$, then  $L_{\delta+1}$ is also $d$-nonescaping and $K_d^{L_\delta} = K_d^{L_{\delta+1}}$, $J_d^{L_\delta} = J_d^{L_{\delta+1}}$, and $E_d^{L_\delta} = E_d^{L_{\delta+1}}$. Moreover, 
		$$z_{\Theta}^{L_\delta} = z_{\Theta*1}^{L_{\delta+1}}$$
		for any enriched angle $\Theta$ of depth $d$. 
	\end{proposition}
	
	\begin{proof}
		The first part of the proposition follows immediately from  $L_{\delta+1}= L_\delta\circ f^q.$
		It then follows the definition of $\upsilon_d$ that 
		$$\upsilon_d^{L_{\delta+1}} = T_d\circ \upsilon_d^{L_{\delta}}.$$
		The rest of the proposition then follows immediately.
	\end{proof}

	\begin{proposition}\label{continuous escaping regions}
		If $L_{\delta_0}$ is $(d+1)$-nonescaping, then for any enriched angle $\Theta$ of depth $d$ there exists a neighborhood $\Delta$ of $\delta_0$ such that the function 
		\begin{equation}\label{continuous escaping regions eq}
			(\delta, w)\mapsto \upsilon_{\Theta}^{L_\delta}(w)
		\end{equation}
		is holomorphic on $\Delta\times\mathbb{H}$ and continuous on $\Delta\times \overline{\mathbb{H}}\cup \{\infty\}.$
	\end{proposition}

	\begin{proof}
		First we consider $\Theta= 0$. As $$\upsilon_0^{L_\delta} = \rho^{-1}\circ T_{-\delta}\circ \Upsilon_0$$
		for the unique branch of $\rho^{-1}$ which agrees with $\varphi_{att}$, the proposition follows immediately. Pulling back by specific branches of $f^{-1}$, we can see that the proposition holds for $d=0.$
		
		Now let us assume that $d>0$ and that the proposition holds for smaller values of $d$. The inductive hypothesis implies that the map $$\delta\mapsto \upsilon_{\Theta}^{L_\delta}(\infty)$$
		is continuous on a neighborhood of $\delta_0$. The branch of $(T_\delta\circ \rho\circ L_\delta^{d})^{-1}$ used to define $\upsilon_\Theta^{L_\delta}$ therefore depends continuously and hence holomorphically on $\delta$, so the function (\ref{continuous escaping regions eq}) is holomorphic on $\Delta\times\mathbb{H}$ for some neighborhood $\Delta$ of $\delta_0$. 
		To extend continuously to $\overline{ \mathbb{H}}$, we observe that as $L_\delta$ is $(d+1)$-nonescaping any continuous branch of $(T_\delta\circ \rho\circ L_\delta^{d})^{-1}$ on $\mathcal{E}_1$ can be analytically extended to a neighborhood of $\overline{\mathcal{E}_1}$.
	\end{proof}

	\subsection{Pre-petals}
	
	In the next section it will be important to consider certain pre-images of the attracting and repelling petals for $f$.
	
	\begin{proposition}\label{pre-petals}
		Fix some $d\geq 0$ and assume that $L$ is $d$-nonescaping. There exist $R, R'\geq 0$ such that for any enriched angle $\Theta$ of depth $d$ with principle integer $j$, there exist  univalent functions $\varphi_{\Theta, att}^L$ and $\varphi_{\Theta, rep}^L$ satisfying:
		\begin{enumerate}
			\item The domains of $\varphi_{\Theta, att}^L$ and $\varphi_{\Theta, rep}^L$ are $T_{R-j}(\Omega_{att}^f)\text{ and }T_{-R'-j}(\Omega_{rep}^f)$ respectively. 
			\item  On the domain of $\varphi_{\Theta, \iota}^L$, $$L^{d}\circ f^{jq}\circ \varphi_{\Theta, \iota}^L = \varphi_{\iota}^f\circ T_j.$$
			\item $\varphi_{\Theta, \iota}^L(z)\to z_\Theta^L$ when $w\to \infty$.
			\item  The image of $\varphi_{\Theta, att}^L$ is contained in $K_{d-1}$.
			\item For any $w$ with $|\emph{Re}\,w|< 1$ and $-\emph{Im}\,w>0$ sufficiently large, $\varphi_{\Theta, rep}^L(w)$ is contained in the image of $\varphi_{\Theta, att}^L$.
			\item If $d=0$, then $R=1$ and $R'=0.$
		\end{enumerate} 
	\end{proposition}

	\begin{proof}
		If $\theta\sim 0$, then we define $\varphi_{\theta, \iota}^L:= \varphi_{\iota}^f.$
		As $cv^f$ avoids the Jordan domain $\varphi_{att}\circ T_1(\Omega_{att})$ and $U_0^f$ is a Jordan domain, there exists a unique continuous branch of $f^{-q}$ defined on $\varphi_{att}\circ T_1(\Omega_{att})$ whose image lies in $U_0^f$ and such that $f^{-q}\circ \varphi_{att}(w)\to z_{-1}^L$ when $w\to \infty.$ Using this branch, we define
		$$\varphi_{-1, att}^L:= f^{-q}\circ \varphi_{att}^f\circ T_1$$
		on $\Omega_{att}.$ Let $X$ be the lower component of $P_{att}\cap P_{rep}$.  It follows from Theorem \ref{flower parabolic} that $cv^f$ avoids $P_{rep}$, so there exists a continuous branch of  $f^{-q}$ on $P_{rep}$ which agrees with the above branch on $X$. Using this new branch of $f^{-q}$, we define
		$$\varphi_{-1, rep}^L:= f^{-q}\circ \varphi_{rep}^f\circ T_1$$
		on $T_{-1}(\Omega_{rep}).$
		
		It follows from Proposition \ref{critical values} that for any integer $\ell$, the critical values of $L^d\circ f^{\ell q}$ are contained in the set
		$$\bigcup_{\substack{0\leq s < q\\0\leq n \leq d\\(n, m)\geq (0, 0)}}L^n\circ f^{mq+s}(cv^f).$$
		As $L$ is $d$-nonescaping, it follows that $L^n\circ f^{mq}(cv^f)\to 0$ when $m\to \pm \infty$ for any $1\leq n \leq d$ and $(n, m)\geq (0, 0).$ Moreover, $L^n\circ f^{mq}\notin J_0^L$ for any $(n, m)\geq (0, 0)$. Hence $L^d\circ f^{\ell q}$ is a covering map over a neighborhood $D$ of $z_{-1}^L$ for any integer $\ell$. 
		It follows from Theorem \ref{flower parabolic} that we can choose $R, R'\geq 0$ so that $$\varphi_{-1, att}^L\circ T_{R-1}(\Omega_{att})\cup \varphi_{-1, rep}^L\circ T_{-R'-1}(\Omega_{rep})\subset D$$
		Thus for any  enriched angle $\Theta$ of depth $d$ with principle integer $j>0$, it follows from Proposition \ref{angle dynamics} that there is a unique continuous branch of $(L^d\circ f^{(j-1)q})^{-1}$ defined on $D$ which sends $z_{-1}^L$ to $z_{\Theta}^L$. Using this branch, we  define
		$$\varphi_{\Theta, \iota}^{L}:= (L^d\circ f^{(j-1)q})^{-1}\circ \varphi_{-1, \iota}^{L}\circ T_{j-1}.$$
		If $d =  0$, then the above implies that the image of $\varphi_{-1, \iota}^L$ avoids the critical values of $f^{\ell q}$ for all $\ell>0$, so we can choose $R=1$ and $R'=0.$
	\end{proof}

	We will call the functions $\varphi_{\Theta, att}^{L}$ and $\varphi_{\Theta, rep}^{L}$  attracting and repelling \textit{pre-petal parameters} respectively. We will call the images $$P_{\Theta, \iota}^L:=\varphi_{\Theta, \iota}^L(Dom(\varphi_{\Theta, \iota}^L))$$
	attracting and repelling \textit{pre-petals} for $L$.
	It follows from the proof of Proposition \ref{pre-petals} that for $\Theta$ an enriched angle of depth $0$, the pre-petal parameter $\varphi_{\Theta, \iota}^L$ does not depend on $L$; in this case we will also denote $\varphi_{\Theta, \iota}^f:=\varphi_{\Theta, \iota}^L$.

	\subsection{Virtually parabolic Lavaurs maps}
	
	Let us denote $\tilde{f}:= \mathcal{R}_\delta f$. We will assume in this subsection that $\delta$ is equal to some rational  $\tilde{p}/\tilde{q}\in [-1/2, 1/2]$, so $\tilde{f}$  has a $\tilde{p}/\tilde{q}$ parabolic fixed point at $0$. Our assumption $f$ has Jordan basin induces the same for $\tilde{f}$.
	
	\begin{proposition}\label{Jordan basin}
		$\tilde{f}$ has Jordan basin.
	\end{proposition}

	\begin{proof}
		When $\tilde{p}/\tilde{q} = 0/1$ this fact is proved in \cite[Theorem 4.6]{Yampolsky}, the same argument can be applied in the general case.
	\end{proof}
	
	The Lavaurs map $L$ is said to have a \textit{virtually parabolic} fixed point at $0$ in this case. 
	We can lift the parabolic dynamics of $\tilde{f}$ to $L$.
	
	\begin{proposition}\label{virtually parabolic parameters}
		There exist and univalent maps $$\tilde{\varphi}_{att}^{L}: \Omega_{att}^{\tilde{f}}\to \mathbb{C}\text{ and }\tilde{\varphi}_{rep}^{L}: \Omega_{rep}^{\tilde{f}}\to \mathbb{C}$$ such that:
		\begin{enumerate}
			\item $\Exp\circ T_\delta\circ \rho^f \circ \tilde{\varphi}_\iota^L= \varphi_\iota^{\tilde{f}}.$
			\item $\tilde{\varphi}_{att}^L(0) = cv^f$.
			\item The image of $\tilde{\varphi}_{\iota}^L$ is a Jordan domain inside $U_0^f$ whose boundary intersects $\partial U_0^f$ only at $0$.
			\item $\tilde{\varphi}_{\iota}^L(w)\to 0$ when $w\to \infty$.
			\item $\tilde{\varphi}_\iota^L\circ T_1 = L^{\tilde{q}}\circ f^{-(\tilde{p}+\mathfrak{c}_+^f\tilde{q})q}\circ \tilde{\varphi}_\iota^L.$
		\end{enumerate}
	\end{proposition}

	\begin{proof}
		As $cv^{\tilde{f}} = \Exp(\delta)$, there exists a unique branch of $T_{-\delta}\circ\Exp^{-1}$ on ${P}_{att}^{\tilde{f}}$ which sends $cv^{\tilde{f}}$ to $0$. As ${P}_{att}^{\tilde{f}}$ is a Jordan domain, it follows from Lemma \ref{covering of rho} that there is a continuous branch of 
		$(\Exp\circ T_\delta\circ \rho^f)^{-1}$ on ${P}_{att}^{\tilde{f}}$ so that
		$$\tilde{\varphi}_{att}^{L}:= (\Exp\circ T_\delta\circ \rho^f)^{-1}\circ \varphi_{att}^{\tilde{f}}$$
		satisfies the first four desired properties. We can extend the restriction of this branch to the lower component of ${P}_{att}^{\tilde{f}}\cap {P}_{rep}^{\tilde{f}}$ to all of ${P}_{rep}^{\tilde{f}}$
		 to similarly define $\tilde{\varphi}_{rep}^{L}$.
		 
		 Let $P$ be the image of $\tilde{\varphi}_{att}^L$ and let $\tilde{P}= T_\delta\circ \rho^f(P).$ The definition of $\tilde{f}$ implies that there is some integer $\tilde{m}$ so that $\Exp$ conjugates $T_{\tilde{m}}\circ \left(T_\delta\circ H^f\right)^{\tilde{q}}$ on $\tilde{P}$ to $\tilde{f}^{\tilde{q}}$ on $P_{att}^{\tilde{f}}.$ Thus $L^{\tilde{q}}\circ f^{\tilde{m}q}(P)$ is contained in a component of $(T_{\delta}\circ \rho^f)^{-1}(\tilde{P})$. It follows from the definition of $L$, more specifically the definition of $\chi_+$, that this component is contained in $U_0^f$ and has $0$ on its boundary; this component is therefore $P$. This proves the last property for $\tilde{\varphi}_{att}^L$ and the same argument can be applied for $\tilde{\varphi}_{rep}^L$; the non-empty intersection of the images of $\tilde{\varphi}_{att}^L$ and $\tilde{\varphi}_{att}^L$ implies that we can use the same integer $\tilde{m}$ for both functions. As $H^{f}(w)- w\to \mathfrak{c}_+^f$ when $\text{Im}\,w\to +\infty$, we can compute 
		 $\tilde{m} +\tilde{q}(\delta+\mathfrak{c}_+^f) = 0.$
	\end{proof}
	
	We denote $\tilde{P}_{\iota}^L$ to be the image of $\tilde{\varphi}_\iota^L$. We also define
	$$U^L:= (\Exp\circ T_\delta\circ \rho^f)^{-1}(U^{\tilde{f}})$$
	and define $U_0^L$ to be the unique component of $U^L$ containing $cv^f$. 
	It follows from Lemma \ref{covering of rho} that $U_0^L$ is a Jordan domain in $U_0^f$ with zero on its boundary and $\Exp\circ T_\delta\circ \rho^f$ conjugates $L^{\tilde{q}}\circ f^{-(\tilde{p}+\tilde{q}\mathfrak{c}_+^f)q}$ on $U_0^L$ to $\tilde{f}^{\tilde{q}}$ on $U_0^{\tilde{f}}.$ 
	Using the virtually parabolic dynamics of $L$, we have the following control of the escaping bubbles of $L$.

	\begin{proposition}\label{parabolic bubble ray}
		There is a unique sequence $(\Theta_d)_{d=0}^\infty$, where $\Theta_d$ is a basic enriched angle of depth $d$ for all $d\geq 0$, such that:
		\begin{enumerate}
			\item For all $d\geq 1$, $\lfloor\Theta_d\rfloor = \Theta_{d-1}$.
			\item For all $d\geq \tilde{q}$, 
			$$L^{\tilde{q}}\circ f^{\tilde{m}q}(B_{\Theta_d}) = B_{\Theta_{d-\tilde{q}}}.$$
			\item For any sufficiently large $d_0\geq 0$, 
			\begin{enumerate}
				\item $\overline{B_{\Theta_d}}\subset \tilde{\varphi}_{rep}^L(\Omega_{rep}^{\tilde{f}})$ for all $d\geq d_0$.
				\item\label{descendents of bubbles} If $$z\in \tilde{P}_{rep}^L\cap \left(\bigcup_{d\geq 0}E_d\right)$$  is sufficiently close to $0$, then $z$ belongs to a descendant of $B_{\Theta_{d_0}}$.
			\end{enumerate}
		\end{enumerate}
	\end{proposition}
	
	\begin{proof}
		For all $0 \leq n < \tilde{q}$ we set $U_n^{\tilde{f}}:= {\tilde{f}}^{n}(U_0^{\tilde{f}}).$ It follows from Proposition \ref{Jordan basin} that $U_n^{\tilde{f}}$ is a Jordan domain for all $n$. As $cv^{\tilde{f}}\in U_0^{\tilde{f}}$ is the unique critical value of $\tilde{f}$, the restriction of $\tilde{f}^n$ to $U_0^{\tilde{f}}$ extends to a homeomorphism from $\overline{U_0^{\tilde{f}}}$ to $\overline{U_n^{\tilde{f}}}$. 
		As $U_0^{\tilde{f}}$ is a Jordan domain, Proposition \ref{basin rigidity} implies that there is a homeomorphism $\xi: \overline{U_0^{Q_1}}\to \overline{U_0^{\tilde{f}}}$ which conjugates $Q_1$ to ${\tilde{f}}^{\tilde{q}}$ and which is analytic on $U_0^{Q_1}$. For all $0\leq n < \tilde{q}$ we define $\gamma_n$ to be the arc ${\tilde{f}}^n\circ \xi([-1, 0])$.
		As $0$ and $\infty$ are the only possible asymptotic values for $\tilde{f}$ and $\tilde{f}^n$ is univalent on $U_0^{\tilde{f}}$ for all $1\leq n< \tilde{q}$, we have
		$$\lim_{t\to -1}\tilde{f}^n\circ\xi(t) = 0$$
		if and only if $n = \tilde{q}$ for  $0\leq n \leq \tilde{q}$. Hence 
		the curve $\gamma_{n}$ contains a point on the boundary of $Dom_0({\tilde{f}})$ if and only if $n = \tilde{q}$, and this point is unique when it exists. As $Dom_0({\tilde{f}})$ is a Jordan domain, we can define  $\Gamma$ to be the union of  $\bigcup_{n=0}^{\tilde{q}-1}\gamma_n$ and a simple curve  which connects the point in $\gamma_{\tilde{q}-1}\cap \partial Dom_0({\tilde{f}})$ to $\infty$,   avoids $Dom_0({\tilde{f}})$, and which is contained in $\mathbb{R}$ near $\infty$. It follows from the construction that the the image of $\Gamma\cap Dom_0({\tilde{f}})$ under ${\tilde{f}}$ is contained in $\Gamma$. We define $\mathcal{V} :=\mathbb{C}\setminus \Gamma$.
		\begin{lemma}
			$\mathcal{V}$ is simply connected.
		\end{lemma}
		\begin{proof}
			By construction, if $\mathcal{V}$ is not simply connected, then there must be a point $z\in \mathbb{C}$ which is contained in both $\gamma_{n_1}$ and $\gamma_{n_2}$ for some $0\leq n_1< n_2< \tilde{ q}.$ Thus there  exists some $t_1, t_2\in [-1, 0)$ such that
			$$z=\tilde{f}^{n_1}\circ \xi(t_1) = \tilde{f}^{n_2}\circ \xi(t_2).$$
			If $t_1\neq -1$ or $t_2 \neq -1$, so $z\in U_{n_1}^{\tilde{f}}$ or $z\in U_{n_2}^{\tilde{f}}$, then it would follow that $n_1 = n_2$ which is a contradiction.  
			As ${\gamma}_n$ contains a point in $\partial Dom_0(\tilde{f})$ if and only if $n = \tilde{q}-1$, we must have $\tilde{n_2}< \tilde{q}-1.$
			As 
			$\tilde{f}^{\tilde{n}}\circ \xi(-1)$ belongs to $Dom_0(\tilde{f})$ for all $0 \leq \tilde{n}< \tilde{q}-1$, we have 
			$$\tilde{f}^{n_2-n_1}(\tilde{f}^{n_1}\circ \xi(-1))= \tilde{f}^{n_2}\circ \xi(-1).$$
			Hence  $\tilde{f}^{n_1}\circ \xi(-1)$
			is periodic under $\tilde{f}$, but this is a contradiction as $$\lim_{t\to -1}\tilde{f}^{\tilde{q}-n_1}(\tilde{f}^{n_1}\circ \xi(t)) = 0.$$
		\end{proof}
		
		It follows from our construction that the intersection of any small neighborhood of zero with $\mathcal{V}$ has exactly $\tilde{q}$ components which zero on their boundary.
		As $\mathcal{V}$ avoids all the critical values of ${\tilde{f}}^{\tilde{q}}$ and ${\tilde{f}}^{\tilde{q}}$ is injective in a neighborhood of $0$, there are exactly $\tilde{q}$ components of ${\tilde{f}}^{-\tilde{q}}(\mathcal{V})$ which have $0$ on their boundary. Let us fix one such component $\mathcal{U}$ and let $\zeta: \mathcal{V}\to \mathcal{U}$ be the corresponding branch of ${\tilde{f}}^{-\tilde{q}}$. As $0$ is on the boundary of $\mathcal{U}$, the Denjoy-Wolff theorem implies that the iterates of $\zeta$ converge locally uniformly to $0$. As $\mathcal{V}$ is connected, there exists some integer $0 \leq n < \tilde{q}$ such that for any $z\in \mathcal{V}$, $\zeta^j(z)$ lies in ${\tilde{f}}^n(P_{rep}^{\tilde{f}})$ for all large $j>0$. Replacing $\mathcal{U}$ with another component of ${\tilde{f}}^{-\tilde{q}}(\mathcal{V})$, we can assume that $n=0$; in particular $P_{rep}^{\tilde{f}}\subset\mathcal{V}.$
		We define $\tilde{\mathcal{V}}$ to be the component of $(\Exp\circ T_\delta\circ \rho^f)^{-1}(\mathcal{V})$ which contains $\tilde{P}_{rep}^L$ and  we define $\tilde{\mathcal{U}}$ to be the component of $(\Exp\circ T_\delta\circ \rho^f)^{-1}(\mathcal{U})$ which contains $\tilde{ \mathcal{V}}.$ 
		Setting $$\tilde{m} =-\tilde{q}(\delta+\mathfrak{c}_+^f),$$
		it follows from Proposition \ref{virtually parabolic parameters} that  $\Exp\circ T_\delta\circ \rho^f$ conjugates the restriction $L^{\tilde{q}}\circ f^{\tilde{m}q}: \tilde{ \mathcal{U}}\to \tilde{ \mathcal{V}}$ to the restriction ${\tilde{f}}^{\tilde{q}}: \mathcal{U}\to \mathcal{V}$. 
		Let $\tilde{\zeta}: \tilde{ \mathcal{V}}\to \tilde{ \mathcal{U}}$ be the corresponding inverse map, so $$\Exp\circ T_\delta\circ \rho^f\circ \tilde{\zeta}= \zeta\circ \Exp\circ T_\delta\circ \rho^f.$$
		
		As $\mathcal{V}$ intersects the compliment of $Dom_0({\tilde{f}})$ in a single component, it follows that $\tilde{\mathcal{V}}$ intersects $E_1$ in a single connected component. 
		For ever integer $s\geq 0$, let $B_{s\tilde{q}}$ be the unique component of $E_{1+s\tilde{q}}$ which contains $\tilde{ \zeta}^s(\tilde{ \mathcal{V}}\cap E_1^L).$
		
		\begin{lemma}\label{bubble lemma}
			Fix some integers $s, s'\geq 0$ and a component $B$ of $E_{s'}$ with parent $z_B$. If $\overline{B}\cap \tilde{\zeta}^s(\tilde{\mathcal{V}})$ is non-empty, then one of the following holds:
			\begin{enumerate}
				\item $s' = s\tilde{q}$ and $B = B_{s\tilde{q}}$.
				\item $s' > 1+s\tilde{q}$ and $\overline{B}\subset \tilde{\zeta}^s(\tilde{\mathcal{V}})$.
				\item $s' = 1+s\tilde{q}$, $\overline{B}\setminus \{z_B\}\subset \tilde{\zeta}^s(\tilde{\mathcal{V}})$, and $z_B\in \partial B_{s\tilde{q}}$.
			\end{enumerate}
		\end{lemma}
		\begin{proof}
			As $B$ and $\tilde{\zeta}^s(\tilde{\mathcal{V}})$ are Jordan domains, if $\overline{B}\cap \tilde{ \zeta}^s(\tilde{\mathcal{V}})$ is non-empty, then ${B}\cap \tilde{\zeta}^s(\tilde{\mathcal{V}})$ is also non-empty. 
			As $L^{\tilde{q}}\circ f^{\tilde{m}q} = \tilde{ \zeta}^{-1}$, we must have $s' \geq s\tilde{q}$. Moreover, the uniqueness in the definition of $B_{s\tilde{q}}$ implies that $B= B_{s\tilde{q}}$ if $s' = s\tilde{q}$.
			Thus we can restrict to the case $s' > s\tilde{q}$.
			
			Assume that there is some point $$x'\in \overline{B}\subset J_{s'}\cup E_{s'}\cup J_{s'-1}$$ which belongs to the boundary of  $ \tilde{\zeta}^s(\tilde{\mathcal{V}})$. As $s'> s\tilde{q}$, the point $x:= (L^{\tilde{q}}\circ f^{mq})^s(x')$ is defined and belongs to $\partial \tilde{\mathcal{V}}$. We have four possible cases:
			\begin{enumerate}
				\item $x\notin U_0^f$,
				\item $\Exp\circ T_\delta\circ \rho^f(x)\notin \overline{Dom_0({\tilde{f}})}$, 
				\item $\Exp\circ T_\delta\circ \rho^f(x)\in \Gamma\cap \partial Dom_0({\tilde{f}})$, 
				or
				\item $\Exp\circ T_\delta\circ \rho^f(x)\in \Gamma\cap Dom_0({\tilde{f}})$.
			\end{enumerate}
			The first case above implies that $x\in J_{s\tilde{q}}^L$, and the second case implies that $x' \in E_{1+s\tilde{q}}$. Both are impossible as $s' > s\tilde{q}$.
			The third case above implies that $x\in J_1\cap \partial B_{0}$, so
			$x'\in J_{1+s\tilde{q}}\cap \partial B_{s\tilde{q}}$. Hence $x' = z_B$ and $s' = 1+s\tilde{q}$.
			The fourth case above is impossible as every point in $\Gamma \cap Dom_0({\tilde{f}})$ has infinite forward orbit under ${\tilde{f}}$. 
			
			Thus we have shown that if $s> s\tilde{q}$, then $\overline{B}$ can intersect the boundary of $\tilde{ \zeta}^s(\tilde{\mathcal{V}})$ only at the parent of $B$ and only when $s = 1+s\tilde{q}$, which completes the proof. 
		\end{proof}
		For all integers $d\geq 0$ not equivalent to $0$ modulo $\tilde{q}$, we inductively define $B_{d}$ to be the unique component of $E_{d+1}$ whose closure contains the parent of $B_{d+1}$. 
		If $d>0$ is equivalent to $0$ modulo $\tilde{q}$, Lemma  \ref{bubble lemma} guarantees that $B_{d+1}$ is a descendant of $B_d$. Lemma \ref{bubble lemma} also ensures that $\overline{B_d}\subset \tilde{\mathcal{V}}$ for all $d>2$, so $\zeta^n(\overline{B_d})$ converges to $0$ inside $\tilde{P}_{rep}^L$ when $n\to \infty$. Hence $\overline{B_d}\subset \tilde{P}_{rep}^L$ for all sufficiently large $d\geq 0.$ 
		Now let us fix any $s\geq 0$ and  let $z$ be a point in $\tilde{P}_{rep}^L\cap(\bigcup_{d\geq 0} E_d)$. As $\tilde{P}_{rep}^L\subset \tilde{ \mathcal{V}}$, if $z$ is sufficiently close to zero then $z\in \tilde{\zeta}^s(\tilde{ \mathcal{V}}).$ If $B$ is the component of $\bigcup_{d\geq 1}E_d$ which contains  $z$, then Lemma \ref{bubble lemma} implies that $B$ is a descendant of $B_{s\tilde{q}}.$
	\end{proof}
	
	We will call the sequence $\langle \Theta_d\rangle_{d=0}^\infty$  in Proposition \ref{parabolic bubble ray} the \textit{parabolic  enriched angles} for $L$.

	\subsection{Parameter spaces} 
	Let us fix some $f$ as in the previous section.
	For any $d\geq 0$, we define the \textit{$d$-nonescaping parameter set}, \textit{$d$-Julia parameter set}, and the \textit{$d$-escaping parameter set} to be 
	\begin{align*}
	\mathcal{K}_{d}^f&:= \left\{\delta\in \mathbb{C}: L_\delta^f \text{ is $d$-nonescaping}\right\},\\
	\mathcal{J}_{d}^f&:= \left\{\delta\in \mathbb{C}: L_\delta^f \text{ is $d$-Julia}\right\}, \text{ and}\\
	\mathcal{E}_{d}^f&:= \left\{\delta\in \mathbb{C}: L_\delta^f \text{ is $d$-escaping}\right\}
	\end{align*}
	respectively. 
	The following proposition shows that for $d=0$, the definition of $\mathcal{E}_0^f$ agrees with ${\mathbb{C}}\setminus \overline{W_+^f}$. 
	\begin{proposition}
		For any $w\in \mathbb{H}$, $\delta = \Upsilon_0^f(w)$ if and only if 
		$T_\delta\circ \rho^f(cv^f) = \Upsilon_0^f(w)$.
	\end{proposition}
	\begin{proof}
		As $\rho^f(cv^f) = 0$, we have $T_\delta\circ \rho^f(cv^f) = \delta$.
	\end{proof}
	
	We set $\mathcal{B}_0^f:= \mathcal{E}_1^f$, and for any basic enriched angle $\Theta$ of depth $d\geq 1$ we define
	$$\mathcal{B}_{\Theta}^f:= \left\{\delta\in \mathcal{E}_{d+1}^f: L_\delta^f\circ f^{-jq}(cv^f)\in B_{\lceil\Theta\rceil-j}^{L_{\delta}^f}\text{ for all large integers }j\geq 0\right\}.$$
	We will say that a parameter $\delta\in \mathcal{J}_{d}^f$ is $d$-asymptotic if $cv^f$ is $d$-asymptotic for $L_{\delta}^f$.  The structure of the escaping parameter sets mimics that of the dynamical escaping sets.

	\begin{figure}
		\begin{center}
			\def\svgwidth{6.5in}
			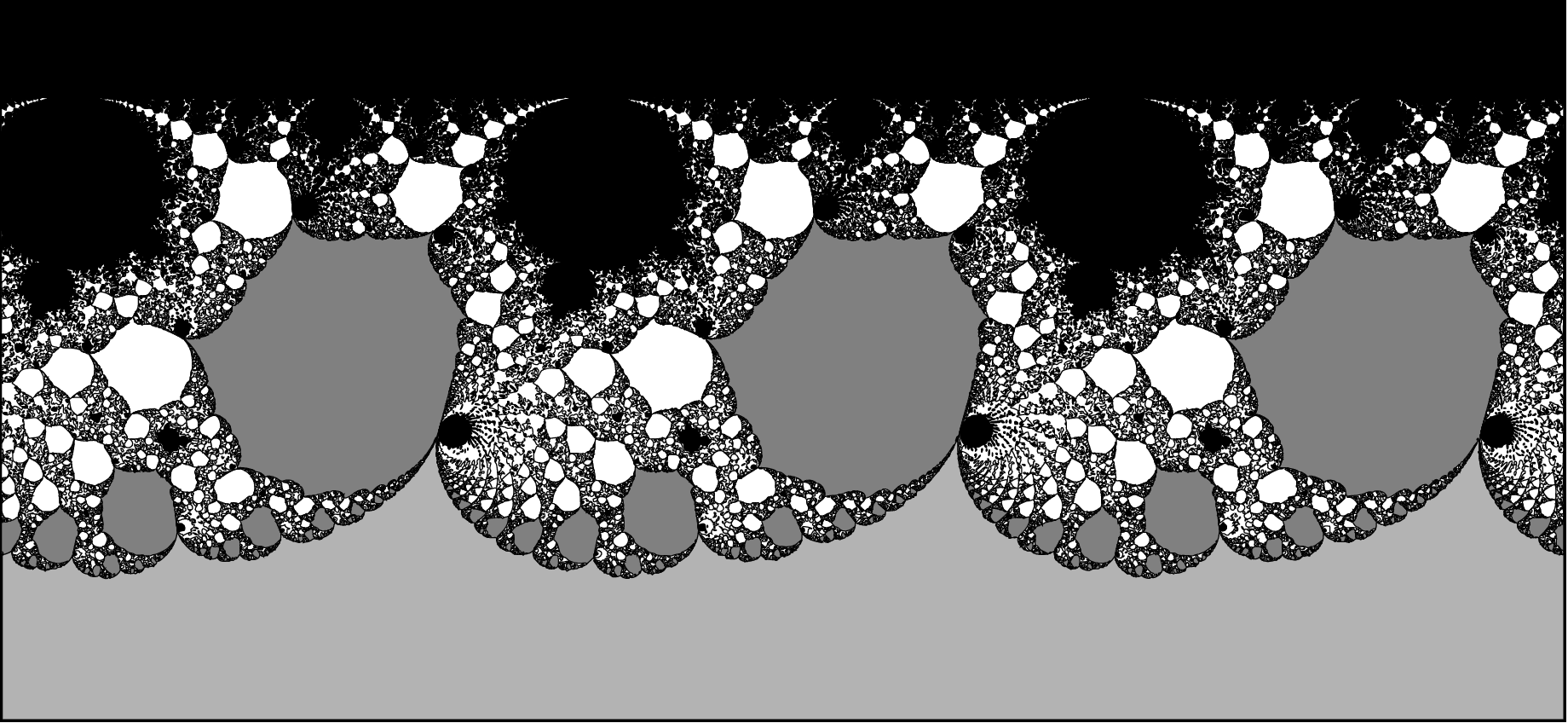
			\caption{The parameter space of Lavaurs maps of $Q_1$. The $d$-escaping parameter set is shown in light gray for $d=1$, dark gray for $d=2$, and white for $d>2$.}
			\label{parameter space pic}
		\end{center}
	\end{figure}

	\begin{proposition}\label{parameter bubbles}
		For any basic enriched angle $\Theta$ of depth $d\geq 1$, 
		$\mathcal{B}_{\Theta}$ is a Jordan domain and there is an analytic isomorphism $\Upsilon_{\Theta}^f: \mathbb{H}\to \mathcal{B}_{\Theta}$
		such that $\delta = \Upsilon_\Theta^f(w)$ if and only if $$L_\delta^f\circ f^{-jq}(cv^f) = \upsilon_{\lceil\Theta\rceil-j}^{L_\delta^f}\circ T_{-j}(w)$$ for all sufficiently large integers $j\geq 0$. Moreover, the components of $\mathcal{E}_{d}$ are all pairwise disjoint.
	\end{proposition}
	
	\begin{proof}
		This proposition is proved in \cite{taming} for the family of horn maps $T_\delta\circ H$. The semi-conjugacies between $T_\delta\circ H$ and $L_\delta$ on $W_+$ implies that the same holds for our family of Lavaurs maps. We will only show here the construction of the map $\Psi_\Theta$, where our argument differs slightly from that in \cite{taming}; we can then prove that  $\mathcal{B}_\Theta$ is a Jordan domain by identical argument to the one given in \cite{taming}. Our argument relies on a quasiconformal surgery, for an introduction to the use of quasiconformal maps in holomorphic dynamics we refer the reader to \cite{branner_fagella}. 
		
		Let us fix some enriched angle $\Theta$ of depth $d\geq 1$ and assume that there is a non-empty connected component $\mathcal{B}$ of $\mathcal{B}_\Theta$, so $\mathcal{B}$ is a connected component of $\mathcal{E}_{d+1}$. 
		 Fix some $\delta\in \mathcal{B}$ and set $$\zeta = T_{\delta}\circ \rho\circ  L_{\delta}^{d}(cv^f).$$
		For all $w\in \mathbb{H}$ we define the quasiconformal map 
		$$\Xi_w: \mathbb{H}\to \mathbb{H}$$ by
		$$\Xi_w(x-iy) = x- y\cdot \frac{w- \text{Re}\,\zeta}{\text{Im}\,\zeta}$$
		for all $x\in \mathbb{R}$ and $y>0$. 
		It follows from this definition that $\Xi_w$ commutes with $T_1$,  extends continuously to the identity on $\mathbb{R}$, and satisfies $\xi_w(\zeta)= w.$
		We set $\xi_w:= \Upsilon_0\circ \Xi_w\circ \Upsilon_0^{-1}$.

		Let us denote by $\Lambda$ the trivial Beltrami differential on $\mathbb{C}.$
		We define a $(T_{\delta}\circ H)$-invariant Beltrami differential on $\mathbb{C}$ by
		$$\tilde{\lambda}_w:= \begin{cases}
		(\xi_w\circ (T_{\delta}\circ H)^d)^*\,\Lambda &\text{ on }(T_{{\delta}}\circ H)^{-d}(\mathcal{E}_1)\text{ for all }d\geq 0,\\
		\Lambda &\text{ on }\mathbb{C}\setminus \bigcup_{d\geq 0}(T_{{\delta}}\circ H)^{-d}(\mathcal{E}_1),
		\end{cases}$$
		where $^*$ denotes the pull-back operation on Beltrami differentials.
		Note that $\tilde{\lambda}_w$ is $T_1$-invariant as $T_\delta\circ H$ commutes with $T_1$.
		
		We define the Beltrami differential $\lambda_w:= (T_{\delta}\circ \rho)^*\,\tilde{\lambda}_w$ on $U_0$; it follows from the semi-conjugacy  of $f^q$ to $T_1$ that $\lambda_w$ is invariant under $f^q$.
		By the measurable Riemann mapping theorem, there exists a quasiconformal homeomorphism $\phi_w: U_0\to U_0$  satisfying $\phi_w^*\,\Lambda = \lambda_w.$
		The map $$\tilde{f}_w:=\phi_w\circ f^q\circ \phi_w^{-1}:U_0\to U_0$$ 
		is therefore holomorphic and quasiconformally conjugate to $f^q$. 
		We have the following classical fact:

		\begin{lemma}\label{rigidity}
			If $g: \mathbb{D}\to \mathbb{D}$ is a holomorphic function satisfying: 
			\begin{enumerate}
				\item zero is the unique critical point of $g$ counting multiplicity; and
				\item there is a bi-infinite sequence $\{z_n\}_{n\in \mathbb{Z}}$ such that $g(z_n)= z_{n+1}$ for all $n$ and $z_n\to 1$ when $n\to \pm \infty$;
			\end{enumerate}
			then $$g(z)= \frac{3z^2+1}{z^2+3}.$$
		\end{lemma}
		\begin{proof}
			See for example \cite[Theroem 2.33]{Yampolsky}.
		\end{proof}

		It follows from Lemma \ref{rigidity} that $\tilde{f}_w$ is conformally conjugate to $f^q$. Hence we can  choose $\phi_w$ so that it depends continuously on $w$ and $\tilde{f}_w = f^q$.
		As $f^q$ has a unique critical point in $U_0$ and $\phi_w$ is orientation-preserving, it follows that $\phi_w$ is the identity on the entire grand orbit of $cv^f$ under $f^q$. As $U_0$ is a Jordan domain,  $\phi_w$ therefore extends continuously to the identity on $\partial U_0$.

		As $T_{\delta}\circ H = T_{\delta}\circ \rho\circ \chi_+$
		on $W_+$,  the Beltrami differential $\chi_+^*\, \lambda_w$ is equal to $\tilde{\lambda}_w$
		on $W_+$.
		We can lift $\phi_w$ by $\chi_+$ to a quasiconformal homeomorphism $\tilde{\phi}_w: W_+\to W_+$; it follows from the above that $\tilde{\phi}_w^*\,\Lambda= \tilde{\lambda}_w$ on $W_+$ and $\tilde{\phi}_w$ extends continuously to the identity on $\partial W_+$. Thus we can extend $\tilde{\phi}_w$ to a quasiconformal automorphism of $\mathbb{C}$ by  $$\tilde{\phi}_w(z) = \begin{cases}
		\tilde{\phi}_w(z) &\text{ if } z\in \overline{W_+},\\
		\xi_w(z) & \text{ if }z\in \overline{\mathcal{E}_1}.
		\end{cases}$$
		It follows from this construction that $\tilde{\phi}_w^*\,\Lambda = \tilde{\lambda}_w$.
		Thus
		$$\rho_w:=\tilde{\phi}_w\circ T_{\delta}\circ \rho\circ \phi_{w}^{-1}$$ is holomorphic and semi-conjugates $f^q$ to $T_1$. It follows from the uniqueness in part (\ref{Fatou uniqueness}) of Theorem \ref{flower parabolic} that there exists some 
		$\delta_w\in \mathbb{C}$ so that 
		$$\rho_w = T_{\delta_w}\circ \rho.$$
		Similarly, the function $${\chi}_w(z) := \begin{cases}
		\phi_w\circ \chi_+\circ \tilde{\phi}_w^{-1}(z)&\text{ if }z\in \overline{W_+}\\
		\chi_+(z) & \text{ if }z\in \overline{\mathcal{E}_1}\cap Dom(\chi_+)
		\end{cases}$$
		is holomorphic and there exists some $\delta_w'\in \mathbb{C}$ so that 
		$${\chi}_w = \chi_+\circ T_{\delta_w'}.$$
		As ${\chi}_w = \chi_+$ on $\mathcal{E}_1\cap \text{Dom}(\chi_+)$, we must have $\delta_w'=0$.
		Thus 
		$$\phi_w\circ L_{\delta} \circ \phi_w^{-1} = \phi_w\circ \chi_+\circ \tilde{\phi}_w^{-1}\circ \tilde{\phi}_w\circ T_{\delta}\circ \rho\circ \phi_w^{-1} = \chi_+\circ T_{\delta_w}\circ \rho = L_{\delta_w}$$
		on $\phi_w(K_1^{L_{\delta}}) = K_1^{L_{\delta_w}}.$
		Moreover, 
		\begin{align*}
		T_{\delta_w}\circ \rho\circ L_{\delta_w}^{d}(cv^f)
		&=\tilde{\phi}_w\circ T_{\delta}\circ \rho\circ \phi_w^{-1}\circ L_{\delta_w}^{d}(cv^f)\\
		&= \xi_w\circ T_{\delta}\circ \rho\circ L_{\delta}^{d}\circ \phi_w^{-1}(cv^f)\\
		&= \xi_w(\delta)\\
		& = w.
		\end{align*}
		As $\phi_w$ depends continuously on $w$, the mapping $\Psi(w):= \delta_w$ is continuous.
		Moreover, $\Psi$ is biholomorphic as it is a continuous inverse  branch of the holomorphic function $$z\mapsto T_z\circ \rho\circ L_z^d(cv^f).$$
		It follows from this  construction that  $\Psi(\mathbb{H})\subset \mathcal{E}_{d+1}$, hence $\Psi(\mathbb{H})\subset \mathcal{B}.$ As any point on the boundary of the image of $\Psi$ cannot belong to $\mathcal{E}_{d+1}$,  $\Psi$ is in fact an analytic isomorphism from $\mathbb{H}$ to $\mathcal{B}.$
	\end{proof}
	
	For any basic enriched angle $\Theta$ of depth larger than $1$, we denote by $\mathcal{Z}_\Theta^f = \mathcal{Z}_\Theta$ the unique parameter in $\partial \mathcal{B}_\Theta$ such that $cv^f= z_\Theta^{L_{\mathcal{Z}_\Theta}}$; we call $\mathcal{Z}_\Theta$ the \textit{parent} of $\mathcal{B}_\Theta.$

	\section{Limits of External rays}\label{rays section}

	For any $N\geq 0$ and $\sigma= (k_n)_{n=1}^N\in \mathbb{Z}^N$ we denote
	\begin{align*}
		\ell(\sigma):=N\text{ and }
		|\sigma|:= \inf_{1\leq n\leq N} |k_n|,
	\end{align*}
	setting $|\sigma|:=\infty$ when $N=0.$ 
	
	For any map $f\in \mathcal{F}^\flower$ with a $p/q$-parabolic $k$-periodic point at zero and any $\Sigma\subset \bigcup_{N\geq 0}\mathbb{Z}^N$, 
	we define a  \textit{Log-multiplier family near $f$  indexed by $\Sigma$} to be a collection of holomorphic maps $\{h_{\sigma, \alpha}\}=\{h_{\sigma, \alpha}\}_{\sigma\in \Sigma, \alpha\in \mathbb{D}}$ 
	such that:
	\begin{enumerate}
		\item $h_{\sigma, \alpha}$ depends continuously on $\sigma$ and holomorphically on $\alpha\in \mathbb{D}$.
		\item $h_{\sigma, \alpha}$ has a $k$-periodic point at zero with multiplier $e^{2\pi i\mu_{p/q}(\alpha)}$.
		\item $h_{\sigma, \alpha}\to f$ when $|\sigma|\to \infty$ and $\ell(\sigma)\to \sup_{\sigma'\in \Sigma}\ell(\sigma').$
	\end{enumerate}

	Let us now fix some $f\in \mathcal{F}^\flower$ with a $p/q$-parabolic fixed point at zero with $\mathfrak{S}(p/q) = +1$ and a Log-multiplier family $\{h_{\sigma, \alpha}\}$ indexed by some $\Sigma$ with $\alpha\in \mathbb{D}$; the cases where zero is a  periodic point of $f$ or $\mathfrak{S}(p/q) = -1$ can be studied similarly.  To simplify our notation we will suppress the dependence on $\sigma$ and $\alpha$ and  write $h = h_{\sigma, \alpha}.$
	
	We will  assume that for all $\mathfrak{C}>1$, if  $\alpha\in A(\mathfrak{C})$ is  sufficiently close to $0$ and if  $h$ is sufficiently close to $f$, then there exists an analytic function
	$$\varphi_{-1, att}^{h, f}:  Dom(\varphi_{-1, att}^{f})\cap T_{1/\alpha}(Dom(\varphi_{-1, rep}^f))\to Dom(h^q)$$ which satisfies:
	\begin{enumerate}
		\item $h^q\circ\varphi_{-1, att}^{h, f}= \varphi_{att}^{h,f}\circ T_1$ on $Dom(\varphi_{-1, att}^{h, f})$.
		\item $\varphi_{-1, att}^{h, f}$ depends continuously and holomorphically on $h$ and converges to $\varphi_{-1, att}^{f}$ when $h\to f$.
		\item $\varphi_{-1, rep}^{h,f}:= T_{-1/\alpha}\circ \varphi_{-1, att}^{h,f}$ converges to $\varphi_{-1, rep}^{f}$  when $h\to f$.
	\end{enumerate}
	We will say that the Log-multiplier family has \textit{pre-petals} in this case. Note that this condition holds automatically if the immediate basin $U_0^f$ is compactly contained in the domain of $f^q$, for example if $f$ is a polynomial.
	
	Let us now fix some Lavaurs map $L = L_\delta^f$ for $f$ which is $d$-nonescaping for some $d\geq 0$. 
	Note that we can choose $\mathfrak{C}>1$ sufficiently large so that if $n-\frac{1}{\alpha}$ is close to $\delta$ for some large integer $n$ then $\alpha\in A(\mathfrak{C})$.
	Using the approximations of $\varphi_{-1, \iota}^L$, we can construct approximations for all of the pre-petal parameters for $L$.
	\begin{proposition}\label{perturbed pre-petals}
		Fix some enriched angle $\Theta$ of depth $d$ with principle integer $j$. If $h$ is sufficiently close to $f$ and if $n-1/\alpha$ is sufficiently close to $\delta$ for some integer $n\geq 0$, then there exists an analytic function $$\varphi_{\Theta, att}^{h, L}:Dom(\varphi_{\Theta, att}^{f})\cap T_{1/\alpha}(Dom(\varphi_{\Theta, rep}^f))\to Dom(h)$$
		 such that:
		\begin{enumerate}
			\item $h^{d(nq+k_+^f)+jq}\circ \varphi_{\Theta, att}^{h,f} = \varphi_{att}^{h,f}\circ T_j$   on the domain of $\varphi_{\Theta, att}^{h, f}$.
			\item $\varphi_{\Theta, att}^{h, f}$ depends continuously and holomorphically on $h$ and converges locally uniformly to $\varphi_{\Theta, att}^{f}$ when $h\to f$ and $n-1/\alpha\to \delta$.
			\item $\varphi_{\Theta, rep}^{h,f}:= T_{-1/\alpha}\circ \varphi_{\Theta, att}^{h,f}$ converges to $\varphi_{\Theta, rep}^{f}$ when $h\to f$ and $n-1/\alpha\to \delta$.
		\end{enumerate}
	\end{proposition}
	
	\begin{proof}
		For $\Theta\sim 0$, we can define $\varphi_{\Theta, \iota}^{h, L}:= \varphi_{\iota}^{h,f}$. If instead $j>0$, 
		then it follows from the proof of Proposition \ref{pre-petals} that $L^{d}\circ f^{(j-1)q}$ univalently maps a neighborhood of $z_{\Theta}^L$ containing the images of $\varphi_{\Theta, att}^L$ and $\varphi_{\Theta, rep}^L$ to a neighborhood of $z_{-1}^L$ containing the images of $\varphi_{-1, att}^L$ and $\varphi_{-1, rep}^L$. As $h^{nq+k_+^f}$ converges to $L$ when $h\to f$ and $n-1/\alpha\to \delta$, we can use the corresponding branches of $(h^{d(nq+k_+^f)+(j-1)q})^{-1}$ to define
		$$\varphi_{\Theta, \iota}^{h, L}:= (h^{d(nq+k_+^f)+(j-1)q})^{-1}\circ \varphi_{-1, \iota}^{h, f}\circ T_{j-1}.$$
	\end{proof}

	As a consequence of the above  and Proposition \ref{angle dynamics}, we have the following relationship between our labeling of the pre-petal parameters and the dynamics of $h$:
	
		\begin{corollary}
			For any enriched angle $\Theta$ of depth $d$, if $h$ is sufficiently close to $f$ and $n-1/\alpha$ is close to $\delta$ then 
			\begin{align*}h^q\circ \varphi_{\Theta, \iota}^{h, L}(w)&= \varphi_{\Theta+1, \iota}^{h, L}\circ T_1(w) \text{ and }\\
			 h^{nq+k_+^f}\circ \varphi_{\Theta, \iota}^{h, f}(w)&= \varphi_{\lceil\Theta\rceil}^{h, L}(w)
			\end{align*}
			for any $w$ where both sides of the respective equations are defined.
		\end{corollary}

	For any topological space $S$, sets $Y_s\subset \hat{\mathbb{C}}$ parameterized by $s\in S$, and any compact set $X\subset \hat{\mathbb{C}}$, we will say that $Y_s \xrightarrow{\subset} X$ when $s\to s_0\in \Lambda$ if 
	$$\sup_{y\in \overline{Y_s}}\inf_{x\in X}d_{\hat{\mathbb{C}}}(x,y)\to 0$$
	when $s\to s_0$, where $d_{\hat{\mathbb{C}}}$ is the spherical metric.
	
	Let $x_0^\sigma>0$ and $\lambda^\sigma\geq 2$ be real numbers which depend only on $\sigma$ and set $$X_m^\sigma:= \{(\lambda^\sigma)^m(x-iy): |x|< x_0^\sigma,|\log y|< \log\lambda^\sigma\}$$
	for all integers $m$. We will now assume that if $h$ is sufficiently close to $f$, then there exists an analytic function $\psi_0^h$ defined on  of $X_0^\sigma$ such that:
	\begin{enumerate}
		\item $\psi_0^h$ depends continuously and holomorphically $h$.
		\item 
		$$h^q\circ\psi_{0}^h (w) = \psi_{0}^h(\lambda^\sigma w)$$
		for all $w$ where both sides of the equation are defined. 
		\item There exists a compact set $\hat{X}_{0, 0}^f\subset P_{rep}^f$ which avoids $U_f$ such that $$\psi_{0}^h(X_0^\sigma)\xrightarrow{\subset}  \hat{X}_{0, 0}^f$$ when $h\to f$.
	\end{enumerate}
	We will say that the log-multiplier family has \textit{external rays} in this case, 
	we will call  $\psi_0^h$ the  \textit{external coordinate} of $h$, and we will call any image of a vertical line under $\psi_0^h$ an \textit{external ray} for $h$. To simplify our notation, we will suppress the dependence on $\sigma$ and write $\lambda = \lambda^\sigma$, $x_0= x_0^\sigma$, and  $X_m = X_m^\sigma$. We will also write $\psi_\theta^{h, L}= \psi_0^h$ for any $\theta\in \mathcal{A}^f$ equivalent to $0.$

	\begin{proposition}\label{ray limits}
		Fix some enriched angle $\Theta$ of depth $d$. There exists an integer $M'$ such that for any $M \leq  M'$, if $h$ and $n-1/\alpha$ are sufficiently close to $f$ and  $\delta$ respectively then there exists an analytic function $\psi_{\Theta}^{h, L}$ defined on $\bigcup_{m= M}^{M'}X_m$ such that
		\begin{align*}
		h^{q}\circ \psi_\Theta^{h, L} (w) &= \psi_{\Theta+1}^{h, L}(\lambda w) \text{ and}\\
		h^{nq+k_+^f}\circ \psi_{\Theta}^{h, L}(w)&= \psi_{\lceil\Theta\rceil}^{h, L}(w)
		\end{align*}
		for any $w$ where both sides of the respective equations are defined. 
		For any $M\leq m\leq M'$ there exists a compact set $\hat{X}_{\Theta, m}^L\subset P_{\Theta, rep}^L$ such that 
		$\psi_\Theta^{h, L}(X_m)\xrightarrow{\subset} \hat{X}_{\Theta, m}^L$ when $h\to f$ and $n-1/\alpha\to \delta$.
	\end{proposition}
	\begin{proof}
		Let us denote by $g^h$ the continuous branch of $(\varphi_{rep}^{h,f})^{-1}$ defined on a neighborhood of $\hat{X}_{0, 0}$ so $$g^h\circ \psi_0^{h, f}(X_0)\xrightarrow{\subset} (\varphi_{rep}^f)^{-1}(\hat{X}_{0, 0})$$ when $h\to f$. Let $j$ be the principle integer of $\Theta$. Proposition \ref{pre-petals} implies that if $M'$ is a negative integer with large magnitude, then
		$$T_{M'}(\varphi_{rep}^f)^{-1}(\hat{X}_{0, 0}^f)\subset Dom(\varphi_{\Theta, rep}^L).$$
		Thus for any $M< M'$, if $h$ is  close to $f$ and $n-\frac{1}{\alpha}$ is close to $\delta$ then we can define
		$$\psi_{\Theta}^{h, L}(w):= \varphi_{\Theta, rep}^{h, L}\circ T_m\circ g^h\circ \psi_0^{h}(\lambda^{-m} w)$$
		for all $w\in X_m$ and $M\leq m\leq M'$. The desired properties all follow immediately.`
	\end{proof}

	Note that in the functions $\psi_\theta^{h, L}$ do not depend on $L$ when for any $\theta\in \mathcal{A}^f$; in this case we will denote $\psi_\theta^h:= \psi_\theta^{h, L}.$ It follows from Proposition \ref{pre-petals} and the proof of Proposition \ref{ray limits} that if $j$ is principle integer of $\theta$, then $\psi_\theta^h$ is defined on $X_{-j}.$

	We will now assume that for any sufficiently large  negative $\theta\in \mathcal{A}^f$ with principle integer $j>0$ there exists some $b_\theta= b_\theta^\sigma$ such that if $h$ is sufficiently close to $f$ then:
	\begin{enumerate}
		\item $\psi_{\theta}^h$ can be analytically extended to $\bigcup_{m=0}^{j}X_{-m}$.
		\item For any $w$ in the non-empty set $X_0\cap T_{-b_\theta}(X_0)$,
		$$\psi_{\theta}^h(w)= h^{k_+^f}\circ \psi_0^h\circ T_{b_\theta}(w).$$
		\item For all $0 \leq m \leq j$, there exists a compact set $\hat{X}_{\theta, -m}\subset f^{k_+^f}(P_{rep}^f)$ which avoids $U^f$ such that $\psi_{\theta}^h(X_{-m})\xrightarrow{\subset} \hat{X}_{\theta, -m}.$
	\end{enumerate}
	In this case we will say that the external rays and pre-petals for the log-multiplier family are \textit{compatible}. 
	\begin{proposition}\label{b theta inductively}
		If $b_\theta$ and $b_{\theta+1}$ are defined for some $\theta\in \mathcal{A}^f$, then $\lambda b_\theta = b_{\theta+1}.$
	\end{proposition}
	\begin{proof}
		Let $j$ be the principle integer of $\theta$, so $j-1$ is the principle integer of $\theta+1$. 
		For  $w\in X_{-j}$, we have
		$$h^{k_+^f}\circ \psi_0^h(\lambda w +\lambda b_\theta) = h^q\circ \psi_\theta^h(w) = \psi_{\theta+1}^h(\lambda w) = h^{k_+^f}\circ \psi_0^h(\lambda w + b_{\theta+1}).$$
	\end{proof}
	Using Proposition \ref{b theta inductively} we inductively define $b_{\theta+1} :=\lambda b_\theta $ for all $\theta\in \mathcal{A}^f$; to simplify our notation we will write $\xi_{\theta} = T_{b_\theta}.$ Using compatibility we can extend Proposition \ref{ray limits}.

	\begin{proposition}\label{ray limits extending}
		If $d>0$, then for any enriched angle $\Theta$ of depth $d$ and  any integers $M\leq M'$, if $h$ and  $n-1/\alpha$ are sufficiently close to $f$ and $\delta$ respectively for some integer $n\geq 0$, then $\psi_{\Theta}^{h, L}$ can be analytically extended to $\bigcup_{m=M}^{M'}(X_m)$. For any $M\leq m\leq M'$ there exists a compact set $\hat{X}^L_{\Theta, m}\subset B_{\lfloor\Theta\rfloor}\setminus \{z_{\lfloor\Theta\rfloor}\}$ such that $\psi_\Theta^{h, L}(X_m)\xrightarrow{\subset} \hat{X}_{\Theta, m}^L$ when $h\to f$ and $n-1/\alpha\to \delta$.
	\end{proposition}

	\begin{proof}
		Let us first assume that $d=1$ and fix some $\theta\in \mathcal{A}_f$ and integers $M\leq M'$.
		Thus there is some integer $s$ so that $\psi_{\langle 0, \theta\rangle}^{h, L}$ is defined on $X_s$ when  $h$ and $n-1/\alpha$ are sufficiently close to $f$ and $\delta$ respectively. Let $M'\geq s$ be a sufficiently large positive integer so that $\hat{X}_{\theta-M', -m}\subset f^{k_+^f}(P_{rep}^f)$ for all $0\leq m \leq M'-s.$ Let $\ell\geq 0$ be a sufficiently large integer so that 
		$$T_{\ell-\delta}\circ (\varphi_{rep}^f)^{-1}\circ f^{-k_+^f}(\hat{X}_{\theta-M', -m})\subset \Omega_{att}^f$$
		for all $0 \leq m \leq M'-s$. 
		It follow Proposition \ref{angle dynamics invertible} and the proof of Proposition \ref{ray limits} that if $h$ and $n-1/\alpha$ are sufficiently close to $f$ and $\delta$ respectively then 
		$$\psi_{\langle 0, \theta+\ell - M'\rangle}^{h, L}(\lambda^\ell w)= \varphi_{att}^{h, f}\circ T_{\ell-n}\circ (\varphi_{att}^{h, f})^{-1}\circ h^{-k_+^f}\circ \psi_{\theta-M'}(w)$$
		for all $w\in X_{-m}$ and $0\leq m\leq M'-s.$
		Applying $h^{(M'-\ell)q}$, it follows that we can analytically extend $\varphi_{\langle 0, \theta\rangle}^{h, L}$ to $X_m$ for all $s\leq m \leq M'$. 
		As $P_{att}^f\cap E_1$ has a unique connected component, the corresponding sets 
		$\hat{X}_{\langle 0, \theta\rangle, m}^L$ for $s\leq m \leq M'$ must belong to $B_0^L$. 
		This completes the proof for basic enriched angles when $d=1$. 
		As $L$ is $d$-nonescaping, $B_\Theta$ is mapped univalently onto $B_{\langle 0, \theta\rangle}$ by $L^{d-1}\circ f^{mq}$ for some $(d-1, m)\geq (0, 0)$ and $\theta\in \mathcal{A}^f$. This allows us to extend the above to the general case.
	\end{proof}

	\begin{proposition}\label{ray limits in petal}
		Fix some enriched angle  $\Theta$ of depth $d$, some $\theta\in \mathcal{A}^f$, and some neighborhood $D$ of $z_\Theta^L$.  There exists an integer $M$ such that if $h$ and $n-1/\alpha$ are sufficiently close to $f$ and $\delta$ respectively, then $\psi_{\Theta}^{h, L}$ can be analytically extended so that it maps
		$$\xi_{\theta-n}\left(\bigcup_{m=M}^{n-M}X_{-m}\right)$$
		into $D$.  If $L$ is $(d+1)$-nonescaping, then $$\psi_{\Theta\oplus \theta}^{h, L}(w) = \psi_{\Theta}^{h, L}( \lambda^{-n}\xi_{\theta}(w))$$
		for all $w$ where both sides of the equation are defined.
	\end{proposition}

		\begin{proof}
			As $b_{\theta-m}\to 0$ when $m\to \infty$, there exists some $M\geq 0$ such that 
		$$X:= \bigcup_{m\geq M}\xi_{\theta-m}(X_{0})\subset X_0\cup \xi_{\theta-M}(X_0).$$
		It therefore follows from the definition of compatible external rays and pre-petals that, after increasing $M$ if necessary,  $\psi_0^h$ is defined on $X$
		when $h$ is close to $f$ and 
		there is a compact set $\hat{X}\subset \varphi_{rep}^f(\Omega_{rep}^f(0))$ such that $\psi_0^h(X)\to \hat{X}$ when $h\to f$.
		For some $R>0$, let $M'\geq 0$ be an integer sufficiently large so that $$T_{M'-\delta}\circ (\varphi_{rep}^f)^{-1}(\hat{X})\subset T_R(\Omega_{att}^{f}).$$
		Note that if we fix $M$ then the above inclusion still holds when we increase $M'$, and conversely if we fix $M'$ then the above inclusion still holds when we increase $M$. Hence we may take $M'= M$. 
		As $\varphi_{rep}^f$ is univalent, when $h$ is sufficiently close to $f$ there exists a unique continuous branch of $(\varphi_{rep}^{h, f})^{-1}$ defined on a neighborhood of $X$ which converges to $(\varphi_{rep}^f)^{-1}$ when $h\to f$. It follows that
		$$T_{M-n}\circ (\varphi_{att}^{h, f})^{-1}\circ \psi_0^h(X)=T_{M-n+1/\alpha}\circ (\varphi_{rep}^{h, f})^{-1}\circ \psi_0^h(X)\subset T_R(\Omega_{att}^{h, f})$$
		when $h$ is close to $f$ and $n-1/\alpha$ is close to $\delta$, using the continuous branch of $(\varphi_{rep}^{h, f})^{-1}$ which converges to $(\varphi_{rep}^{f})^{-1}$ locally uniformly on a neighborhood of $\hat{X}$ when $h\to f$.
		Thus we can define 
		$$\psi_0^{h, L}(\lambda^{-m}w):= \varphi_{att}^{h, f}\circ T_{-m}\circ(\varphi_{att}^{h, f})^{-1}\circ \psi_0^h(w)$$
		for all $0\leq m\leq n-M$ and $w$ in a neighborhood of $X$. In particular, this gives us an analytic extension of $\psi_0^{h, L}$ to 
		$$\xi_{\theta-n}\left(\bigcup_{m=0}^{n-M}X_{-m}\right).$$
		As $\varphi_{\iota}^{f}(w)\to 0$ when $w\to \infty$, we can choose $R$ sufficiently large so that $$\overline{\varphi_{att}^{ f}\circ T_R(\Omega_{att}^f)\cup \varphi_{att}^{f}\circ T_{-R}(\Omega_{rep}^f)}\subset D.$$
		Increasing $M$ if necessary, it follows that 
		$$\psi_{0}^{h, L}\circ \xi_{\theta-n}\left(\bigcup_{m=M}^{n-M}X_{-m}\right)\subset \overline{\varphi_{att}^{h, f}\circ T_R(\Omega_{att}^{h, f})\cup \varphi_{att}^{h, f}\circ T_{-R}(\Omega_{rep}^{h,f})}\subset D.$$
		
		 If $L$ is $1$-nonescaping, then it follows from the proof of proposition \ref{ray limits extending} that $$\psi_{\theta-m}^{h, L}(\lambda^{-m}w) =h^{k_+^f}\circ\varphi_{att}^{h, f}\circ T_{n-m}\circ(\varphi_{att}^{h,f})^{-1}\circ \psi_{0\oplus \theta}^{h, L}(w)$$
		for any $w$ where both sided of the equation are defined. 
		As 
		$$h^{k_+^f}\circ\varphi_{att}^{h, f}\circ T_{n-m}\circ(\varphi_{att}^{h,f})^{-1}\circ \psi_{0}^{h, L}( \xi_{\theta-n}(w))= h^{k_+^f}\circ \psi_0^{h, L}(\lambda^{n-m} \xi_{\theta-n}(w))= \psi_{\theta-m}^{h, L}(\lambda^{n-m}w)$$
		by the definition of compatible external rays and pre-petal coordinates, we can conclude that 
		$$\psi_{0\oplus \theta}^{h, L}(w) =\psi_0^{h, L}\circ \xi_{\theta-n}(\lambda^{-n}w)= \psi_0^{h, L}(\lambda^{-n}\xi_{\theta}(w))$$
		for any $w$ where both sided of the equation are defined.
		
		The above proves the proposition for $\Theta\sim 0$, for general choices of $\Theta$ the proposition follows by considering further pre-images under $h^q$ and $h^{nq+k_+^f}.$
	\end{proof}

	For any $\theta \in \mathcal{A}^f$, we define $\psi_{\theta, n}^h:= \psi_\theta^h$ for all integers $n$. If $d>0$, then for any  enriched angle $\Theta= \langle \theta_m\rangle_{m=0}^d$ of depth $d$ and any integer $n$, we define the function $$\psi_{\Theta, n}^{h, f}(w):= \psi_{\lfloor\Theta\rfloor, n}^{h, f}\left(\lambda^{-n} \xi_{\theta_d}(w)\right).$$
	Proposition \ref{ray limits in petal} implies that $\psi_{\Theta, n}^{h, f} = \psi_{\Theta}^{h, L}$ wherever both sets are defined when $h$ and $n-1/\alpha$ are  close to $f$ and  $\delta$ respectively.

	\subsection{Parameter rays}
	
	Now we will lift the external rays in the previous subsection to the parameter space of the Log-multiplier family. Our main tool in this task is Lemma \ref{walz} below, which follows directly from the argument principle. Given a simple closed curve $\gamma\subset \mathbb{C}$, sets $Y_t\subset \mathbb{C}$ which depend on $t\in\gamma$, and a continuous function $g_1: \gamma\to \mathbb{C}$,  we will say that $g_1(t)$ \textit{winds once around $Y_t$ when $t$ winds once around $\gamma$} if 
	$$\frac{1}{2\pi i}\oint_{\gamma}\frac{dt}{g_1(t)-g_2(t)} = 1$$
	for any continuous function $g_2:\gamma\to \mathbb{C}$ satisfying $g_2(t)\in Y_t$.
	\begin{lemma}\label{walz}
		Fix some Jordan domain  $D$ and some $X\subset \mathbb{C}$. Let $g_1:\overline{D}\to \mathbb{C}$ be a continuous function which is holomorphic on ${D}$ and let $g_1: X\times \overline{{D}}\to \mathbb{C}$ be a map such that $t\mapsto g_1(z, t)$ is continuous on $\overline{D}$ and holomorphic on $D$ for all $z\in X$. If $g_1(t)$ winds once around $g_2(X,t)$ when $t$ winds once around $\partial D$,  then there exists a function $\phi: X\to D$ such that $\phi(z)$ is the unique point in $D$ satisfying 
		$$g_1(\phi(z))=g_2(z, \phi(z)).$$
		If the function $z\mapsto g_2(z, t)$ is injective, continuous, or holomorphic for all $t\in {D}$, then $\phi$ is injective, continuous, or holomorphic respectively.
	\end{lemma}

	Using Lemma \ref{walz}, we can construct analogues of the external coordinates in parameter space.
	\begin{proposition}\label{parameter rays}
		For any sufficiently small neighborhood $\mathcal{D}$ of $(\varphi_{rep}^f)^{-1}(\hat{X}_{0, 0}^f)$ there exists some $n_0\geq 0$ such that if $|\sigma|$ and $\ell(\sigma)$ are sufficiently large then there exists a holomorphic function 
		$$\Psi^\sigma: \bigcup_{n\geq n_0}X_{-n}^\sigma\to \mathbb{D}$$
		such that for any $w\in X_0$ and $n\geq n_0$, $\alpha= \Psi^\sigma(\lambda^{-n}w)$ is the unique choice of $\alpha\in \mathbb{D}$ satisfying $n-\frac{1}{\alpha}\in \mathcal{D}$ and
		\begin{equation*}\label{parameter ray equation}
			\varphi_{rep}^{h, f}\left(n-\frac{1}{\alpha}\right) =h^{nq}(cv^{h, f}) = \psi_0^{h, f}(\lambda^n w).
		\end{equation*}
		For any large $R>0$, if both $M\geq 0$ and $n\geq 0$ are sufficiently large then for any $w\in X_{-M}$, $\alpha = \Psi^\sigma(\lambda^{-n}w)$ is the unique choice of 
		$\alpha\in \mathbb{D}$ satisfying (\ref{parameter ray equation}) and
			$$\left|n+M-\frac{1}{\alpha}\right|< R.$$ 
	\end{proposition}

	\begin{proof}
		Let $\mathcal{D}\subset \Omega_{rep}^f$ be a Jordan domain containing 
	 	$(\varphi_{rep}^f)^{-1}(\hat{X}_{0, 0}^f)$. 
	 	Thus 
	 	$\varphi_{rep}^f(\delta)$ winds once around $\hat{X}_{0, 0}^f$ when $\delta$ winds once around the $\partial \mathcal{D}$. For all $n\geq 0$ we set 
		$$\mathcal{D}_n:= \left\{\alpha\in \mathbb{C}: n- \frac{1}{\alpha}\in \mathcal{D}\right\}.$$
		As $\varphi_{rep}^{h, f}\to \varphi_{rep}^f$ and $\psi_{0}^h(X_{0})\xrightarrow{\subset}  \hat{X}_{0, 0}$ when $h\to f$, if $|\sigma|$, $\ell(\sigma)$, and $n\geq 0$ are sufficiently large then 
		$$\varphi_{rep}^{h, f}\left(n-\frac{1}{\alpha}\right)= h^{nq}(cv^{h, f})$$ winds once around $\psi_{0}^{h,f}(X_{0}^\sigma)$ when $\alpha$ winds once around $\partial \mathcal{D}_n.$
		As the function $$(\alpha, w)\mapsto \psi_0^{h, f}(w)$$ is holomorphic on $\mathcal{D}_n\times X_{0}$, 
		it follows from Lemma \ref{walz}   that there is some $n_0\geq 0$ and a holomorphic function ${\Psi}_n^\sigma: X_{0}\to \mathcal{D}_n$ for all $n\geq n_0$ such that 
		$\alpha={\Psi}_n^\sigma(w)$ is the unique choice of $\alpha\in\mathcal{D}_n$ satisfying 
		$$\psi_0^{h, f}(w)= h^{nq}(cv^{h, f})$$
		for any $w\in X_{0}.$
		Setting 
		$\Psi^\sigma(\lambda^{-n}w) = \Psi_n^\sigma(w)$ for all $n\geq n_0$ and $w\in X_{0}$, it follows from the uniqueness above that $\Psi^\sigma$ is a holomorphic function on $\bigcup_{n\geq n_0}X_{-n}$.
		
		Let $R>0$ be large enough so that $(\varphi_{rep}^f)^{-1}(\hat{X}_{0, 0}^f)\subset \mathbb{D}_R$, and let $M\geq 0$ be large enough so that 
		$$(\varphi_{rep}^f)^{-1}(\hat{X}_{0, -M}^f)\subset T_{-M}(\mathbb{D}_R)\subset \Omega_{rep}^f.$$
		Thus 
		$\varphi_{rep}^f(\delta)$ winds once around $\hat{X}_{0, -M}^f$ when $\delta$ winds once around the boundary of $T_{-M}(\mathbb{D}_R)$. Setting 
		$$\mathcal{D}_n:= \left\{\alpha\in \mathbb{C}: n+M- \frac{1}{\alpha}\in \mathbb{D}_R\right\}$$
		for all $n\geq 0$, we can repeat the above argument to complete the proof.
	\end{proof}

	We will call the map $\Psi^\sigma$ an \textit{external parameter} for the Log-multiplier family. For any $t\in \mathbb{R}$ we will call
	$$\Psi^\sigma\left(\{t-iy: y>0\}\cap Dom(\Psi^\sigma)\right)$$
	a \textit{parameter ray at angle $t$}; we will say that the parameter ray \textit{lands} at some $\alpha\in \mathbb{D}$ when $\Psi^\sigma(t-iy)$ is defined for all small $y>0$ and converges to $\alpha$ when $t\to 0$. We note that it follows from Proposition \ref{parameter rays} that the parameter ray at angle $0$ lands at $0$ when $|\sigma|$ and $\ell(\sigma)$ are large.  Our goal for the remainder of this section is to analytically extend the external parameters $\Psi^\sigma$ and control how some of the parameter rays land.
	
	For any basic enriched angle $\Theta= \langle \theta_n\rangle_{n=0}^d$ of depth $d\geq 1$,  integer $M$, and large integer $j\geq 0$, we define 
	\begin{align*}
		\hat{\mathcal{X}}_{\Theta, M}^f:=&\{\delta \in \mathcal{K}_{d}^f:  L_\delta^f\circ f^{-jq}(cv^f)\in \hat{X}_{\lceil\Theta\rceil-j, M-j}^{L_\delta^f} \}\\
		 =& \{\delta \in \mathcal{K}_{d}^f:  f^{k_+^f}\circ \varphi_{rep}^{f}\circ T_{-j}(\delta) \in \hat{X}^{L_\delta^f}_{\lceil\Theta\rceil-j, M-j} \}\\
		 =& \Upsilon_{\lfloor\Theta\rfloor}\circ T_j\circ \Upsilon_0^{-1} \circ (\varphi_{rep}^f)^{-1}\circ f^{-k_+^f}(\hat{X}_{\theta_d-j, M-j}^{f} );
	\end{align*}
	 it follows from Proposition \ref{ray limits} that for any $\Theta$ and $M$ this set does not depend on $j$ when $j$ is sufficiently large.
	Applying Lemma \ref{walz}, we have the following parameter analogues of Propositions \ref{ray limits} and \ref{ray limits in petal}.

	\begin{proposition}\label{ray limits parameter}
		For any basic enriched angle $\Theta$ of depth $d\geq 1$,  integer $M$, and any sufficiently large integer $j\geq 0$, there exists an integer $n_0\geq 0$ such that if $|\sigma|$ and $\ell(\sigma)$ are sufficiently large then there exists a a neighborhood $\mathcal{D}$ of $\hat{\mathcal{X}}_{\Theta, M}^f$ and a  holomorphic map 
		$$\Psi_{\Theta, n}^\sigma: X_M\to \mathbb{D}$$
		for all $n\geq n_0$ such   that for any $w\in X_M$,  $\alpha= \Psi_{\Theta, n}^\sigma(w)$ is the unique choice of $\alpha$ with $n-\frac{1}{\alpha}\in \mathcal{D}$ satisfying
			$$h^{(n-j)q+k_+^f}(cv^{h, f}) = \psi_{\lceil \Theta\rceil -j, n}^{h, f}(\lambda^{-j}w).$$
	\end{proposition}

	\begin{proof}
	For $\Theta = \langle \theta_m\rangle_{m=0}^d$ and a large integer $j\geq 0$ we set $$Y:= T_j\circ \Upsilon_1\circ (\varphi_{rep}^f)^{-1}\circ f^{-k_+^f}(\hat{X}_{\theta_d-j, M-j}^{L_\delta^f} ),$$
	so $\Upsilon_{d}(\mathcal{X}_{\Theta,M}) = Y$ and 
		$$T_j\circ \upsilon_{d-1}^{L_\delta}(\hat{X}_{\lceil\Theta\rceil-j, M-j}^{L_\delta}) = Y.$$
		Let us suppose for now that $Y\subset \mathbb{H}$, and let $\mathcal{D}\subset \mathcal{B}_{\lfloor \Theta\rfloor}$ be a Jordan domain so that $Y\subset \Upsilon_{d}(\mathcal{D})$. Thus $\Upsilon_{d}(\delta) $ winds once around $Y$ when $\delta$ winds once around $\partial \mathcal{D}$.
		As $\upsilon_{d-1}^{L_\delta}$ is injective and does not depend on $\delta$ when $d=1$, it follows that in this case $L_\delta \circ f^{-jq}(cv^f)$ winds once around  $\hat{X}_{\lceil\Theta\rceil-j, M-j}^{f}$ when $\delta$ winds once around $\partial \mathcal{D}$. When $d>1$, the fact that $\upsilon_{d-1}$ is injective on $B_{\lceil \lfloor\Theta\rfloor\rceil}^{L_\delta}$ allows us to reach the same conclusion; we can inductively use the fact that $z_\theta^{L_\delta}$ does not depend on $\delta$ for any $\theta\in \mathcal{A}^f$ to conclude that composing by $(T_j\circ \upsilon_{d-1}^{L_\delta})^{-1}$ does not affect the winding number.

		For all $n\geq 0$, set 
		$$\mathcal{D}_n := \left\{\frac{1}{n-\delta}:\delta\in \mathcal{D}\right\}.$$
		It  follows from the above that there exists some $n_0\geq 0$ such that if $|\sigma|$ and $\ell(\sigma)$ are sufficiently large and $n\geq n_0$ then
		$h^{(n-j)q+k_+^f}(cv^{h, f})$ winds once around $\psi_{\lceil\Theta\rceil-j, n}^{h, f}(X_{M-j})$
			when $\alpha$ winds once around $\partial \mathcal{D}_n$. As the functions $$\alpha\mapsto h^{(n-j)q+k_+^f}(cv^{h, f})\text{ and }(\alpha, w)\mapsto \psi_{\lceil\Theta\rceil-j, n}^{h, f}(\lambda^{-j}w)$$ are holomorphic on $\mathcal{D}_n$ and $\mathcal{D}_n\times X_M$ respectively for any $\sigma$ and $n$, it follows from Lemma \ref{walz} that when $|\sigma|$ and 
			$\ell(\sigma)$ are large and $n\geq n_0$ there exists a holomorphic  function 
			$\Psi_{\Theta, n}^\sigma: X_{M}^\sigma\to \mathcal{D}_n$ such that $\alpha = \Psi_{\Theta, n}^\sigma(w)$ is the unique choice of $\alpha\in \mathcal{D}_n$ satisfying
			$$ h^{(n-j)q+k_+^f}(cv^{h, f})=\psi_{\lceil\Theta\rceil-j, n}^{h, f}(\lambda^{-j}w).$$
			
			In the case where $Y\cap \partial {\mathbb{H}}$ is non-empty, we can instead take  $\mathcal{D}$ to be a Jordan domain which compactly contains $\hat{\mathcal{X}}_{ \Theta, M}$. If $\mathcal{D}$ is chosen in a sufficiently small neighborhood of
				 $\hat{\mathcal{X}}_{ \Theta, M}$, 				
				then we can still guarantee that $L_\delta \circ f^{-jq}(cv^f)$				
				winds once around $\hat{X}_{ \lceil\Theta\rceil-j, M-j}^{L_\delta}$ 				
				when $\delta$ winds once around $\partial \mathcal{D}$; indeed otherwise $\hat{\mathcal{X}}_{\Theta, M}$ would not be homeomorphic to $Y$. The rest of the argument proceeds identically.
	\end{proof}

	\begin{proposition}\label{parameter ray limits implosive}
		Let $\Theta$ be a basic enriched angle of depth $d\geq 1$  and fix some neighborhood $\mathcal{D}$ of $\mathcal{Z}_\Theta$. For any $\theta\in \mathcal{A}^f$ there exist integers $M$ and $n_0\geq 0$ such that if $|\sigma|$ and $\ell(\sigma)$ are sufficiently large then for all $n\geq n_0$,  $\Psi_{\Theta, n}^\sigma$ can be analytically extended so that it maps
		$$\xi_{\theta-n}\left(\bigcup_{m=M}^{n-M}X_{-m}\right)$$
		into
		$$ \left\{\frac{1}{n-\delta}: \delta\in \mathcal{D}\right\}.$$  
		Moreover, $$\Psi_{\Theta\oplus \theta, n}^\sigma (w) = \Psi_{\Theta, n}^\sigma(\lambda^{-n}\xi_{\theta}(w))$$
		for all any $w$ where both sides of the equation are defined.
	\end{proposition}

	\begin{proof}
		First we observe that for all sufficiently large integers $j\geq 0$, the holomorphic function 
		$$\delta\mapsto L_\delta\circ f^{-jq}(cv^f) - z_{\lceil\Theta\rceil-j}^{L_{\delta}}$$
		is defined all $\delta$ close to $\mathcal{Z}_\Theta$ and has a unique zero at $\mathcal{Z}_\Theta$. Moreover, this zero is simple as otherwise there would be more than one component of $\mathcal{E}_{d}^f$ attached at $\mathcal{Z}_\Theta$. Thus for any sufficiently small Jordan domain $\mathcal{D}$ which contains $\mathcal{Z}_\Theta$, $L_\delta\circ f^{-jq}(cv^f)$ winds once around $z_{\lceil\Theta\rceil-j}^{L_{\delta}}$ when $\delta$ winds once around $\partial \mathcal{D}$ by the argument principle. 
		For some $s>0$, let $D_\delta$ be the disk centered at $z_{\lceil\Theta\rceil-j}^{L_{\delta}}$ of radius $s$. By the compactness of $\partial {\mathcal{D}}$,  we can choose $s$ small enough so that  $L_\delta\circ f^{-jq}(cv^f)$ winds once around $D_\delta$ when $\delta$ winds once around $\partial \mathcal{D}.$
		
		For all $n\geq 0$, we set 
		$$\mathcal{D}_n:= \left\{\frac{1}{n-\delta}: \delta\in \mathcal{D}\right\}.$$
		Thus there exists some $n_0\geq 0$ such that
		$h^{(n-j)q+k_+^f}(cv^f)$ winds once around $D_{n-1/\alpha}$ when $\alpha$ winds once around $\partial \mathcal{D}_n$ for all $n\geq n_0$ when $|\sigma|$ and $\ell(\sigma)$ are large.
		By the compactness of $\overline{\mathcal{D}}$ and Proposition \ref{ray limits in petal}, there exists some integer $M$ so that 
		$$\psi_{\lceil\Theta\rceil, n}^{h, f}\circ \xi_{\theta-n}\left(\bigcup_{m=M}^{n-M}X_{-m}\right)\subset D_{n-1/\alpha}$$
		for all $\alpha \in \overline{\mathcal{D}_n}$ when $|\sigma|$ and $n$ are large. As the functions
		$$\alpha \mapsto h^{(n-j)q+k_+^f}(cv^f) \text{ and }(\alpha, w)\mapsto \psi_{\lceil\Theta\rceil-j, n}^{h, f}(w)$$
		are holomorphic on $\mathcal{D}_n$ and $\mathcal{D}_n\times \xi_{\theta_n}\left(\bigcup_{m=M}^{n-M}X_{-m}\right)$ respectively for any $\sigma$ and $n$, it follows from Lemma \ref{walz} that for all $n\geq n_0$ there exists a holomorphic function $$\Psi: \xi_{\theta-n}\left(\bigcup_{m=M}^{n-M}X_{-m}\right)\to \mathcal{D}_n$$ when $|\sigma|$ and $\ell(\sigma)$ are large such that if $\alpha = \Psi(w)$ for some $w\in \xi_{\theta-n}\left(\bigcup_{m=M}^{n-M}X_{-m}\right)$, then 
		$$h^{(n-j)q+k_+^f}(cv^f) =\psi_{\lceil\Theta\rceil-j, n}^{h, f}(w).$$
		It follows from  Proposition \ref{ray limits in petal} and the uniqueness in Lemma \ref{walz}   that 
		$$\Psi(w)= \Psi_{ \Theta, n}^\sigma(w) \text{ and }\Psi_{\Theta\oplus \theta, n}^\sigma (w) = \Psi_{\Theta, n}^\sigma(\lambda^{-n}\xi_{\theta}(w))$$
		for any $w$ where both sides of the respective equations are defined.
	\end{proof}

	Just as the external coordinates $\psi_{\theta}^{h, f}$ can be related to $\psi_{0}^{h, f}$ for $\theta\in \mathcal{A}^f$, it follows from the uniqueness in Proposition \ref{parameter rays} that we can relate $\Psi_{\Theta,n}^{\sigma}$ to $\Psi_n^\sigma$. 
	\begin{corollary}\label{comparing first parameter rays}
		For any sufficiently large negative $\theta\in \mathcal{A}^f$, if $n$, $|\sigma|$, and $\ell(\sigma)$ are sufficiently large then 
		$$\Psi_{\langle 0, \theta\rangle, n}^\sigma(w) = \Psi_n^\sigma\circ \xi_\theta^\sigma(w)$$
		for all $w\in X_0^\sigma\cap (\xi_\theta^\sigma)^{-1}(X_0^\sigma).$
	\end{corollary}

	The maps $\Psi_{\Theta, n}^\sigma$ therefore provide analytic extensions of the external parameters $\Psi^\sigma$. In order to control the landing of parameter rays we need more than Propositions \ref{ray limits parameter} and \ref{parameter ray limits implosive}; we need to consider Log-multiplier families induced by near-parabolic renormalization. 

	\subsection{Induced Log-multiplier families}\label{induced families}
	
	Let us fix some rational $\tilde{p}/\tilde{q}\in [-{1}/{2}, {1}/{2}]$ and set $\tilde{f}= \mathcal{R}_{\tilde{p}/\tilde{q}} f$. 
	For any  $\tilde{\alpha}\in \mathbb{D}$, if $|\sigma|$ and $n\geq 0$ are sufficiently large then we can define
	$$\tilde{h}_{\tilde{\sigma}, \tilde{\alpha}}:= \mathcal{R}_fh_{\sigma, \alpha},$$
	where $\tilde{\sigma}= (\sigma, n)$ and $$\alpha = \frac{1}{n-\mu_{\tilde{p}/\tilde{q}}(\tilde{\alpha})}.$$
	Note that we can choose large $n$ in the definition above independent of $\sigma$ and $\tilde{\alpha}.$ 
	We define  $\tilde{\Sigma}$ to be a subset of $\mathbb{Z}^{N+1}$ such that 
	$\tilde{h}_{\tilde{\sigma}, \tilde{\alpha}}$ is defined for all $\tilde{\sigma}\in \tilde{\Sigma}$ and $\tilde{\alpha}\in \mathbb{D}$.
	We will simplify our notation by writing $L = L_{\tilde{p}/\tilde{q}}^f$ and $\tilde{h} = \tilde{h}_{\tilde{\sigma}, \tilde{\alpha}},$
	so $\tilde{h}= \mathcal{R}_fh.$ Also, we denote 
	$$\textbf{p}/\textbf{q}=\textbf{p}_n/\textbf{q}_n:= \mu_{p/q}\left(\frac{1}{n-\tilde{p}/\tilde{q}}\right).$$
	We can compute 
	$$\textbf{q}= \tilde{q}(nq+k_+^f)-q(\tilde{p}+\tilde{q}\mathfrak{c}_+^f),$$
	so our assumption that $\mathfrak{S}(p/q) = +1$ implies that $\tilde{h}^{\textbf{q}}\to L^{\tilde{q}}\circ f^{-\tilde{p}q}$ when $h\to f$ and $\tilde{h}\to \tilde{f}$.
	
	\begin{proposition}
		The maps $\tilde{h}$ form a log-multiplier family near $\tilde{f}$ indexed by $\tilde{\Sigma}.$
	\end{proposition}

	\begin{proof}
		It follows from the definition that  $\tilde{h}$ depends continuously on $\tilde{\sigma}\in \tilde{\Sigma}^*$ and holomorphically on $\tilde{\alpha}\in \mathbb{D}$.
		It follows from the definition of near-parabolic renormalization that $0$ is a fixed point of  $\tilde{h}$ and $\tilde{h}'(0) = e^{2\pi i \mu_{\tilde{p}/\tilde{q}}(\tilde{\alpha})}$, and it follows from Proposition \ref{near parabolic renormalization} that $\tilde{h}\to \tilde{f}$ when $|\tilde{\sigma}|\to \infty$ and $\tilde{\alpha}\to 0$.
	\end{proof}

	We will say that the Log-multiplier family $\{\tilde{h}\}$ is \textit{induced} by the Log-multiplier family $\{h\}$.	Just as we can lift the parabolic dynamics of $\tilde{f}$ to $L$, it is shown in \cite{nondegenerate} that we can lift the near-parabolic dynamics of $\tilde{h}$ to $h$.

	\begin{theorem}\label{lifting petals}
		For any $\tilde{\mathfrak{C}}>0$, there is a compact set $Y\subset\mathbb{C}$ and  integers $M, M'\geq 0$ such that if $\tilde{h}$ is sufficiently close to $\tilde{f}$ and $\tilde{\alpha}\in A(\tilde{\mathfrak{C}})$, then there exist functions
		\begin{align*}
			\Phi_{att}^{\tilde{h}, \tilde{L}}&: \Omega_{att}^{\tilde{h}, \tilde{f}}\cap T_{-M'}(\Omega_{att}^{\tilde{h}, \tilde{f}})\to \Omega_{att}^{h, f}\cup Y\\
			\tilde{\varphi}_{att}^{h, L}&: \Omega_{att}^{\tilde{h}, \tilde{f}}\cap T_{-M'}(\Omega_{att}^{\tilde{h}, \tilde{f}})\to \hat{\mathbb{C}}\\
		\end{align*}
		such that:
		\begin{enumerate}
			\item $\tilde{\varphi}_{att}^{h, L}$ is analytic and 
			$$\tilde{\varphi}_{att}^{h, L}(w) = \begin{cases}
			\varphi_{att}^{h, f}\circ \Phi_{att}^{\tilde{h}, \tilde{f}}(w) & \text{ if } \Phi_{att}^{\tilde{h}, \tilde{f}}(w)\in \Omega_{att}^{h, f},\\
			h^{-Mq}\circ\varphi_{att}^{h, f}\circ T_{M}\circ \Phi_{att}^{\tilde{h}, \tilde{f}}(w)& \text{ if } \Phi_{att}^{\tilde{h}, \tilde{f}}(w)\in Y,\\
			\end{cases}$$
			for some continuous branch of $h^{-Mq}$.
			\item $\tilde{\varphi}_{att}^{h, L}(0) = cv^{h, f}$.
			\item $\tilde{\varphi}_{att}^{h, L}\circ T_1=h^{\emph{\textbf{q}}} \circ  \tilde{\varphi}_{att}^{h, L}.$
			\item $\tilde{\varphi}_{att}^{h, L}\to \tilde{\varphi}_{att}^L$ and $\tilde{\varphi}_{rep}^{h, L}:= \tilde{\varphi}_{att}^{h, L}\circ T_{-\frac{1}{\tilde{\alpha}}}\to \tilde{\varphi}_{att}^L$ when $h\to f$ and $\tilde{h}\to \tilde{f}$.
			\item For any $w_0\in \Omega_{att}^{\tilde{f}}$ or $w_0'\in \Omega_{rep}^{\tilde{f}}$, if $\tilde{h}$ is sufficiently close to $\tilde{f}$ then 
			$$\Exp\circ T_{-\frac{1}{\alpha}}\circ \rho^{h, f}\circ \tilde{\varphi}_{att}^{h, L}(w) = \varphi_{att}^{\tilde{h},\tilde{f}}(w)\text{ or }\Exp\circ T_{-\frac{1}{\alpha}}\circ \rho^{h, f}\circ \tilde{\varphi}_{rep}^{h, L}(w) = \varphi_{rep}^{\tilde{h},\tilde{f}}(w)$$
			for all $w$ or $w'$ close to to $w_0$ or $w_0'$ respectively.
		\end{enumerate}
	\end{theorem}

	\begin{proof}
		See \cite[Theorem 4.1]{nondegenerate}.
	\end{proof}

	As a consequence of the Theorem \ref{lifting petals}, we can make the following observation:

	\begin{proposition}
		The Log-multiplier family $\{\tilde{h}_{\tilde{\sigma}, \tilde{\alpha}}\} $ has pre-petals.
	\end{proposition}

	\begin{proof}
		While we defined $\varphi_{-1, \iota}^{\tilde{f}}$ to be the $\varphi_\iota^f$ post-composed with a specific branch of $\tilde{f}^{-\tilde{q}}$, let us now show that we can instead define $\varphi_{-1, \iota}^{\tilde{f}}$ using the pre-petal parameters of $L$. We will only consider $\varphi_{-1, att}^{\tilde{f}}$, the repelling case is similar.
		
		Let us denote $g = \Exp\circ T_\delta \circ \rho$. Using the branch of $g^{-1}$ which, after continuous extension to the closures, maps  $\overline{U_0^{\tilde{f}}}$ to $\overline{U_0^L}$, for any $\theta\in \mathcal{A}^{\tilde{f}}$ we define the points
		$$\tilde{z}_\theta^L:= g^{-1}(z_\theta^{\tilde{f}})$$
		and Jordan domains
		$$\tilde{P}_{\theta, att}^L:= g^{-1}(P_{\theta, att}^{\tilde{f}}).$$
		Let $U_{-1}^L$ be the unique component of $f^{-q}(U_0^L)$ contained in $U_0^f$. As $cv^f\in U_0^L$ and $L^{\tilde{q}}\circ f^{-\tilde{p}q}$ maps $U_0^L$ to itself as a branched double covering, $L^{\tilde{q}}\circ f^{(-\tilde{p}-1)q}$ must map $U_0^L$ univalently onto $U_{-1}^L$.
		It follows from definition that the continuous extension of $L^{\tilde{q}}\circ f^{(-\tilde{p}-1)q}$ to $\overline{U_0^L}$ fixes zero. Hence 
		$$L^{\tilde{q}}\circ f^{(-\tilde{p}-1)q}(\tilde{z}_{-1}^{L}) = z_{-1}^f.$$		
		Similarly, $L^{\tilde{q}}\circ f^{(-\tilde{p}-1)q}$ must map $\tilde{P}_{-1, att}^L$ into the unique component of $f^{-q}(\tilde{P}_{0, att}^L)$ contained in $U_0^f$ which has $z_{-1}^f$ on its boundary. As $L^{\tilde{q}}\circ f^{(-\tilde{p}-1)q}$ is univalent in a neighborhood of $z_\Theta^L$, it follows that $\tilde{P}_{-1, att}^L$ is part of the unique component of $(L^{\tilde{q}}\circ f^{-\tilde{p}q})^{-1}(\tilde{P}_{0, att}^L)$ contained in $K_{\tilde{q}}^L$ which has $z_\Theta^L$ on its boundary.
		
		It follows from the above that there is some enriched angle $\Theta$ for $f$ of depth $\tilde{q}$ and with principle integer $-\tilde{p}$ such that $\tilde{z}_{-1}^L= z_\Theta^L$. 		
		Let us assume for now that 
		the image of $$\tilde{P}_{0,att}^L\subset L^{\tilde{q}}\circ f^{-\tilde{p}q}(P_{\Theta, att}^L)\subset P_{att}^f,$$
		so it follows from the above that 
		so we can define $$\tilde{P}_{att}^L= \varphi_{\Theta, att}^L\circ T_{\tilde{p}} \circ (\varphi_{att}^f)^{-1}(\tilde{P}_{att}^L).$$
		Thus for an integer $m\geq 0$ sufficiently large so that $f^{mq}(\tilde{P}_{att}^L)\subset P_{att}^f$,  we have
		$$\varphi_{-1, att}^{\tilde{f} }= \Exp\circ T_{\delta}\circ f^{mq}\circ (\varphi_{att}^f)^{-1}\circ \varphi_{\Theta, att}^L\circ T_{\tilde{p}} \circ (\varphi_{att}^f)^{-1}\circ \tilde{\varphi}_{att}^L.$$
		Thus we can similarly define 
		$$\varphi_{-1, att}^{\tilde{h}, \tilde{f}} := \Exp\circ T_{-\frac{1}{\alpha}}\circ (\varphi_{att}^{h, f})^{-1}\circ h^{mq}\circ \varphi_{\Theta, att}^{h, L}\circ T_{\tilde{p}} \circ \Phi_{att}^{\tilde{h}, \tilde{f}},$$
		and this map satisfies all of the desired properties when $h\to f$ and $\tilde{h}\to \tilde{f}$.
		In the general case, where the image of $T_{\tilde{p}} \circ \Phi_{att}^{\tilde{h}, \tilde{f}}$ is not contained in the domain of $\varphi_{\Theta, att}^{h, L}$, we need to analytically extend   first $\varphi_{att}^{h, f}$ and then $\varphi_{\Theta, att}^{h, L}$ by applying branches of iterates of $h^{-q}$; Theorem \ref{lifting petals} ensures that we need to do this extension only on a uniform compact subset of $\mathbb{C}$ and only using a uniform number of iterates of $h^{-q}$. 
	\end{proof}

	We set 
	$\tilde{\lambda}^{\tilde{\sigma}}:= (\lambda^\sigma)^{n\tilde{q}-\tilde{p}}$ and $\tilde{x}_0^{\tilde{\sigma}}:= {x_0^{\sigma}}/({2\tilde{\lambda}^{\tilde{\sigma}}})$; for all integers $m$ we define
	$$\tilde{X}_m^{\tilde{\sigma}}:= \{(\tilde{\lambda}^{\tilde{\sigma}})^m(x-iy)\in \mathbb{H}: |x|< \tilde{x}^{\tilde{\sigma}}, |\log y|< \log \tilde{\lambda}^{\tilde{\sigma}}\}.$$
	We will usually suppress the dependence on $\tilde{\sigma}$ in our notation; for example we will write $\tilde{\lambda} = \tilde{\lambda}^{\tilde{\sigma}}.$
	Let $\Theta_\infty= \langle \theta_m\rangle_{m=0}^\infty$ be the parabolic enriched angle for $L$ and set $\Theta_d= \langle \theta_m\rangle_{m=0}^d$ for all $d\geq 0$.
	By Proposition \ref{parabolic bubble ray} there exists some $d_0\gg 0$ such that $\overline{B_{\Theta_d}}\subset \tilde{P}_{rep}^L$ for all $d\geq (d_0-2)\tilde{q}$; we define the map 
	$$\xi_{\Theta_\infty}^{\tilde{\sigma}}(w):= w + \dfrac{\lambda^{n\tilde{q}}\sum_{m=0}^{\tilde{q}-1}\lambda^{n(m-\tilde{q})}b_{\theta_{d_0+m}}}{\tilde{\lambda}-1}.$$
	\begin{proposition}\label{induced external rays}
		If $h$ and $\tilde{f}$ are close to $f$ and $\tilde{f}$ respectively, then $\psi_{\Theta_{d_0}, n}^{h, f}\circ \xi_{\Theta_\infty}^{\tilde{\sigma}}$ is defined on $\tilde{X}_{0}$ and satisfies
		$$h^{\emph{\textbf{q}}}\circ \psi_{\Theta_{d_0}, n}^{h, f}\circ \xi_{\Theta_\infty}^{\tilde{\sigma}}(w) = \psi_{\Theta_{d_0}, n}^{h, f}\circ \xi_{\Theta_\infty}^{\tilde{\sigma}}(\tilde{\lambda}w)$$
		wherever both sides of the equation are defined. There is a compact set $\widehat{\tilde{X}}_{0, 0}^L\subset \tilde{P}_{rep}^L$ which avoids $ U^L$ such that $\psi_{\Theta_{d_0}, n}^{h, f}\circ \xi_{\Theta_\infty}^{\tilde{\sigma}}(\tilde{X}_{0})\xrightarrow{\subset} \widehat{\tilde{X}}_{0, 0}^L$ when $h\to f$ and $\tilde{h}\to \tilde{f}$.
	\end{proposition}
	\begin{proof}
		Proposition \ref{parabolic bubble ray} implies that $L^{\tilde{q}}\circ f^{-\tilde{p} q}(B_{\Theta_{d_0+\tilde{q}}})= B_{\Theta_{d_0}}.$
		Proposition \ref{ray limits} therefore implies if $h$ is  close to $f$ and $\tilde{h}$ is close to $\tilde{f}$ then $$h^{\textbf{q}}\circ \psi_{\Theta_{d_0}, n}^{h, f}(w) = \psi_{\Theta_{d_0-\tilde{q}, n}}(\lambda^{-\tilde{p}} w)$$
		for all $w\in X_0$.
		Additionally, it follows from Propositions \ref{parabolic bubble ray} and \ref{ray limits in petal} that 
		$$\psi_{\Theta_{d, n}}^{h, f}(w) = \psi_{\Theta_{d-1, n}}(\lambda^{-n}\xi_{\theta_d}(w))$$
		for all $w\in X_0$ and $d\geq 1$. 
		Setting $$\tilde{b}_n:=\sum_{m=0}^{\tilde{q}-1}\lambda^{n(m-\tilde{q})}b_{\theta_{d_0+m}}$$ and proceeding inductively, as $\xi_{\theta_d}(w) = w+b_\theta^h$ it follows that 
		$$\psi_{\Theta_{d_0, n}}^{h, f}(w) = \psi_{\Theta_{d_0-\tilde{q}, n}}^{h, f}\left(\lambda^{-n\tilde{q}}w+\tilde{b}_n\right)$$
		for all $w\in X_0$. 
		It follows by a straightforward computation that 
		$$h^{\textbf{q}}\circ \psi_{\Theta_{d_0}, n}^{h, f}\circ \xi_{\Theta_\infty}^{\tilde{\sigma}}(w) = \psi_{\Theta_{d_0}, n}^{h, f}\circ \xi_{\Theta_\infty}^{\tilde{\sigma}}(\tilde{\lambda}w)$$ 
		for all $w\in (\xi_{\Theta_\infty}^{\tilde{\sigma}})^{-1}(X_0)$. 
		Additionally, it follows from the above construction that $\psi_{\Theta_{d_0}, n}^{h, f}$ is defined on the set
		$$Y_n:=\{w\in \mathbb{H}: w+it\in T_{\lambda^{\tilde{p}}\tilde{b}_n}(X_{-n\tilde{q}+\tilde{p}}) \text{ and }x-it'\in T_{\tilde{\lambda}\tilde{b}_n}(X_{n\tilde{q}-\tilde{p}}) \text{ for some }t, t'\geq 0 \}.$$
		As $\tilde{b}_n\to 0$ when $n\to \infty$, another straightforward computation implies 
		$$\xi_{\Theta_\infty}^{\tilde{\sigma}}(\tilde{X}_0)\subset Y$$
		 when $h$ and $\tilde{h}$ are close to $f$ and $\tilde{f}$ respectively.
		Hence $\psi_{\Theta_{d_0}, n}\circ \xi_{\Theta_\infty}^{\tilde{\sigma}}$ is defined on $\tilde{X}_0$ and Propositions \ref{ray limits} and \ref{ray limits in petal} imply that
		$$\psi_{\Theta_{d_0}, n}^{h, f}\circ \xi_{\Theta_\infty}^{\tilde{\sigma}}\xrightarrow{\subset}\bigcup_{d={d_0-\tilde{q}}}^{d+\tilde{q}}\left(\overline{\bigcup_{m\in \mathbb{Z}}\hat{X}_{\Theta_d, m}^L}\right)\subset \bigcup_{d={d_0-\tilde{q}}}^{d+\tilde{q}}\overline{B_{\Theta_d}^L}\subset (P_{rep}^L\setminus U^L)$$
		when $h\to f$ and $\tilde{h}\to \tilde{f}.$
	\end{proof}

	\begin{corollary}\label{inductive compatibility}
		The Log-multiplier family $\{\tilde{h}_{\tilde{\sigma}, \tilde{\alpha}}\}$ has compatible external rays and pre-petals.
	\end{corollary}

	\begin{proof}
		It follows from Theorem \ref{lifting petals} that when $\tilde{h}$ is sufficiently close to $\tilde{f}$ we can define the external coordinate of $\tilde{h}$ by
		$$\psi_0^{\tilde{h}, \tilde{f}}:=\Exp \circ T_{-1/\alpha}\circ \rho^{h, f}\circ \psi_{\Theta_{d_0}, n}^{h, f}\circ \xi_{\Theta_\infty}^{\tilde{\sigma}}: \tilde{X}_0\to \mathbb{C}.$$
		It follows from Proposition \ref{induced external rays} that these external coordinates satisfy the hypotheses so that the Log-multiplier family has pre-petals. 
		
		We observe that for any $\tilde{\theta}\in \mathcal{A}^{\tilde{f}}$ with principle integer $j\geq 0$,  there exists some enriched angle $\Theta$ for $f$ of depth $d_0+j\tilde{q}$ so that 
		$$\psi_{\tilde{\theta}}^{\tilde{h}, \tilde{f}}(w) = \tilde{h}^{k_+^{\tilde{f}}}\circ\Exp \circ T_{-1/\alpha}\circ \rho^{h, f}\circ \psi_{\Theta, n}^{h, f}\circ \xi_{\Theta_\infty}^{\tilde{\sigma}}(\lambda)$$
		when $\tilde{h}$ is sufficiently close to $\tilde{f}$. 
		It follows from Proposition \ref{parabolic bubble ray} that if $\tilde{\theta}$ is sufficiently large and negative, which forces $B_\Theta$ to be close to zero, then $B_\Theta$ is a descendant of $B_{\Theta_{d_0}}.$ We can therefore extend $\psi_0^{\tilde{h}, \tilde{f}}$ using Propositions \ref{ray limits extending} and \ref{ray limits in petal} so that the hypotheses of compatible external rays and pre-petals are satisfied. 
	\end{proof}

	Thus there is some $\tilde{n}_0> 0$ such that we can  define the parameter coordinate $$\Psi^{\tilde{\sigma}}: \bigcup_{\tilde{n}\geq \tilde{n}_0}\tilde{X}_{-\tilde{n}}\to \mathbb{D}$$   when $|\tilde{\sigma}|$ and $\ell(\tilde{\sigma})$ are sufficiently large. 
	Just as the external coordinates of $\tilde{h}$ are defined in terms of the external coordinates for $h$, we can define $\Psi^{\tilde{\sigma}}$ in terms of $\Psi^\sigma$.
	\begin{proposition}\label{relating parameter coordinates}
		If $|\tilde{\sigma}|$ and $\ell(\tilde{\sigma})$ are sufficiently large, then
		$$\frac{1}{n- \mu_{\tilde{p}/\tilde{q}}\circ \Psi^{\tilde{\sigma}}\left(\tilde{\lambda}^{-\tilde{n}_0}w\right)}= \Psi_{\Theta_{d_0+\tilde{n}_0\tilde{q}}, n}^{\sigma}\left(\lambda^{\tilde{n}_0\tilde{p}}\cdot \xi_{\Theta_\infty}^{\tilde{\sigma}}(w)\right)$$
		for all $w\in \tilde{X}_0\cap (\xi_{\Theta_\infty}^{\tilde{\sigma}})^{-1}(X_0).$
	\end{proposition}
	
	\begin{proof}
		Let $\mathcal{D}\subset Dom(\tilde{\varphi}_{rep}^L)$ be a Jordan domain containing $$\bigcup_{d=d_0-\tilde{q}}^{d_0+\tilde{q}}(\tilde{\varphi}_{rep}^{L})^{-1}\left(\overline{B_{\Theta_{d}}^L}\right).$$
		It follows from the definition of the external rays of $\tilde{h}$ and the proof of Proposition \ref{parameter rays} that for any $w\in \tilde{X}_0$,  when $|\tilde{\sigma}|$ and $\ell(\tilde{\sigma})$ are sufficiently large
		$\tilde{\alpha}= \Psi^{\tilde{\sigma}}(\tilde{\lambda}^{-\tilde{n}_0}w)$ is the unique choice of $\tilde{\alpha}$ with $\tilde{n}_0-\frac{1}{\tilde{\alpha}}\in \mathcal{D}$ satisfying 
		$$\tilde{h}^{\tilde{n}_0\tilde{q}}(cv^{\tilde{h}, \tilde{f}}) = \varphi_{rep}^{\tilde{h}, \tilde{f}}\left(\tilde{n}_0-\frac{1}{\tilde{\alpha}}\right) = \psi_0^{\tilde{h}, \tilde{f}}(w),$$
		or equivalently
		$$h^{\tilde{n}_0\textbf{q}}(cv^{h, f})= \tilde{\varphi}_{rep}^{h, L}\left(\tilde{n}_0-\frac{1}{\tilde{\alpha}}\right) = \psi_{\Theta_{d_0}, n}^{h, f}\circ \xi_{\Theta_\infty}^{\tilde{\sigma}}(w).$$
		For all $z\in \tilde{X}_0$ and $0 \leq m \leq \tilde{n}_0$, let us define
		$$\phi(\tilde{\lambda}^{-m}z)= T_{-m+1/\tilde{\alpha}}\circ (\tilde{\varphi}_{rep}^{h, L})^{-1}\circ \psi_{\Theta_{d_0}, n}^{h, f}\circ \xi_{\Theta_\infty}^{\tilde{\sigma}}(z),$$
		taking the branch of $ (\tilde{\varphi}_{rep}^{h, L})^{-1}$ which approximates  $(\tilde{\varphi}_{rep}^{ L})^{-1}$. The map $\phi$ is holomorphic on $\tilde{X}_{-m}$ for all $0\leq m \leq \tilde{n}_0$. 
		
		Let us suppose for now that the image of $\phi$ is contained in the domain of $\tilde{\varphi}_{att}^{h, f}$, so 
		$$\tilde{\varphi}_{att}^{h, f}\circ \phi = \psi_{\Theta_{d_0}, n}^{h, f}\circ \xi_{\Theta_\infty}^{\tilde{\sigma}}.$$
		Thus for any integer $j\geq 0$, if $|\tilde{\sigma}|$ and $\ell(\tilde{\sigma})$ are sufficiently large then $$h^{(n-j)q+k_+^f}(cv^{h, f} )= \psi_{\lceil\Theta_{d_0+\tilde{n}_0\tilde{q}}\rceil-j}(\lambda^{\tilde{n}_0\tilde{p}-j} \xi_{\Theta_\infty}^{\tilde{\sigma}}(w)).$$
		Keeping $w$ fixed, $n-\frac{1}{\alpha}$ converges to a parameter $\delta$ when  $|\sigma|\to \infty$ such that 
		$$L_\delta \circ f^{-jq}(cv^f)\in \hat{X}^{L_\delta}_{\lceil\Theta_{d_0+\tilde{n}_0\tilde{q}}-j, \tilde{n}_0\tilde{p}\rceil-j},$$
		so $\delta\in \mathcal{X}_{\Theta_{d_0+\tilde{n}_0\tilde{q}}, \tilde{n}_0\tilde{p}}.$
		It therefore follows from the uniqueness in Proposition \ref{ray limits parameter} 
		that when $|\tilde{\sigma}|$ and $\ell(\tilde{\sigma})$ are sufficiently large, 
		$$\frac{1}{n-\mu_{\tilde{p}/\tilde{q}}(\tilde{\alpha})}=\alpha = \Psi_{\Theta_{d_0+\tilde{n}_0\tilde{q}}, n}^{\sigma}(\lambda^{\tilde{n}_0\tilde{p}}\xi_{\Theta_\infty}^{\tilde{\sigma}}(w)).$$
		The general case follows similarly by the following lemma:
		\begin{lemma}
			If $|\tilde{\sigma}|$ and $\ell(\tilde{\sigma})$ are sufficiently large, then there exists a simply connected set $Y\subset Dom(\phi)$ containing both $w$ and $\tilde{\lambda}^{-{\tilde{n}_0}}(w)$  such that $\tilde{\varphi}_{rep}^{h, L}$ can be analytically extended to $\phi(Y).$
		\end{lemma}
		\begin{proof}
			Fixing some $z\in \overline{B_{\Theta_{d_0}}^L}\setminus \{z_{\Theta_{d_0}}^L\}$, let us first show that the set
			 $$V:=T_{-z}\left(\bigcup_{d=d_0-\tilde{q}}^{d_0}(\tilde{\varphi}_{rep}^{L})^{-1}\left(\overline{B_{\Theta_{d}}^L}\setminus\{z_{\Theta_d}^L\}\right)\right)$$
			is homotopic in $\mathbb{C}$ rel $\mathbb{Z}$ to the line segment $[0, 1]$.
				As $B_{\Theta_d}^L$ is a Jordan domain whose closure is intersects the closure of $B_{\Theta_{d'}}^L$ only when $d'=d\pm1$ and only at  
			$z_{\Theta_{d'}}^L$ for all $d\geq 1$, it is clear that $V$ is homotopic in $\mathbb{C}$ to simple curve $\gamma:[0, 1]\to $ which satisfies $\gamma(0) = 0$ and $\gamma(1)=1$. 
			If the image of $\gamma$ is not homotopic rel $\mathbb{Z}$ to $[0,1]$  then there is some $t_0\in (0,1)$ so that either $\gamma(t_0)\geq2$ or $\gamma(t_0)\leq-1$
			Let us assume that both $x$ and $x+j$ belong to $V$ for some $x\in \mathbb{C}$ and $j\in \mathbb{Z}$. 
			As $\tilde{\varphi}_{rep}^{L}(x)\in E_d^L\text{ and }\tilde{\varphi}_{rep}^{L}(x)\in E_{d-j\tilde{q}}^L$ for some $d\geq j\tilde{q}$, it follows from the definition of $V$ that $|j|\leq 1.$ Thus we can choose $\gamma$ so that if $\gamma(t)-\gamma(t')=j\in \mathbb{Z}$ for some $t, t'\in [0,1]$, then $|j|=1$. Moreover, setting $x= z$ above it follows that we can choose $\gamma$ so that the image of $\gamma$ avoids $-1.$ 
			It therefore follows from the intermediate value theorem that the image of $\gamma$ is homotopic $[0,1].$
			
			For any integer $m$, we set $$Y_m= \{\tilde{\lambda}^{-m}(x-iy)\in \tilde{X}: y> \lambda^{-1}\}.$$
			Thus $w\in Y_0$ and $\phi(Y_0)\xrightarrow{\subset} V$ when $|\tilde{\sigma}|\to \infty$.
			Setting $Y= \bigcup_{m =0}^{\tilde{n}_0}Y_{-m}$, it follows from the above that $\phi(Y)$ is homotopic rel $\mathbb{Z}$ to the line segment $[0, \tilde{n}_0+1]$. Similarly to using $(\rho^{h, f})^{-1}$ to extend $\varphi_{att}^{h ,f}$, we can analytically extend $\tilde{\varphi}_{att}^{h, L}$ to $\phi(Y)$; it follows from our homotopy condition on $\phi(Y)$ that this extension is well-defined.
		\end{proof}
	\end{proof}
	
	The analysis in this section gives us precise control on the geometry of parameter rays when $|\sigma|\to \infty$. As a consequence of that control we have the following proposition which will be one of the main tools we use in proving the Theorem \ref{main}.

		\begin{proposition}\label{cage}
			There exists $C>0$ such that  if $|\tilde{\sigma}|$ and $\ell(\tilde{\sigma})$ are sufficiently large then there exists a curve
			$$\gamma^{\tilde{\sigma}}: [-1, 1]\to \mathbb{H}$$
			satisfying:
			\begin{enumerate}
				\item $\gamma^{\tilde{\sigma}}(t)$ is contained in the image of $\Psi^{{\sigma}}$ for all $t\in (-1, 1).$
				\item $\gamma^{\tilde{\sigma}}(\pm 1)=\frac{1}{n\pm 1}$
				\item For all $t\in [-1, 1]$, 
				$$\left|n-\frac{1}{\gamma^{\tilde{\sigma}}(t)}\right|< {C}.$$
			\end{enumerate}
		\end{proposition}

		\begin{proof}
			It follows from Proposition \ref{relating parameter coordinates} that when $|\tilde{\sigma}|$ and $\ell(\tilde{\sigma})$ are sufficiently large we can analytically extend 
			$\Psi^{\sigma}_{\Theta_{d_0}+\tilde{n}_0\tilde{q}, n}$ so that there exists a  curve $\gamma^{{\sigma}}_{n, 0}: [0, 1]\to \mathbb{C}$ satisfying:
			\begin{enumerate}
				\item $\gamma^{{\sigma}}_{n, 0}(t)$ is contained in the image of $\Psi^{\sigma}_{\Theta_{d_0}+\tilde{n}_0\tilde{q}, n}$ for all $t\in (0, 1]$
				\item  $\gamma^{\sigma}_{n, 0}(0)= \frac{1}{n}$ and $\gamma^{{\sigma}}_{n, 0}(1)\in \Psi^{\sigma}_{\Theta_{d_0}+\tilde{n}_0\tilde{q}, n}(X_0)$.
				\item For all $t\in [0, 1]$, 
				$$n-\frac{1}{ \gamma^{\sigma}_{n, 0}(t)}\in \mathbb{D}.$$
			\end{enumerate}
			It follows from finitely many applications of Propositions \ref{ray limits parameter} and \ref{parameter ray limits implosive} 
			that there is some $R_1>0$ so that when $|\tilde{\sigma}|$ and $\ell(\tilde{\sigma})$ are sufficiently large there is a curve $\gamma^{\sigma}_{n, 1}:[1, 2]\to\mathbb{C}$ satisfying:
			\begin{enumerate}
				\item $\gamma^{\sigma}_{n, 1}(t)$ is contained in the image of $\Psi_{\Theta_1, n}^\sigma$ for all $t\in [1, 2].$
				\item $\gamma^{\sigma}_{n, 1}(1) = \gamma^{{\sigma}}_{n, 0}(1)$ and $\gamma^{\sigma}_{n, 1}(2)\in \Psi_{\Theta_1, n}^\sigma(X_0)$.
				\item For all $t\in [1, 2]$,
				$$n-\frac{1}{\gamma^{\sigma}_{n, 1}(2)}\in \mathbb{D}_{R_1}.$$
			\end{enumerate}
			It follows from Proposition \ref{parameter rays} and Corollary \ref{comparing first parameter rays} that there is some $R_2>0$ so that when $|\tilde{\sigma}|$ and $\ell(\tilde{\sigma})$ are sufficiently large there is a curve $\gamma^{\sigma}_{n, 2}:[-1, 1]\to \mathbb{C}$ satisfying:
			\begin{figure}
				\begin{center}
					\def\svgwidth{2.5in}
					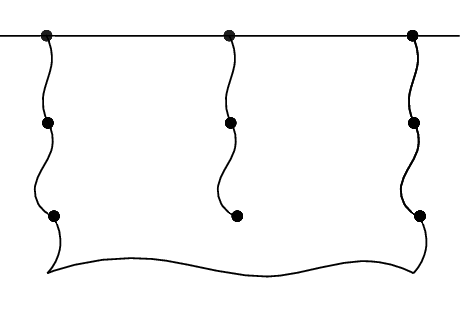
					\caption{The curves which define $\gamma^{\tilde{\sigma}}$.}
					\label{cage_fig}
				\end{center}
			\end{figure}
			\begin{enumerate}
				\item $\gamma^{\sigma}_{n, 2}(t)$ is contained in the image of $\Psi^\sigma$ for all $t\in [-1, 1]$.
				\item $\gamma^{\sigma}_{n, 2}(\pm 1) =  \gamma^{\sigma}_{n\pm 1, 1}(2)$.
				\item For all $t\in [-1, 1]$, 
				$$n-\frac{1}{\gamma^{\sigma}_{n, 2}(t) }\in \mathbb{D}_{R_2}.$$
			\end{enumerate}
			We form $\gamma^{\tilde{\sigma}}$ by composing $\gamma_{n\pm1, 0}^\sigma$, $\gamma_{n\pm1, 1}^\sigma$, and $\gamma_{n, 2}^\sigma$.
		\end{proof}

	Let us now note that while in this section we have only considered upper renormalizations, we can repeat all of the above with lower renormalizations. The only distinction is that for $\tilde{f} = \mathcal{R}_{\tilde{p}/\tilde{q}}^-f$ and $\tilde{h}_{\tilde{\sigma}, \tilde{\alpha}}= \mathcal{R}_f^-h_{\sigma, \alpha}$, we need to slightly adjust the relationship between $\tilde{\alpha}$ and $\alpha$. In particular, we would like to have
	$$\alpha = \frac{1}{n+ \mathfrak{c}_-^{h_{\sigma, \alpha}, f}+ \mu_{\tilde{p}/\tilde{q}}(\tilde{\alpha})},$$
	so $\tilde{h}_{\tilde{\sigma}, \tilde{\alpha}}'(0) = \Exp\circ \mu_{\tilde{p}/\tilde{q}}(\tilde{\alpha})$ by Proposition \ref{horn map perturbed}. As $\mathfrak{c}_-^{h_{\sigma, \alpha}, f}\to \mathfrak{c}_-^f$, the above equation implicitly defines $\alpha$ in terms of $\tilde{\alpha}$; for details see the proof of Proposition \ref{change of coordinates} below. Equipped with the functions $\tilde{h}_{\tilde{\sigma}, \tilde{\alpha}}$, the rest of the analysis above can be done with minimal change. In particular, we have the following version of Proposition \ref{cage} in the lower renormalization case:
	\begin{proposition}\label{cage lower}
		There exists $C>0$ such that  if $|\tilde{\sigma}|$ and $\ell(\tilde{\sigma})$ are sufficiently large then there exists a curve
		$$\gamma^{\tilde{\sigma}}: [-1, 1]\to \mathbb{H}$$
		satisfying:
		\begin{enumerate}
			\item $\gamma^{\tilde{\sigma}}(t)$ is contained in the image of $\Psi^{{\sigma}}$ for all $t\in (-1, 1).$
			\item $\gamma^{\tilde{\sigma}}(\pm 1)=\frac{1}{n\pm 1+ \mathfrak{c}_-^{h, f}}$, where $h= h_{\sigma, \gamma^{\tilde{\sigma}}(\pm 1)}$
			\item For all $t\in [-1, 1]$, 
			$$\left|n-\frac{1}{\gamma^{\tilde{\sigma}}(t)}\right|< {C}.$$
		\end{enumerate}
	\end{proposition}
	
	\section{Bounding limbs in $\mathcal{M}$}\label{bounds section}
	
	We are now ready to apply the analysis of the previous sections to quadratic polynomials. First let us recall some classical facts about quadratic polynomials; for more details on the iteration of polynomials we refer the reader to \cite{Orsay1} and \cite{Orsay2}.
	
	We consider the family of quadratic polynomials $$Q_\zeta(z) := \Exp(\zeta)\cdot z+z^2$$
	with $\zeta\in \mathbb{C}$. The set $\mathcal{M}$ is all the parameters $\zeta$ such that $Q_\zeta^n(cv^{Q_\zeta})\not\to \infty$ when $n\to \infty$. A parameter $\zeta\in \mathcal{M}$ is \textit{parabolic} when $Q_\zeta\in \mathcal{F}^\flower$. The set $\mathcal{M}\setminus \{\zeta\}$ has exactly two connected components in this case; we denote by $\mathcal{L}_{\zeta}$ the unique bounded component. Moreover there is a unique component of the interior of $\mathcal{L}_{\zeta}$ which has $\zeta$ on its boundary; this component is called the \textit{main hyperbolic component} of $\mathcal{L}_\zeta$ and we will denote it by $\mathcal{H}_\zeta$.
	For any parabolic parameter $\zeta$, there exists a set $\mathcal{M}_\zeta\subset \mathcal{L}_\zeta$ and a covering map $\tau_\zeta: \mathcal{M}\setminus \mathbb{Z}\to \mathcal{M}_\zeta$ such that $\tau_\zeta(w)\to \zeta$ when $w$ converges to $\infty$ or an integer. The tuning map $\tau_\zeta$ is dynamical in the following sense: $\tau_\zeta$ restricts to an analytic covering map from the upper half-plane $-\mathbb{H}$ to $\mathcal{H}_\zeta$ and for $\alpha = \tau_\zeta(w)$ is the unique parameter in $\mathcal{H}_\zeta$ such that $Q_\alpha$ has an attracting periodic cycle with multiplier $\Exp(w).$ We will also set $\tau_\infty(z) = z$.
	
	We define \textit{combinatorics} to be sequences $\omega= \langle\alpha_n\rangle_{n=1}^N$ of 
	parabolic parameters in $\mathcal{M}\setminus \mathbb{Z}$.  We define the \textit{depth} of the combinatorics to be the length of the sequence. As for other sequences, we will use $\oplus$ to denote the concatenation operation on combinatorics. 
	We will say that $\omega$ is the combinatorics of a parabolic parameter $\zeta\in \mathcal{M}$ if 
	$\mathcal{M}_\zeta$ is the image of $\tau_{\alpha_1}\circ \cdots \circ \tau_{\alpha_N}$, or equivalently $\tau_{\alpha_1}\circ \cdots \circ \tau_{\alpha_N}(0) = \zeta$.
	We will also say that the combinatorics of $\infty$ is the empty sequence. 
	If $\alpha_n=p_n/q_n$ is a reduced rational number for all $n$, then the combinatorics $\omega$ are called \textit{satellite} and we denote
	$$\|\omega\| = \prod_{n=1}^Nq_n^2.$$
	As any satellite combinatorics $\omega$ of depth $d$ never include the parameter $0/1$, we note that
	$\|\omega\|\geq 2^{2d}$.

	\subsection{Satellite towers}
	For any integer $N\geq 0$ or $N=\infty$, we define a \textit{height $N$ quadratic parabolic tower} $\mathcal{T}$ to be
	a sequence $( f_n)_{n=0}^N$ such that
	\begin{enumerate}
		\item $f_0= Q_\zeta$ for some $\zeta\in \mathbb{C}$,
		\item $f_{n}\in \mathcal{F}^\flower$ for all $0\leq n \leq N$, and
		\item  $f_{n+1}$ is a parabolic renormalization of either $f_n$ or $f_n^*$ for all $0\leq n < N$.
	\end{enumerate}
	 We will  say that $\mathcal{T}$ is a \textit{satellite tower} if $f_{n}$ has a parabolic fixed point at zero for all $0\leq  n\leq N$. 
	 
	 Let us now fix a  satellite tower $\mathcal{T}= ( f_n)_{n=0}^N$ with $0\leq N< \infty$, so there is a sequence $( p_n/q_n)_{n=0}^{N}$ of rational numbers in $(-1/2, 1/2]$ so that $f_n$ has a $p_n/q_n$-parabolic fixed point at zero for all $n$.

	 \begin{proposition}\label{change of coordinates}
	 	There exists $K^\mathcal{T}\geq 2$ and such that for $\{\sigma \in \mathbb{Z}^N: |\sigma|\geq K^{\mathcal{T}}\}:= \Sigma^\mathcal{T}$ and any $\sigma= \langle k_n\rangle_{n=1}^N\in \Sigma^\mathcal{T}$, there exists an injective holomorphic or anti-holomorphic univalent function $\mu_\sigma^\mathcal{T}:\mathbb{D}\to \mathbb{C}$ satisfying:
	 	\begin{enumerate}
	 		\item The collection of maps $$h_{\sigma, \alpha}^\mathcal{T}:=\mathcal{R}_{f_{N-1}}^\pm\circ\cdots\circ \mathcal{R}_{f_{1}}^\pm {Q}_{\mu_\sigma^\mathcal{T}(\alpha)},$$
	 		where the upper or lower near-parabolic renormalization relative to $f_n$ is taken when $f_{n+1}$ is an upper or lower parabolic renormalization respectively of $f_n$ or $f_n^*$  for all $n$, 
	 		form a Log-multiplier family near $f_N$ parameterized by $\Sigma^\mathcal{T}.$
	 		\item\label{combinatorics of near tower} There exists some $p_\sigma^\mathcal{T}/q_\sigma^\mathcal{T}$ and $\alpha_\sigma^\mathcal{T}$ with combinatorics $\omega_{\sigma}^\mathcal{T}$ such that either
	 		$$\mu_\sigma^\mathcal{T}(z) = \tau_{\alpha_\sigma^\mathcal{T}}\circ \mu_{p_\sigma^\mathcal{T}/q_\sigma^\mathcal{T}}(z) \text{ or }\mu_\sigma^\mathcal{T}(z) =\tau_{\alpha_\sigma^\mathcal{T}}\circ \mu_{p_\sigma^\mathcal{T}/q_\sigma^\mathcal{T}}(-\overline{z})$$
	 		for all $z\in\mathbb{H}\cap \mathbb{D}.$ 
	 		If $N= 0$ then $\alpha_\sigma= \infty$ and $p_\sigma^\mathcal{T}/q_{\sigma}^{\mathcal{T}} = p_0/q_0$, otherwise 
	 		$$\omega_{\sigma}^\mathcal{T} = \begin{cases}
	 		\omega_{\sigma'}^{\mathcal{T}'}&\text{ if }f_N= \mathcal{R}_{p_N/q_N}^+f_{N-1} \text{ or }f_N= \mathcal{R}_{p_N/q_N}^+f_{N-1}^*,\\
	 		\omega_{\sigma'}^{\mathcal{T}'} \oplus \langle p_{\sigma'}^{\mathcal{T}'}/q_{\sigma'}^{\mathcal{T}'}\rangle 
	 		&\text{ if }f_N= \mathcal{R}_{p_N/q_N}^-f_{N-1} \text{ or }f_N= \mathcal{R}_{p_N/q_N}^-f_{N-1}^*,\\
	 		\end{cases}$$
	 		and 
	 		$$p_\sigma^{\mathcal{T}}/ q_\sigma^{\mathcal{T}} =\begin{cases}
	 		p_{\sigma'}^{\mathcal{T}'}/ q_{\sigma'}^{\mathcal{T}'}\oplus \mathfrak{S}(p_{\sigma'}^{\mathcal{T}'}/q_{\sigma'}^{\mathcal{T}'})/k_N\oplus -p_N/q_N & \text{ if }f_N= \mathcal{R}_{p_N/q_N}^+f_{N-1},\\
	 		p_{\sigma'}^{\mathcal{T}'}/ q_{\sigma'}^{\mathcal{T}'}\oplus -\mathfrak{S}(p_{\sigma'}^{\mathcal{T}'}/q_{\sigma'}^{\mathcal{T}'})/k_N\oplus -p_N/q_N
	 		&\text{ if }f_N= \mathcal{R}_{p_N/q_N}^+f_{N-1}^*,\\
	 		-1/k_N\oplus {p_N/q_N} & \text{ if }f_N= \mathcal{R}_{p_N/q_N}^-f_{N-1},\\
	 		1/k_N\oplus {p_N/q_N} & \text{ if }f_N= \mathcal{R}_{p_N/q_N}^-f_{N-1}^*.
	 		\end{cases} $$
	 		\item\label{distortion tower}  
	 		For all $\sigma \in \Sigma^\mathcal{T}$ and $z, w\in \mu_{p_\sigma^\mathcal{T}/q_{\sigma}^{\mathcal{T}}}(\mathbb{D})$ we have
	 		$$\frac{2^{-d}|z-w|}{\|\omega_\sigma^\mathcal{T}\|}\leq \left|\mu_\sigma^\mathcal{T}\circ (\mu_{p_\sigma^\mathcal{T}/q_{\sigma}^{\mathcal{T}}})^{-1}(z)- \mu_\sigma^\mathcal{T}\circ (\mu_{p_\sigma^\mathcal{T}/q_{\sigma}^{\mathcal{T}}})^{-1}(w)\right|\leq \frac{2^{d}|z-w|}{\|\omega_{\sigma}^{\mathcal{T}}\|},$$
	 		where $d$ is the depth of $\omega_\sigma^\mathcal{T}$.
	 	\end{enumerate}
	 \end{proposition}
 
 	\begin{proof}
 		If $N= 0$, so $\sigma= \infty \in \mathbb{Z}^0$, then we can take $\mu_{\sigma}^\mathcal{T}= \mu_{p_0/q_0}$. So we now assume that $N>0$ and that the proposition holds for smaller values of $N$. Denoting $\mathcal{T}'= (f_n)_{n=0}^{N-1}$, it follows from the inductive hypothesis that there is some 
 		$\Sigma^{\mathcal{T}'}\subset \mathbb{Z}^{N-1}$ so that the univalent function $\mu_{\sigma'}^{\mathcal{T}'}:\mathbb{D}\to \mathbb{C}$ is defined for all $\sigma'\in \Sigma^{\mathcal{T}'}.$ Let us  fix some $\sigma = \langle k_n\rangle_{n=1}^N$ with $\sigma':= \langle k_n\rangle_{n=1}^{N-1}\in \Sigma^{\mathcal{T}'}$.
 		
 		We assume for now that $f_N$ is an upper parabolic renormalization of $f_{N-1}$.
 		For any $\alpha\in \mathbb{D}$ and $k_N\geq 2$ we set
 		$$\alpha'=
 		\frac{1}{k_N- \mu_{p_N/q_N}(\alpha)}= \mu_{1/k_N}\circ \mu_{p_N/q_N}(\alpha)$$
 		and we define $\mu_\sigma^\mathcal{T}(\alpha) := \mu_{\sigma'}^{\mathcal{T}'}(\alpha')$.
 		Thus there exists  some $K\geq 2$ so that $\mu_\sigma^\mathcal{T}$ is defined on $\mathbb{D}$ for all $\sigma \in \Sigma^{\mathcal{T}'}\times \{k\in \mathbb{Z}: k\geq K\}$, and the functions
 		$h_{\sigma, \alpha}^\mathcal{T}:= \mathcal{R}_{f_{N-1}}^+h_{\sigma', \alpha'}^{\mathcal{T}'}$ form a Log-multiplier family near $f_N$. 
 		We set $\alpha_{{\sigma}}^\mathcal{T}= \alpha_{{\sigma}'}^{\mathcal{T}'}$ and
 		$$p_{\sigma}^\mathcal{T}/q_\sigma^\mathcal{T} = p_{\sigma'}^{\mathcal{T}'}/ q_{\sigma'}^{\mathcal{T}'}\oplus \mathfrak{S}(p_{\sigma'}^{\mathcal{T}'}/q_{\sigma'}^{\mathcal{T}'})/k_N\oplus -p_N/q_N,$$
 		so it follows from Proposition \ref{adding fractions} that $\mu_{\sigma}^\mathcal{T} = \mu_{\sigma'}^{\mathcal{T}'}\circ (\mu_{p_{\sigma'}^{\mathcal{T}'}, q_{\sigma'}^{\mathcal{T}'}})^{-1}\circ \mu_{p_\sigma^\mathcal{T}/q_\sigma^\mathcal{T}}$. We note that even when $\mu_{\sigma'}^{\mathcal{T}'}$ is anti-holomorphic, our inductive hypothesis that $h_{\sigma', \alpha'}^{\mathcal{T}'}$ depends holomorphically on $\alpha$ ensures that $h_{\sigma, \alpha}^\mathcal{T}$ depends holomorphically on $\alpha$.  
 		
 		If instead $f_N$ is an upper parabolic renormalization of $f_{N-1}^*$, then for $\alpha\in \mathbb{D}$ and $k_N\geq  2$ we can define
 		$${\alpha'}=
 		\frac{-1}{k_N+\mu_{-p_N/q_N}(-\overline{\alpha})}= \mu_{-1/k_N}\circ \mu_{-p_N/q_N}(-\overline{\alpha})$$
 		and $\mu_{\sigma}^\mathcal{T} (\alpha) = \mu_{\sigma'}^{\mathcal{T}'}(\alpha').$
 		Thus there exists some $K\geq 2$ so that $\mu_{\sigma}^\mathcal{T}$ is defined on $\mathbb{D}$ for all  $\sigma \in\Sigma^{\mathcal{T}'}\times \{k\in \mathbb{Z}: k\geq K\}$
 		and the functions $h_{\sigma, \alpha}^{\mathcal{T}}:= \mathcal{R}_{f_{N-1}}^+h_{\sigma', \alpha'}^{\mathcal{T}'}$ form a Log-multiplier family near $f_{N}$.
 		We set $\alpha_{{\sigma}}^\mathcal{T}= \alpha_{{\sigma}'}^{\mathcal{T}'}$ and
 		$$p_{\sigma}^\mathcal{T}/q_\sigma^\mathcal{T} = p_{\sigma'}^{\mathcal{T}'}/ q_{\sigma'}^{\mathcal{T}'}\oplus -\mathfrak{S}(p_{\sigma'}^{\mathcal{T}'}/q_{\sigma'}^{\mathcal{T}'})/k_N\oplus -p_N/q_N,$$
 		so it follows from Proposition \ref{adding fractions} that $\mu_{\sigma}^\mathcal{T}(z)= \mu_{\sigma'}^{\mathcal{T}'}\circ (\mu_{p_{\sigma'}^{\mathcal{T}'}, q_{\sigma'}^{\mathcal{T}'}})^{-1}\circ \mu_{p_\sigma^\mathcal{T}/q_\sigma^\mathcal{T}}(-\overline{z})$. Note that as $h_{\sigma, \alpha}^{\mathcal{T}}$ depends anti-holomorphically on
 		$h_{\sigma', \alpha'}^{\mathcal{T}'}$ in this case, our definitions ensure that $h_{\sigma, \alpha}^{\mathcal{T}}$ depends holomorphically on $\alpha$ regardless of whether $\mu_{\sigma}^{\mathcal{T}}$ is holomorphic or anti-holomorphic.

 		If $f_N$ is a lower parabolic renormalization of $f_{N-1}$, then the situation is slightly more complicated. For $k_N\geq 2$ we would like to define $\mu_\sigma^\mathcal{T}$ so that 
 		$$\mu_\sigma^\mathcal{T}(\alpha)= \mu_{\sigma'}^{\mathcal{T}'}(\alpha'),$$
 		where $\alpha'$ satisfies the equation
 		\begin{equation}\label{defining mu equation}
 			\alpha' = \frac{1}{k_N+\mathfrak{c}_-^{h_{\sigma', \alpha'}^{\mathcal{T}'}, f_{N-1}}+\mu_{p_N/q_N}(\alpha)}.
 		\end{equation}
 		Fix some $R>1$.
 		As the equation 
 		$$z' = \frac{1}{k_N+\mathfrak{c}_-^{f_{N-1}}+\mu_{p_N/q_N}(z)}$$
 		has a unique solution in ${\mathbb{D}}_R$ for any $z\in \mathbb{C}$  and any large $k_N>0$, and as $\mathfrak{c}_-^{h_{\sigma', \alpha'}^{\mathcal{T}'}, f_{N-1}}$ depends holomorphically on $\alpha'$ and  tends to $\mathfrak{c}_-^{f_{N-1}}$ when $k_N\to \infty$ and $\alpha'\to 0$, it follows from  Rouch\'e's theorem that there is some $K>0$ so that  the equation (\ref{defining mu equation}) has a unique solution in $\mathbb{D}_R$ for any $z\in \mathbb{D}$ when $k_N\geq K.$ 
 		It follows from the implicit function theorem that the resulting map $\alpha\mapsto \alpha'$ is univalent on $\mathbb{D}_R$, allowing us to define $\mu_\sigma^\mathcal{T}$. We set $p_{\sigma}^\mathcal{T}/ q_{\sigma}^\mathcal{T}:= -1/k_N\oplus p_N/q_N$, so 
 		$$\mu_{p_\sigma^\mathcal{T}/q_\sigma^\mathcal{T}}(z) = \frac{-1}{k_N+\mu_{p_N/q_N}(z)}.$$
 		For $\alpha\in \mathbb{D}$ and $\mu_{\sigma}^\mathcal{T}(\alpha)= \mu_{\sigma'}^{\mathcal{T}'}(\alpha')$, 
 		 the maps $h_{\sigma, \alpha}^\mathcal{T}:= \mathcal{R}_{f_{N-1}}^-h_{\sigma', \alpha'}^{\mathcal{T}'}$ form a Log-multiplier family near $f_{N}$.
 		It follows from Proposition \ref{horn map perturbed} that the map $h_{\sigma', \alpha'}^{\mathcal{T}'}$ has a  periodic cycle in $\mathbb{C}^*$ whose multiplier  is $\Exp \circ \mu_{p_\sigma^\mathcal{T}/q_\sigma^\mathcal{T}}(\alpha).$ Proposition \ref{lifting dynamics} therefore implies that $Q_{\mu_{\sigma}^{\mathcal{T}}(\alpha)}$ has a periodic cycle with either the same multiplier or complex conjugate multiplier. Thus  $\mu_{\sigma}^\mathcal{T}$ maps $\mathbb{H}\cap \mathbb{D}_R$ into a component $\mathcal{H}_{\alpha_\sigma^\mathcal{T}}$ for some parabolic parameter $\alpha_\sigma^\mathcal{T}$ and either 
 		$$\mu_{\sigma}^\mathcal{T}(z)= \tau_{\alpha_\sigma^\mathcal{T}}\circ \mu_{p_\sigma^\mathcal{T}/q_\sigma^\mathcal{T}}\text{ or }\mu_{\sigma}^\mathcal{T}(z)= \tau_{\alpha_\sigma^\mathcal{T}}\circ \mu_{p_\sigma^\mathcal{T}/q_\sigma^\mathcal{T}}(-\overline{z})$$
 		for all $z\in \mathbb{H}\cap \mathbb{D}_R.$
 		By taking $R$ sufficiently large, it follows from the uniqueness in the definition of $\mu_\sigma^\mathcal{T}$ that  $\alpha_\sigma^\mathcal{T}$ depends only $\sigma'$. For fixed $\sigma'$, as $\mu_\sigma^\mathcal{T}(0)\to \mu_{\sigma'}^{\mathcal{T}'}(0)$ when $k_N\to \infty$ we have $\alpha_\sigma^\mathcal{T}= \mu_{\sigma'}^{\mathcal{T}'}(0).$
 		It follows from the inductive hypothesis that the combinatorics  of $\alpha_\sigma^\mathcal{T}$ are therefore $\omega_{\sigma'}^{\mathcal{T}'}\oplus (p_{\sigma'}^{\mathcal{T}'}/q_{\sigma'}^{\mathcal{T}'}).$
 		As
 		$$\mu_{\sigma}^\mathcal{T}\circ (\mu_{p_\sigma^\mathcal{T}/q_\sigma^{\mathcal{T}}})^{-1}(z) = \mu_{\sigma'}^{\mathcal{T}'}\left(\frac{-z}{1-z\cdot \mathfrak{c}_-^{h_{\sigma', \alpha'}^{\mathcal{T}'}, f_{N-1}} }\right)$$ and $ \mu_{p_\sigma^\mathcal{T}/q_\sigma^\mathcal{T}}(\mathbb{D}_R)$ converges uniformly to zero when $k_N\to \infty$ for any $R>0$, it follows from the inductive hypothesis and Proposition \ref{Mobius geometry} that
 		\begin{align*}
 			|\mu_{\sigma}^\mathcal{T}\circ (\mu_{p_\sigma^\mathcal{T}/q_\sigma^{\mathcal{T}}})^{-1}(z)-\mu_{\sigma}^\mathcal{T}\circ (\mu_{p_\sigma^\mathcal{T}/q_\sigma^{\mathcal{T}}})^{-1}(w)|\leq \frac{2^{d-1}}{\|\omega_{\sigma'}^{\mathcal{T}'}\|}\cdot \frac{\sqrt{2}}{(q_{\sigma'}^{\mathcal{T}'})^2}\cdot {\sqrt{2}|z-w|} = \frac{2^{d}|z-w|}{\|\omega_{\sigma}^\mathcal{T}\|}
 		\end{align*}
 		and
 		$$|\mu_{\sigma}^\mathcal{T}\circ (\mu_{p_\sigma^\mathcal{T}/q_\sigma^{\mathcal{T}}})^{-1}(z)-\mu_{\sigma}^\mathcal{T}\circ (\mu_{p_\sigma^\mathcal{T}/q_\sigma^{\mathcal{T}}})^{-1}(w)|\geq \frac{2^{-d+1}}{\|\omega_{\sigma'}^{\mathcal{T}'}\|}\cdot \frac{1}{\sqrt{2}(q_{\sigma'}^{\mathcal{T}'})^2}\cdot \frac{|z-w|}{\sqrt{2}}= \frac{2^{-d}|z-w|}{\|\omega_\sigma^\mathcal{T}\|}$$
 		when $k_N$ is sufficiently large, where $d$ is the depth of $\omega_\sigma^\mathcal{T}$.
 		Note that our choice of constant $2^d$ is arbitrary, we could instead inductively choose $\lambda^d$ for any $\lambda>1$ by making $K^\mathcal{T}$ larger.
 		
 		If $f_N$ is a lower parabolic renormalization of $f_{N-1}^*$, then we want to define $\mu_\sigma^\mathcal{T}(\alpha) = \mu_{\sigma'}^{\mathcal{T}'}(\alpha')$
 		for 
 		$$\alpha' = \frac{-1}{k_N+ \overline{\mathfrak{c}_-^{({h_{\sigma', \alpha'}^{\mathcal{T}'}})^*, f_{N-1}^*}}+ \mu_{p_N/q_N}(\overline{\alpha})}.$$
 		By the same reasoning above, we can define the map $\mu_\sigma^\mathcal{T}$ and Log-multiplier family of maps $h_{\sigma, \alpha}^\mathcal{T}= \mathcal{R}_{f_{N-1}}^-h_{\sigma', \alpha'}^{\mathcal{T}'}.$ Setting $p_{\sigma}^\mathcal{T} / q_{\sigma}^{\mathcal{T}}:= 1/k_N\oplus p_N/q_N$, in this case it follows that $h_{\sigma', \alpha'}^{\mathcal{T}'}$ has a periodic cycle whose multiplier is $\Exp\circ \mu_{p_{\sigma}^{\mathcal{T}}/q_{\sigma}^{\mathcal{T}}}(-\overline{\alpha})$. The rest of the proof follows similarly. 
 	\end{proof}
 
 	For any $\sigma \in {\Sigma}^\mathcal{T}$, we will say that a parameter $\zeta$ is \textit{$\sigma$-close} to $\mathcal{T}$ if $(\mu_{\sigma}^\mathcal{T})^{-1}(\zeta)\in (-1/2, 1/2].$ We will also call $\mu_\sigma^\mathcal{T}(0)$  the \textit{$\sigma$-approximation} of $\mathcal{T}$.
	 It follows from Corollary \ref{inductive compatibility} that the Log-multiplier family $\{h_{\sigma, \alpha}^\mathcal{T}\}$ has compatible pre-petals and external rays. Applying Proposition \ref{cage} yields:

	 \begin{theorem}\label{main theorem}
	 	Assume $N>0$. There is a constant $C>1$ which depends only on $f_{N-1}$ such that if
	 	$\zeta$ is a parabolic parameter which is  $\sigma$-close to $\mathcal{T}$ for some $\sigma\in \Sigma^\mathcal{T}$ with $|\sigma|$ sufficiently large,  
	 	then 
	 	$$\sup_{z\in \mathcal{L}_\zeta}\left|\mu_{p_N/q_N}\circ (\mu_\sigma^{\mathcal{T}})^{-1}(z) \right| < C.$$  Moreover if $\alpha$ is the $\sigma$-approximation of $\mathcal{T}$, then there exists some $\tilde{C}>1$ which depends only on $f_{N}$ such that
	 	 $$\inf_{z\in \mathcal{L}_\zeta}\left| (\mu_\sigma^{\mathcal{T}})^{-1}(z) \right| > \tilde{C}.$$
	 \end{theorem}
	 
	 \begin{proof}
	 	Without  loss of generality we will assume that ${f}_{N}$ is an upper satellite parabolic renormalization of $f_{N-1}$; the other cases can be handled similarly. 
	 	
	 	We set  $\sigma= \langle k_{n}\rangle_{n=1}^N$,  $\sigma' = \langle k_{n}\rangle_{n=1}^{N-1}$, and $\mathcal{T}' = \langle f_n\rangle_{n=0}^{N-1}.$
	 	We also set $\mu_\sigma^\mathcal{T}(\alpha)= \mu_{\sigma'}^{\mathcal{T}'}(\alpha')= \zeta$, so the proof of Proposition \ref{change of coordinates} implies that 
	 	$$\alpha' = \frac{1}{k_N-\mu_{p_N/q_N}(\alpha)}.$$
	 	Let $\Psi^{{\sigma}'}$ be the external parameter as in Proposition \ref{parameter rays} for the Log-multiplier family $h_{\sigma', \alpha'}^{\mathcal{T}'}$ near $f_{N-1}$	 	
	 	and let $\gamma_{{\sigma}}:{[-1, 1]}\to \mathbb{H}$ be the curve as in Proposition \ref{cage}. We define the set 
	 	$$\Gamma:= \left[\frac{1}{k_N+1}, \frac{1}{k_N-1}\right]\cup \gamma^{{\sigma}}([-1, 1]).$$
	 	It follows from Proposition \ref{cage} that $\Gamma$ is a closed curve in $\mathbb{C}$, we  denote by  $X$ the union of the bounded components of $\mathbb{C}\setminus \Gamma$.
	 	Proposition \ref{cage}  and the maximum modulus principle imply that there is a constant $C>1$ which depends only on $f_{N-1}$ such that  
	 	$$\left|k_N- \frac{1}{z}\right|< C$$
	 	for all $z\in X.$

	 	As $\alpha$ is $\sigma$-close to $\mathcal{T}$, we have $\alpha\in (-1/2, 1/2],$ so $\alpha'\in (\frac{1}{k_N+1}, \frac{1}{k_N-1})$.
	 	As  every point in $\mu_{\sigma'}^{\mathcal{T}'}(-\mathbb{H})$ is contained in $\mathcal{M}$, thus there must be some component $\mathcal{L}$ of $\mathcal{M}\setminus \{\zeta\}$ which avoids $\mu_{\sigma'}^{\mathcal{T}'}(-\overline{\mathbb{H}})$.
	 	Proposition \ref{relating parameter coordinates} implies that
	 	$\mathcal{M}$ avoids  $\mu_{\sigma'}^{\mathcal{T}'}(\Gamma)$, 
	 	hence $\mathcal{L}\subset \mu_{\sigma'}^{\mathcal{T}'}(X).$ Thus $\mathcal{L}$ is bounded, so $\mathcal{L}= \mathcal{L}_{\zeta}$. It therefore follows from the above that $$\sup_{z\in \mathcal{L}_\zeta}\left| \mu_{p_N/q_N}\circ (\mu_\sigma^\mathcal{T})^{-1}(z)\right|=\sup_{z\in \mathcal{L}_\zeta}\left|k_N-\frac{1}{(\mu_{\sigma'}^{\mathcal{T}'})^{-1}(z)}\right|< C$$ when $|\sigma|$ is sufficiently large.

	 	We now consider the satellite tower $\tilde{\mathcal{T}}:= \mathcal{T}\oplus \langle \mathcal{R}_0^-f_N \rangle$ and set $\tilde{\sigma}:= \sigma\oplus \langle k\rangle$ for some $k\gg 0$. Choosing $\tilde{C}>1$ so that $C^{-1}\ll \frac{1}{k+\mathfrak{C}_-^{f_N}}$, it follows from Proposition \ref{change of coordinates} that $\mu_{\tilde{\sigma}}^{\tilde{\mathcal{T}}}(0)\in \mathcal{L}_\alpha$ and 
	 	$$\left(\mu_\sigma^\mathcal{T}\right)^{-1}\circ \mu_{\tilde{\sigma}}^{\tilde{\mathcal{T}}}(0) \approx \frac{1}{k+\mathfrak{c}_-^{f_N}}> \tilde{C}^{-1}.$$
	 \end{proof}
 
 	We now observe that Theorem \ref{main} from the introduction follows from Theorem \ref{main theorem}.
 
 	\begin{corollary}\label{main corollary}
 		For any positive integers $N$ and $M$ there exists $\epsilon>0$ and $C>1$ such that if $\alpha$ has combinatorics $\omega= \langle p_n/q_n\rangle_{n=1}^{N}$ with $p_n/q_n\in \bigcup_{m=1}^{M}\mathbb{Q}_{m}$ for all $n\geq 1$ and $|p_n/q_n|< \epsilon$ for all $n\geq 2$, 
 		then 
 		$$\frac{C^{-1}}{\|\omega\|}< \emph{Diam}\, \mathcal{L}_{\alpha} < \frac{C}{\|\omega\|}.$$
 	\end{corollary}
 	\begin{proof}
 		If the proposition does not hold, then 
 		there is a sequence of parabolic parameters $\langle\alpha_j\rangle_{j=1}^\infty$ such that 
 		${\alpha_j}$ has combinatorics $\omega_j:=\langle p_{n, j}/q_{n, j}\rangle_{n=1}^{N}$, where $p_{n, j}/q_{n, j}\in \bigcup_{m=1}^{M}\mathbb{Q}_{m}$, and either 
 		$${\|\omega_j\|\cdot \text{Diam}\, \mathcal{L}_{\alpha_j}}\to \infty\text{ or }{\|\omega_j\|\cdot \text{Diam}\, \mathcal{L}_{\alpha_j}}\to 0$$ when $j\to \infty$. Moreover, $p_{n, j}/q_{n, j}\to 0$ when $j\to \infty$ for all $2\leq n \leq N$.

 		For any $1\leq n \leq N$, up to a subsequence we can write express $p_{n, j}/q_{n, j}$ as
 		$$\frac{p_{n, j}}{q_{n, j}}= \frac{\tilde{p}_{n, 1}}{\tilde{q}_{n, 1}}\oplus \frac{\tilde{\varepsilon}_{n, 1}}{k_{n, 1, j}}\oplus \cdots  \oplus \frac{\tilde{p}_{n, M_n}}{\tilde{q}_{n, M_n}}\oplus \frac{\tilde{\varepsilon}_{n, M_n}}{k_{n, M_n, j}}\oplus \frac{\tilde{p}_{n, M_{n}+1}}{\tilde{q}_{n, M_n+1}},$$
 		where $\tilde{p}_{n, m}/\tilde{q}_{n, m}$ are rational numbers in $(-1/2, 1/2]$, $\varepsilon_{n, m}= \pm 1$, and $k_{n, m, j}$ are positive integers which tend to $\infty$ when $j\to \infty$ for all $m.$ Note that $p_{n, j}/q_{n, j}\to 0$ when $j\to \infty$ implies that $\tilde{p}_{n, 1}/\tilde{q}_{n, 1}= 0/1$ for all $n\geq 2$. We observe that part (\ref{combinatorics of near tower}) of Proposition \ref{change of coordinates} gives us an explicit algorithm for defining a satellite tower $\mathcal{T}= \langle f_s\rangle_{s=0}^S$ so that for all large $j$ there exists $\sigma_j\in \Sigma^\mathcal{T}$ such that  $\alpha_j=\mu_{\sigma_j}^{\mathcal{T}}(0)$ and $|\sigma_j|\to \infty$ when $j \to \infty$. Theorem \ref{main theorem} therefore implies that there is a constant $C>1$ such that 
 		$$C^{-1}< \inf_{z\in\mathcal{L}_{\alpha_j}}\left|(\mu_{\sigma}^\mathcal{T})^{-1}(z)\right|\leq \sup_{z\in\mathcal{L}_{\alpha_j}}\left|(\mu_{\sigma}^\mathcal{T})^{-1}(z)\right|< {C}$$
 		for infinitely many $j$. Part (\ref{distortion tower}) of Proposition \ref{change of coordinates} therefore implies that, after increasing $C$, we have
 		$$ \frac{C^{-1}}{\|\omega_j\|}< \text{Diam}\, \mathcal{L}_{\alpha_j} < \frac{C}{\|\omega_j\|}$$
 		for infinitely many $j$, but this is a contradiction.
 	\end{proof}
 
 	Let us note that while Corollary \ref{main corollary} is the main result of this paper, all of the real work is done by Theorem \ref{main theorem}. For some applications other consequences of Theorem \ref{main theorem} may prove more useful, the following is one such example which is used in \cite{taming}:
 	\begin{corollary}
 		For any rational $p_0/q_0$ there exists a constant $C>0$ such that if $\alpha$ has combinatorics $\langle p_0/q_0, p_1/q_1\rangle$ for any non-zero rational $p_1/q_1\in (-1/2, 1/2]$, then
 		$$\sup_{z\in \mathcal{L}_\alpha}\left| \frac{1}{q_0 z - p_0} - \frac{1}{q_0 \alpha - p_0} \right|< C.$$
 	\end{corollary}
 
 	\begin{proof}
 		First we observe that the PLY inequality implies that  for any $\epsilon>0$ there exists a neighborhood $\Delta$ of $p_0/q_0$ such that  if $|p_1/q_1|>\epsilon$ then 
 		$\mathcal{L}_{\alpha}$ avoids $\Delta$.  In this case it follows that there exists $C_\epsilon>0$  such that $$\sup_{z\in \mathcal{L}_\alpha} \left|\frac{1}{q_0 z - p_0} \right|< {C_\epsilon}.$$ 		
 		Let us now assume that $0<-p_1/q_1< \epsilon$, so $p_1/q_1 = \frac{-1}{n+\theta}$ for some integer $n\gg0$ and rational $\theta\in (-1/2, 1/2].$
 		Thus $\alpha$ is $\langle n\rangle$-close to the satellite tower $\mathcal{T}=\langle Q_{p_0/q_0}, \mathcal{R}_0^-Q_{p_0/q_0}\rangle$ by Proposition \ref{change of coordinates}, so Theorem implies that there is a constant $C>0$ which depends only on $Q_{p_0/q_0}$ such that if $\epsilon$ is sufficiently small then 
 		$$\sup_{z\in \mathcal{L}_\alpha}\left|(\mu_{\langle n\rangle}^\mathcal{T})^{-1}(z)\right|< {C}.$$
 		The definition of $\mu_{\langle n \rangle}^\mathcal{T}$ in the proof of Proposition \ref{change of coordinates} implies that
 		$$\sup_{z\in \mathcal{L}_\alpha}\left| n+\mathfrak{c}_-^{Q_{p_0/q_0}}- \frac{1}{\mu_{p_0/q_0}^{-1}(z)}\right| \leq  \sup_{z\in \mathcal{L}_\alpha}\left|(\mu_{\langle n\rangle}^\mathcal{T})^{-1}(z)\right| +1$$
 		when $n$ is sufficiently large, so 
 		$$\sup_{z\in \mathcal{L}_\alpha}\left|\frac{1}{\mu_{p_0/q_0}^{-1}(\alpha)}- \frac{1}{\mu_{p_0/q_0}^{-1}(z)}\right|< 2\left(C+1\right)$$
 		It follows from the definition of $\mu_{p_0/q_0}$ that 
 		$$\frac{1}{q_0w-p_0} = \frac{q_0}{\mu_{p_0/q_0}^{-1}(w)}+\mathfrak{S}(p_0/q_0)q_0'$$
 		for any $w\in \mathbb{C}$, where $p_0'/q_0'$ is the parent of $p_0/q_0$. 
 		Thus $$\sup_{z\in \mathcal{L}_\alpha}\left| \frac{1}{q_0 z - p_0} - \frac{1}{q_0 \alpha - p_0} \right|< 2\left(C+1\right)q_0.$$
 		The argument for the case where $0 < p_1/q_1< \epsilon$ is the same, using instead the tower
 		$\langle Q_{p_0/q_0}, \mathcal{R}_0^-Q_{p_0/q_0}^*\rangle$.
 	\end{proof}

 	\subsection{Generalizations of Theorem \ref{main theorem}}\label{generalizations}
 	
 	We conclude with some remarks on possible extensions of Theorem \ref{main theorem} to more general parabolic towers. 
 	
 	First let us consider the case where $\mathcal{T}= \langle f_n\rangle_{n=0}^\infty$ is an infinite height satellite tower. 
 	For all $N\geq 0$ we denote $\mathcal{T}_N = \langle f_n\rangle_{n=0}^N$, and for all  $\sigma\in \Sigma^{\mathcal{T}_N}$ we set $\mu_{\sigma}^\mathcal{T}:= \mu_{\sigma}^{\mathcal{T}_N}$ and $h_{\sigma, \alpha}^\mathcal{T}:= h_{\sigma, \alpha}^{\mathcal{T}_N}$. 
 	Theorem \ref{main theorem} provides us with constants $C_N$ and $\tilde{C}_N$ for $\mathcal{T}_N$ which depend only on $f_{N-1}$ and $f_N$ respectively. 
 	To generalize Theorem \ref{main theorem} to $\mathcal{T}$, we need to control these constants for when $N$ tends to $\infty$. 
 	Let us say that the tower $\mathcal{T}$ is \textit{convergent} when there exists some $f_\infty\in \mathcal{F}^\flower$ such that $f_N\to f_\infty$ when $N\to \infty$.  	
 	As $\mathcal{R}_{f_N}h = \mathcal{R}_{f_\infty}h$ when $N$ is sufficiently large, the convergence $\mathcal{T}$ implies that the maps $h_{\sigma, \alpha}^{\mathcal{T}}$  form a Log-multiplier family near $f_\infty$. 
 	If $f_\infty$ has Jordan basin and this Log-multiplier family has compatible external rays and pre-petals, then we can apply Proposition \ref{cage} to this family and get uniform control over the constants $C_N$ and $\tilde{C}_N$.
 	Some examples of convergent towers are given in \cite{shishikura} and \cite{nondegenerate}, such as when $f_N = \mathcal{R}_0^+f_{N-1}$ for all $N$ or $f_N = \mathcal{R}_0^-f_{N-1}$ for all $N$. By Proposition \ref{change of coordinates}, uniform control for these two examples would have the following consequences respectively:
 	\begin{conjecture}
 		There is an $M\geq 2$ and $C>1$ such that if $p/q$ has modified continued fraction $\langle (a_n, \varepsilon_n)\rangle_{n=1}^N$ with $a_n\geq M$ for all $n$, then 
 		$$\frac{C^{-1}}{q^2}<\emph{Diam}\,\mathcal{L}_{p/q}< \frac{C}{q^2}.$$
 	\end{conjecture}
 	\begin{conjecture}
 		There is an $M\geq 2$ and $C>1$ such that if $\zeta$ has combinatorics $\langle 1/q_n\rangle_{n=1}^N$ with $q_n\geq M$ for all $n$, then 
 		$$\frac{2^{-N}C^{-1}}{\prod_{n=1}^Nq_n^2}<\emph{Diam}\,\mathcal{L}_{\zeta}< \frac{2^NC}{\prod_{n=1}^Nq_n^2}.$$
 	\end{conjecture}
 	\noindent We may also expect uniform control of the constants $C_N$ and $\tilde{C}_N$  when  $f_N = \mathcal{R}_0^\pm f_{N-1}$ for all $N$, this would imply Conjecture \ref{main conjecture} from the introduction.

 	Now let us consider when $\mathcal{T}= \langle f_n\rangle_{n=0}^N$ is a quadratic parabolic tower but is not a satellite tower. 
 	We will say that $\mathcal{T}$ is \textit{persistently parabolic} if the parabolic periodic cycle for some $f_n$ does not bifurcate under perturbation. More precisely, if there exists some $0< n\leq N$ such that $f_n = \mathcal{R}_{\delta_{n}}^\pm f_{n-1}$ has a $p/q$-parabolic $k$-periodic cycle $\langle z_j\rangle_{j=0}^{q-1}$ and $\mathcal{R}_{\delta}^\pm f_{n-1}$ has a $p/q$-parabolic $k$-periodic cycle close to $\langle z_j\rangle_{j=0}^{q-1}$ for all $\delta$ close to $\delta_n$.
 	As an example, there are no persistently parabolic satellite towers; we conjecture that in fact $\mathcal{T}$ is  never persistently parabolic. 
 	If $\mathcal{T}$ has finite height and is not persistently parabolic, then we can recreate Proposition 
 	\ref{change of coordinates} and similarly construct $\Sigma^\mathcal{T}, \mu_\sigma^\mathcal{T}$, and $ h_{\sigma, \alpha}^\mathcal{T}.$ However some minor modifications may be necessary: we may need to conjugate by a M\"obius transformation and consider an iterate so that  $h_{\sigma, \alpha}^\mathcal{T}$ has a fixed point at zero, and if the multiplier of the parabolic cycle of $f_N$ is one then $\mu_{\sigma}^\mathcal{T}$ may not be injective. After these modifications, if $f_N$ has Jordan basin and $h_{\sigma, \alpha}^\mathcal{T}$ has compatible external rays and pre-petals which are induced by the external rays of $Q_{\mu_\sigma^\mathcal{T}(\alpha)}$ as in the satellite tower case, then an analogous version of Theorem \ref{main theorem} holds by identical argument. We conjecture that these conditions are always satisfied, so we expect that Theorem \ref{main theorem} can be generalized to all quadratic parabolic towers. 
 	
 	Finally, we observe that in our analysis we make no use of the fact that all the critical points of maps in $\mathcal{F}$ have local degree two; we only need a unique critical value. By modifying the definition of $\mathcal{F}$, namely allowing the critical points to have arbitrary degree, we can produce similar results for more general classes of maps. In particular,  the analogues of Theorem \ref{main theorem} and Corollary \ref{main corollary} hold for the family of unicritical  polynomials  $$z\mapsto \left(\frac{\Exp(\alpha)\cdot z+d}{d}\right)^d-1$$
 	for any degree $d\geq2$.

	\bibliographystyle{alpha}
	\bibliography{bibliography}
	\ \\
	\textsc{Department of Mathematics, Harvard University, Cambridge, Massachusetts}\par\nopagebreak
	\medskip
	\noindent\textit{E-mail address}: \href{mailto:kapiamba@math.harvard.edu}{\texttt{kapiamba@math.harvard.edu}}

\end{document}